\documentclass[12pt]{amsart}
\usepackage{amssymb}
\usepackage{amsmath}
\usepackage{versions}
\usepackage[usenames,dvipsnames]{color}

\usepackage{systeme}

\def\d{\mathbb{D}}
\def\c{\mathbb{C}}

\def\Z{\mathbb{Z}}

\def\be{\begin{equation}}
\def\ee{\end{equation}}

\def\tetrainn{{\mathbb{\overline{E}}}\text{-inner}}
\def\s0{s_0}
\def\p0{p_0}
\let\phi\varphi
\let\epsilon\varepsilon

\newtheorem{theorem}{Theorem}[section]
\newtheorem{corollary}[theorem]{Corollary}
\newtheorem{lemma}[theorem]{Lemma}
\newtheorem{proposition}[theorem]{Proposition}

\numberwithin{equation}{section}

\newtheorem{example}[theorem]{Example}

\newtheorem{definition}[theorem]{Definition}

\newtheorem{remark}[theorem]{Remark}

\newtheorem{fact*}[theorem]{Fact}

\newcommand\e{\mathrm{e}}
\newcommand{\eps}{\varepsilon}

\newcommand{\T}{\mathbb{T}}

\newcommand{\D}{\mathbb{D}}
\newcommand{\C}{\mathbb{C}}

\newcommand{\E}{\mathbb{E}}

\newcommand\la{\lambda}

\newcommand\id{\rm id}

\newcommand\Aut{\mathrm{Aut}~}

\newcommand\beq{\begin{equation}}

\newcommand\eeq{\end{equation}}

\newcommand\bbm{\begin{bmatrix}}
\newcommand\ebm{\end{bmatrix}}
\newcommand\bpm{\begin{pmatrix}}
\newcommand\epm{\end{pmatrix}}

\let\phi\varphi

\numberwithin{equation}{section}
\begin{document}
\title[Rational $\overline{\E}$-inner functions]{Rational  tetra-inner  functions and the special variety of the tetrablock}
\author{Omar M. O. Alsalhi and Zinaida A. Lykova}
\date{26th July 2021}

\begin{abstract} 
The set
\[
\overline{\E}= \{ x \in \C^3: \quad 1-x_1 z - x_2 w + x_3 zw \neq 0 \mbox{  whenever }
|z| < 1, |w| < 1 \}
\]
is called the tetrablock and has intriguing complex-geometric properties. It  is  polynomially convex, nonconvex and starlike about $0$.
It has a group of automorphisms parametrised by $\Aut \D \times \Aut \D \times \Z_2$ and its distinguished boundary $b\overline{\E}$ is homeomorphic to the  solid torus $\overline{\D} \times \T$. It has a special subvariety
\[\mathcal{R}_{\mathbb{\overline{E}}} = \big\{ (x_{1}, x_{2}, x_{3}) \in \overline{\E}   : x_{1}x_{2}=x_{3} \big\},
\] 
called {\em the 
royal variety} of $\overline{\E}$. $\mathcal{R}_{\mathbb{\overline{E}}}$  is a complex geodesic of ${\E}$ and it is invariant under all automorphisms of ${\E}$. We make use of  these geometric properties of $\overline{\E}$ to develop an explicit structure theory for the rational maps from the unit disc $\D$  to $\overline{\E}$ that map the
unit circle  $\T$ to the distinguished boundary $b\overline{\E}$ of $\overline{\E}$. Such maps are called rational $\mathbb{ \overline{ E}}$-inner functions.
 We call the points $\la\in \overline{\D}$ such that $x(\la)\in \mathcal{R}_{\mathbb{\overline{E}}}$ the {\em royal nodes} of $x$.
 We describe the construction of rational ${\mathbb{\overline{E}}}$-inner functions  of prescribed degree from the zeros  of $x_{1}$ and $x_{2}$ and  the royal nodes of $x$. The proof of this theorem is constructive: it gives an algorithm for the construction of a $3$-parameter family of such functions $x$ subject to the computation of  Fej\'er-Riesz factorizations of certain non-negative functions on the circle.
We show that, for each nonconstant  rational $\tetrainn$ function $x$, either $x(\overline{\mathbb{D}}) \subseteq \mathcal{R}_{\mathbb{\overline{E}}} \cap \overline{\E}$ or $x(\overline{\mathbb{D}})$ meets $\mathcal{R}_{\mathbb{\overline{E}}}$ exactly $\deg(x)$ times.
 We study convex subsets of the set $\mathcal{J}$ of all rational $\mathbb{ \overline{ E}}$-inner functions and extreme points of $\mathcal{J}$.
We show that whether a  rational $\mathbb{ \overline{ E}}$-inner function $x$  is an extreme point of $\mathcal{J}$ depends on how many royal nodes of $x$ lie on $\mathbb{ T}$.
\end{abstract}



\thanks{The first author  was supported by the Government of Saudi Arabia. The second author was partially supported by the UK Engineering and Physical Sciences Research Council grants  EP/N03242X/1. }


\maketitle
\tableofcontents
\section{Introduction}\label{intro}

Some unsolved problems in $H^\infty$ control theory  require a deep understanding
of ``inner mappings" from $\D$ to certain domains in $\C^d$ with $d>1$.  For example, a special case of the problem of robust stabilization under structured uncertainty, or the  $\mu$-synthesis problem \cite{Do,dullerud,AY99,ALY2015},  leads naturally to a class of $\mu$-synthesis domains.
A typical member of this class of domains  is
\[
{\E}= \{ x \in \C^3: \quad 1-x_1 z - x_2 w + x_3 zw \neq 0 \mbox{  whenever }
|z| \le 1, |w| \le 1 \}
\]
the {\em open tetrablock, } as was observed by Abouhajar, White and  Young in \cite{AWY}. The complex geometry of $\E$ was further developed in \cite{EKZ2013,KZ2015,KZ2016,You08} and associated operator theory in \cite{TB,BS}. The solvability of the $\mu$-synthesis problem connected to $\E$  can be expressed in terms of the existence of rational inner functions from the open unit disc $\D$ in the complex plane $\C$ to the closure of $\E$ subject to interpolation conditions \cite{BLY16}.

Recall that a classical {\em rational inner function} is a rational map $f$ from the unit disc $\D$  to its closure $\overline{\D}$ with the property that $f$ maps the unit circle $\T$ into itself. See \cite{CGP} for a survey of results, linking inner functions and operator theory.
We denote the closure of $\E$ by $\overline{\E}$ and we define a {\em rational $\overline{\E}$-inner or {\rm (}tetra-inner{\rm )}  function} to be a rational analytic map $x:\D\to \overline{\E}$ such that  $x$ maps $\T$ into the distinguished boundary $b\overline{\E}$ of $\overline{\E}$.   Here, $b\overline{\E}$  is the smallest closed subset of $\overline{\E}$ on which every continuous function on $b\overline{\E}$ that is analytic in $\E$ attains its maximum modulus. Of course, rational $\overline{\E}$-inner functions  have  many similarities with rational inner functions. On the other hand, the complex geometry of  $\overline{\E}$ is  richer  than that  of $\overline{\D}$,  and so 
rational $\overline{\E}$-inner functions  have some striking differences from  rational inner functions.

Here are some points of difference between well-studied domains in $\C^3$,
such as the tridisc $\D^3$ and the Euclidean ball $\mathbb{B}_3$ 
 on the one hand and $\E$ on the other.
Firstly, whereas $\D^3$ and  $\mathbb{B}_3$ are homogeneous (so that the holomorphic automorphisms of these domains act transitively),
$\E$ is inhomogeneous \cite{You08}.  

Secondly, the distinguished boundary of $\E$ differs markedly in its topological properties from those of $\D^3$ and $\mathbb{B}_3$.
The distinguished boundaries of $\D^3$ and $\mathbb{B}_3$ are the $3$-dimensional torus and the $5$-sphere respectively. They are smooth manifolds without boundary. The distinguished boundary $b\overline{\E}$ of $\E$
is homeomorphic to the solid torus $\overline{\mathbb{D}} \times \T$,  which 
has a boundary. 
 For a rational $\overline{\E}$-inner function $x$ the curve $x(\e^{it}), 0\leq t<2\pi$, lies in 
$b\overline{\E}$ and may or may not touch the boundary of $b\overline{\E}$, and so the algebraic and geometric properties of $x$  are dependent on that.

We call the set
\be\label{eq1.30}
\mathcal{R}_{\mathbb{\overline{E}}} =
 \big\{ (x_{1}, x_{2}, x_{3}) \in \C^3  : x_{1}x_{2}=x_{3} \big\}
\ee
the {\em royal variety} of the tetrablock. The  complex geodesic
$\mathcal{R}_{\mathbb{\overline{E}}}\cap \E$ is invariant under the group of biholomorphic automorphisms of $\E$ \cite{You08}. This paper shows that 
the variety $\mathcal{R}_{\mathbb{\overline{E}}}$ plays a central role in the function theory of $\E$. The intersection
 $\mathcal{R}_{\mathbb{\overline{E}}}\cap b\overline{\E}$ is exactly the boundary of  $b\overline{\E}$, that is, $\{(x_1, x_2, x_1 x_2) \in \C^3: |x_1| =|x_2| =1 \}$  \cite[Theorem 7.1]{AWY}, which is homeomorphic to the $2$-torus $\T \times \T$.

These  geometric properties of $\E$ lead to very interesting facts  in the theory of rational $\overline{\E}$-inner functions  that do not have analogues in the theory of classical inner functions.

One of the main theorems of this paper is the following.

\begin{theorem}\label{prop1.20} 
Let  $x$ be a  nonconstant  rational $\tetrainn$ function. Then either $x(\overline{\mathbb{D}}) \subseteq \mathcal{R}_{\mathbb{\overline{E}}} \cap \overline{\E}$ or $x(\overline{\mathbb{D}})$ meets $\mathcal{R}_{\mathbb{\overline{E}}}$ exactly $\deg(x)$ times.
\end{theorem}
It is  Theorem \ref{times}.
Here $\deg(x)$ is the degree of $x$. In Subsection \ref{degree_E-inner} we 
define $\deg(x)$ in a natural way by means of fundamental groups. 
In Proposition \ref{deg} we show that, for any rational $\overline{\E}$-inner function $x= (x_1,x_2,x_3)$, $\deg (x)$ is equal to the degree $\deg (x_3)$ of the finite Blaschke product $x_3$.
The precise way of counting the number of times that $x(\overline{\D})$ meets $\mathcal{R}_{\mathbb{\overline{E}}}$ is also described in Section \ref{construct}.  We call the points $\la\in \overline{\D}$ such that $x(\la)\in \mathcal{R}_{\mathbb{\overline{E}}}$ the {\em royal nodes} of $x$ and, for such $\la$, we call $x(\la)$ a {\em royal point} of $x$.

Another main result is the construction of rational ${\mathbb{\overline{E}}}$-inner functions of prescribed degree from the zeros of $x_{1}$ and $x_{2}$ and  the royal nodes of $x$. One can consider this result as an analogue of the expression for a finite Blaschke product in terms of its zeros. This result is proved in Theorem \ref{cons}.

\begin{theorem}\label{thm1.15}  Let $n$ be a positive integer. Suppose 
that $\alpha_{1}^{1},...,\alpha_{k_{1}}^{1} \in \overline{\mathbb{D}}$ and $\alpha_{1}^{2},...,\alpha_{k_{2}}^{2} \in\overline{\mathbb{D}}$, where $k_{1}+k_{2}=n$. Suppose that $\sigma_{1},..., \sigma_{n} \in \overline{\mathbb{D}}$ are distinct from the points of the set $\{\alpha_{j}^{i}, j=1,...,k_{i} , i=1,2\} \cap \mathbb{T}$. Then there exists a rational~ $\;\overline{\mathbb{E}}$-inner function $x=(x_{1},x_{2},x_{3}) : \mathbb{D} \rightarrow \overline{\mathbb{E}}$ such that
	\begin{enumerate}
		\item the zeros of $x_{1}$ in $\overline{ \mathbb{D}}$, repeated according to multiplicity, are $\alpha_{1}^{1},...,\alpha_{k_{1}}^{1}$;
		\item the zeros of $x_{2}$ in $\overline{ \mathbb{D}}$, repeated according to multiplicity, are $\alpha_{1}^{2},...,\alpha_{k_{2}}^{2}$;
		\item the royal nodes of $x$ are $\sigma_{1},..., \sigma_{n} \in \overline{\mathbb{D}}$, with repetition according to the multiplicity of the nodes.
	\end{enumerate}
This function $x$ can be constructed as follows. Let $t_{+}>0$ and let $t \in \mathbb{C} \backslash \{0\}$. Let $R$ be defined by
	\begin{equation*}
	R(\lambda)= t_{+} \prod_{j=1}^{n}(\lambda-\sigma_{j})(1-\overline{\sigma_{j}}\lambda).
	\end{equation*} 
	Let $E_{1}$ be defined by
	\begin{equation*}
	E_{1}(\lambda)=t \prod_{j=1}^{k_{1}}(\lambda-\alpha_{j}^{1})\prod_{j=1}^{k_{2}}(1-\overline{\alpha}_{j}^{2}\lambda) .
	\end{equation*}

	Then the following statements hold:
	\begin{enumerate}
		\item[(i)] There exists an outer polynomial $D$ of degree at most $n$ such that
		\begin{equation*}
		\lambda^{-n}R(\lambda)+|E_{1}(\lambda)|^{2}=|D(\lambda)|^{2}
		\end{equation*}
		for all $\lambda \in \mathbb{T}$.
		
		\item[(ii)] The function $x$ defined by 
		\begin{equation*}
		x=\bigg( \frac{E_{1}}{D}, \frac{E_{1}^{\sim n}}{D},\frac{D ^{\sim n}}{D}\bigg) \; \text{where, for any polynomial $D$,} \;
D^{\sim n}(\lambda) = \lambda^n \overline{ D\big(\frac{1}{\bar{\lambda}}\big) },
		\end{equation*}
		is a rational $\overline{\mathbb{E}}$-inner function such that the degree of $x$ is equal to $n$ and conditions {\em(1), (2) and (3)} hold. The royal polynomial $R_{x}$ of $x$ is equal to $R$.
	\end{enumerate}
\end{theorem}
Here the royal polynomial of $x$ is defined as  $R_{x}= \big[D^{\sim n}D-E_{1}E_{1}^{\sim n}\big]$.
The proof of this theorem is constructive: it gives an algorithm for the construction of a $3$-parameter family of such functions $x$.

In Section \ref{convexity} we study convex subsets of the set $\mathcal{J}$ of all rational $\mathbb{ \overline{ E}}$-inner functions and extremality.
We  show that the set $\mathcal{J}$ is not convex. On the other hand, 	the subset of $\mathcal{J}$ with a fixed inner function $x_{3}$ {\em is} convex (Theorem \ref{conv13}). 
Recall that the distinguished boundary of the tridisc $\mathbb{ D}^{3}$ contain no line segments. Thus every inner function in the set of analytic functions $ {\rm Hol}(\mathbb{ D}, \mathbb{ D}^{3})$ from $\D$ to $\mathbb{ D}^{3}$ is an extreme point of $ {\rm Hol}(\mathbb{ D}, \mathbb{ D}^{3})$. However, this property is in clear contrast with the situation in the tetrablock.
We show that whether a rational inner function $x$  is an extreme point of $\mathcal{J}$ depends on how many royal nodes of $x$ lie on $\mathbb{ T}$. 
\begin{theorem} \label{parex-intr}
	Let $x$ be a rational $\mathbb{ \overline{ E}}$-inner function and let $x$ have $n$ royal nodes where $k$ of them are in $\mathbb{T}$.  If $2k \leq n$, then  $x$  is not  an extreme point of $\mathcal{J}$.  
\end{theorem}
The way of counting the number of royal nodes was introduced in Section \ref{construct}.
In Proposition \ref{exetremE} we provide a class of extreme functions of the set $\mathcal{J}$.

In \cite{AlsLyk} there is a construction of the general rational $\overline{\mathcal{E}}$-inner function $x= (x_1,x_2, x_3)$ of degree $n$, in terms of different data, namely, the royal nodes of $x$ and royal values of $x$.  The algorithm  \cite{AlsLyk} for the construction of $x$ exploits a known construction of the finite  Blaschke products of given degree which satisfy some interpolation conditions with the aid of  a Pick matrix formed from the interpolation data.

The authors are grateful to Nicholas Young for some helpful suggestions.

\section{The tetrablock $\mathbb{E}$} \label{tetrablock}

\begin{definition} \cite{AWY}
The {\em open tetrablock}  is the domain defined by 
\begin{equation*}
 \mathbb{E} = \{x\in \mathbb{C}^3  : 1-x_{1}z -x_{2}w + x_{3} zw \neq 0    \enspace  \text{\rm for} \enspace  |z| \leq 1, |w| \leq 1\}.
\end{equation*}
\end{definition}
Despite the fact that $\mathbb{E }$ is not convex, its intersection with $\mathbb{R}^{3}$ {\em is}. It is proved in \cite{AWY} that $\mathbb{E} \ \cap \ \mathbb{R}^{3}$ is the open tetrahedron with the vertices $(1,1,1), (1,-1,-1),  (-1,1,-1)$  and $(-1,-1,1)$. The following function plays an important role in the study of the tetrablock.

\begin{definition} \label{fracfn}
For $ x=(x_{1},x_{2},x_{3}) \in \mathbb{C}^3$ and $z \in \mathbb{C}$ we define
\begin{eqnarray*}
\Psi(z,x) &=&\frac{x_{3}z-x_{1}}{x_{2}z-1} , \qquad \text{whenever} \ x_{2}z \neq 1. 
\end{eqnarray*}
\end{definition}
\begin{remark} {\em 
In the case that $x_{3}=x_{1}x_{2}$, $z \in \mathbb{C}$,
\begin{equation*}
\Psi(z,x) =\frac{x_{1}x_{2}z-x_{1}}{x_{2}z-1}= \frac{x_{1}(x_{2}z-1)}{x_{2}z-1}=x_{1}.
\end{equation*} 
}
\end{remark}
\begin{theorem} {\em \cite[Theorem 2.2]{AWY}} \label{tetbnd}
Let $x \in \mathbb{C}^3$. The following are equivalent
\begin{enumerate}
  \item $x \in \mathbb{E}$;
  \item $ \|\Psi(.,x) \|_{H^{\infty} }< 1$ and if $x_{1}x_{2} = x_{3}$, then, in addition, $|x_{2}|<1;$
       \item $|x_{1}-\overline{x}_{2}x_{3}|+|x_{2}-\overline{x}_{1}x_{3}| < 1 - |x_{3}|^{2}$;
        \item there exists a $2\times2$ matrix $A=[a_{ij}]$ such that $\|A\| < 1$ and $x=(a_{11},a_{22}, \det(A))$; 
         \item $|x_{3}| < 1$ and there exist $\beta_{1}, \beta_{2} \in \mathbb{C}$ such that $|\beta_{1}|+|\beta_{2}|<1$ and $$x_{1}=\beta_{1}+\overline{\beta}_{2}x_{3},\quad x_{2}=\beta_{2} + \overline{\beta}_{1}x_{3}.$$
\end{enumerate}
\end{theorem}

\begin{theorem} {\em \cite[Theorem 2.4]{AWY}} \label{tetbnd2}
Let $x \in \mathbb{C}^3$. The following are equivalent
\begin{enumerate}
  \item $x \in \overline{\mathbb{E}}$;
  \item $ \|\Psi(.,x) \|_{H^{\infty} } \leq 1$ and if $x_{1}x_{2} = x_{3}$, then, in addition, $|x_{2}| \leq 1;$ 
       \item $|x_{1}-\overline{x}_{2}x_{3}|+|x_{2}-\overline{x}_{1}x_{3}| \leq 1 - |x_{3}|^{2}$  and if  $|x_3| =1$ then, in addition,  $|x_{1}| \leq 1$;
        \item there exists a $2\times2$ matrix $A=[a_{ij}]$ such that $\|A\| \leq 1$ and $x=(a_{11},a_{22}, \det(A))$; 
         \item $|x_{3}| \leq 1$ and there exist $\beta_{1}, \beta_{2} \in \mathbb{C}$ such that $|\beta_{1}|+|\beta_{2}|\leq1$ and $$x_{1}=\beta_{1}+\overline{\beta}_{2}x_{3},\quad x_{2}=\beta_{2} + \overline{\beta}_{1}x_{3}.$$ 
\end{enumerate}
\end{theorem}

    \subsection{The tetrablock and the $\mu_{\rm Diag}$-synthesis problem}\label{mu-diag}

The tetrablock is associated with the  $\mu_{\rm Diag}$-synthesis problem from $\mathbb{D}$ to $\mathbb{C}^{2 \times 2}$.  The structured singular value in this case is defined by
\begin{equation} \label{mu}
\mu_{\rm Diag}(A)= \frac{1}{\inf\{\|X\|: X\in {\rm Diag},\: \det(I-AX)=0\}},
\end{equation}
where $${\rm Diag}:=\Bigg\{ \begin{bmatrix} z &0 \\ 0&  w \end{bmatrix} \ : z,w \in \mathbb{C}\Bigg\}.$$  
We set $\mu_{\rm Diag}(A) =0$ if $(I-AX)$ is non-singular for all $X \in {\rm Diag}$. 

\begin{definition} \label{pi12}
 We define the map $\pi: \mathbb{C}^{2 \times 2} \rightarrow \mathbb{C}^{3}$  for a matrix $A=\begin{bmatrix} a_{11} &a_{12} \\ a_{21}&  a_{22} \end{bmatrix}$ in $\mathbb{C}^{2 \times 2}$ to be
\begin{equation*}
\pi(A)= (a_{11},a_{22},\det(A)),
\end{equation*}
and $\Sigma$ to be
\begin{equation*} \index{$\Sigma$}
\Sigma:= \{A \in \mathbb{C}^{2\times 2} : \mu_{\rm Diag}(A) < 1\},
\end{equation*}
where $\mu_{\rm Diag}(A)$ is defined by equation \eqref{mu}.
\end{definition}

\begin{theorem} {\em \cite[Theorem 9.2] {AWY}} \label{AWYthm}
Suppose that $\lambda_{1},...,\lambda_{n} \in \mathbb{D}$ are distinct points and $A_{k}=[a_{ij}^{k}] \in \Sigma$  are such that $a_{11}^{k}a_{22}^{k}\neq \det(A_{k})$, $1 \leq k \leq n$. The following conditions are equivalent.
\begin{enumerate}
         \item There exists an analytic function $ F: \mathbb{D} \rightarrow  \Sigma$ such that $F(\lambda_{k}) = A_{k}$, $1 \leq   k \leq n$;
              \item There exists an analytic function $\varphi:\mathbb{D} \rightarrow \mathbb{E}$ such that $\varphi (\lambda_{k})=\pi (A_{k})$, that is,

$$\varphi(\lambda_{k}) = (a_{11}^{k}, a_{22}^{k}, \det(A_{k})), \quad k=1,2,...,n .$$ 
            
\end{enumerate}
\end{theorem}

In the following theorem the authors give a necessary and sufficient condition for the solvability of a $\mu_{\rm Diag}$-synthesis problem by a rational $\mathbb{\overline{E}}$-inner function. 

\begin{theorem}{\em \cite [Theorem 1.1 and Theorem 8.1] {BLY16}} \label{solv}
Let $\lambda_{1},..,\lambda_{n}$ be distinct points in $\mathbb{D}$ and let $A_{k}=[a_{ij}^{k}] \in \mathbb{ C}^{2 \times 2}$ be such that $a_{11}^{k}a_{22}^{k}\neq \det(A_{k})$, $1 \leq k \leq n$. Let $$(x_{1}^{k}, x_{2}^{k}, x_{3}^{k})=\big(a_{11}^{k},a_{22}^{k}, \det (A_{k})\big), \quad 1 \leq k \leq n.$$ The following two conditions are equivalent.
\begin{enumerate}
\item There exists an analytic $2 \times 2$ matrix function $F$ in $\mathbb{D}$ such that 

\begin{equation*}
F(\lambda_{k})= A_{k} \quad \mbox{for} \quad k=1,..,n, 
\end{equation*}
and
\begin{equation*}
\mu_{\rm Diag}(F(\lambda)) \leq 1 \quad \text{for all}  \quad \lambda \in \mathbb{D};
\end{equation*}

\item there exists a {\em rational} $\mathbb{\overline{E}}$-inner function $x: \mathbb{D} \rightarrow \overline{\mathbb{E}}$ such that
\begin{equation*}
x(\lambda_{k})=(x_{1}^{k}, x_{2}^{k}, x_{3}^{k}) \quad \mbox{for} \quad k=1,..,n.
\end{equation*}
\end{enumerate}
\end{theorem}

\noindent Therefore, the understanding of rational $\mathbb{\overline{E}}$-inner functions will be useful for such $\mu$-synthesis problems.

 \subsection{The distinguished boundary of the tetrablock}\label{dist-bound-E}

\begin{theorem} {\em \cite[Theorem 2.9]{AWY}} \label{polcon}
$\mathbb{\overline{E}}$ is polynomially convex. 
\end{theorem}
Therefore, there exists a distinguished boundary $b\mathbb{\overline{E}}$ \index{$b\mathbb{\overline{E}}$} of $\mathbb{E}$. \index{$A(\mathbb{E})$} Let $A(\mathbb{E})$ be the algebra of continuous scalar functions on $\mathbb{\overline{E}}$ that are holomorphic on $\mathbb{E}$ endowed with the supremum norm. If there is a function $f \in A(\mathbb{E})$ and a point $p$ in $\overline{\mathbb{E}}$ such that $f(p) =1 $ and $|f(x)| <1$ for all $x \in \overline{\mathbb{E}}  \backslash \{p\}$, then $p \in b\mathbb{\overline{E}}$ and is called a \textit{peak point} \textit{\index{peak point}}of $\overline{\mathbb{E}}$, and the function $f$ is called a \textit{peaking function}  for $p$. 

\begin{theorem} {\em \cite[Theorem 7.1]{AWY}} \label{imp}
For $x \in \mathbb{C}^3$ the following are equivalent.
\begin{enumerate}
  \item $x_{1}=\overline{x}_2x_{3}, |x_{3}|=1$ and $ |x_{2}| \leq1$;
  \item either $x_{1}x_{2} \neq x_{3}$ and $\Psi(.,x)$ is an automorphism of $\mathbb{D}$ or $x_{1}x_{2} = x_{3}$ and $|x_{1}|=|x_{2}|=|x_{3}|=1$;
   \item $x$ is a peak point of $\overline{\mathbb{E}}$;
    \item there exists a $2\times2$ unitary matrix $U$ such that $x=\pi(U)$;
     \item there exists a symmetric $2\times2$  unitary matrix $U$ such that $x=\pi(U)$;
       \item $x \in  b\mathbb{\overline{E}}$;
        \item  $x \in \overline{\mathbb{E}}$ \;and $ |x_{3}|=1$.
\end{enumerate}
\end{theorem}

\begin{lemma} \label{ddb}
Let $x=(x_{1}, x_{2}, x_{3}) \in \mathbb{C}^{3}$. Then $x \in b \mathbb{\overline{E}}$ if and only if  $$x_{2}=\overline{x}_1x_{3}, \enspace |x_{3}|=1 \enspace \text{and} \quad  |x_{1}| \leq1.$$
\end{lemma}

\begin{proof}
By Theorem \ref{imp} (1),
\begin{equation*}
x \in b\mathbb{\overline{E}} \quad \Leftrightarrow \quad x_{1}=\overline{x}_{2}x_{3}, \ |x_{3}|=1 \ \text{and} \  |x_{2}| \leq1.
\end{equation*}
Since $|x_{3}|=1$ this implies $\overline{x}_{3}x_{3}=1$. Now, since $x \in b\mathbb{\overline{E}}$,
\begin{eqnarray*}
x_{1}&=&\overline{x}_{2}x_{3}, \enspace \text{and so} \quad 
\overline{x}_{1}=x_{2}\overline{x}_{3}.
\end{eqnarray*}
\text{Thus} \enspace $\overline{x}_{1}x_{3}=x_{2}\overline{x}_{3}x_{3}=x_{2}$.
Note, by Theorem \ref{tetbnd2}, $|x_{1}| \leq1$. \newline 
\indent Conversely, if 
$$x_{2}=\overline{x}_{1}x_{3}, \quad |x_{3}|=1 \quad \text{and} \quad  |x_{1}| \leq1$$
then, as in the previous steps, one can show that 
$x \in b\mathbb{\overline{E}}$. Therefore,
 $$x \in b\mathbb{\overline{E}} \enspace \text{if and only if} \enspace  x_{2}= \overline{x_{1}}x_{3}, \enspace |x_{3}|=1 \enspace \text{and} \enspace |x_{2}| \leq 1.$$
\end{proof}

\section{The symmetrised bidisc and $\Gamma$-inner functions} \label{G-inner}

In Section \ref{sece}
we show that there exist useful relations  between $\Gamma$-inner functions and $\mathbb{ \overline{ E}}$-inner functions. Recall the definition of the symmetrised bidisc $\Gamma$.
\begin{definition} 
	The symmetrised bidisc \index{symmetrised bidisc} is the set 
	\begin{equation*}
	\mathbb{G } \ {\stackrel{\text{def}}{=}} \  \big\{ (z+w,zw) : |z|<1, |w|<1 \big\},
	\end{equation*}
	and its closure is  
	\begin{equation*}
	\Gamma \ {\stackrel{\text{def}}{=}} \ \big\{ (z+w,zw) : |z| \leq1, |w|\leq 1 \big\}.
	\end{equation*}
\end{definition}
 In 1995 Jim Agler and Nicholas Young started the study of the symmetrised bidisc with the aim of solving a robust control problem in $H^{\infty}$ control theory. Although, the aim has not yet been achieved, it turned out that the symmetrised bidisc has a rich structure and it has attracted 
the attention of specialists in 
the several complex variables and in operator theory. 
                                               
We will use the co-ordinates $(s,p)$ for points in the symmetrized bidisc $\mathbb{G }$, chosen to suggest `sum' and `product'. 
The following result  \cite[Proposition 3.2] {ALY13} provides practical  criteria for membership of $\mathbb{G }$, of the distinguished boundary $b\Gamma$ of $\Gamma$ and of the topological boundary $\partial \Gamma$ of $\Gamma$.

\begin{proposition} \label{sp} {\em \cite[Proposition 3.2] {ALY13} }
	Let $(s,p)$ belong to $\mathbb{C}^{2}$. Then 
	\begin{enumerate} 
		\item $(s,p)$ belongs to $\mathbb{G}$ if and only if \index{$\mathbb{G}$} 
		\begin{equation*}
		| s-\overline{s}p| < 1-|p|^{2};
		\end{equation*} 
		\item $(s,p) $ belongs to $\Gamma$ if and only if 
		\begin{equation*}
		|s| \leq 2   \quad \text{and} \quad | s-\overline{s}p| \leq 1-|p|^{2};
		\end{equation*}
		\item $(s,p)$ lies in $b\Gamma$ if and only if
		\begin{equation*}
		|p|=1 , \quad |s|\leq 2  \quad \text{and } \quad s-\overline{s}p = 0;
		\end{equation*}
		\item $(s,p) \in \partial \Gamma$ if and only if \begin{equation*} \index{$b\Gamma$}
		|s| \leq 2   \quad \text{and} \quad | s-\overline{s}p| = 1-|p|^{2}.
		\end{equation*}
	\end{enumerate}
	
\end{proposition}

$\Gamma$-inner functions were defined and studied in \cite{ALY13}.
\begin{definition}
A {\em $\Gamma$-inner function} \index{$\Gamma$-inner function}  is an analytic function $h: \mathbb{D} \rightarrow \Gamma$ such that the radial limit 
\begin{equation} \label{eqn312}
 \lim_{r \rightarrow 1^{-}} h (r \lambda)
\end{equation}
exists and belongs to $b\Gamma$ for almost all $\lambda \in \mathbb{T}$ with respect to Lebesgue measure.
\end{definition}
\noindent By Fatou's Theorem, the limit \eqref{eqn312} exists for almost all $\lambda \in \mathbb{ T}$.

\begin{definition} \label{symmetry}
	Let f be a polynomial of degree less than or equal to $n$, where $n\geq 0$. Then we define the polynomial $f^{\sim n}$ by 
	\begin{equation*}
	f^{\sim n}(\lambda) = \lambda^{n} \overline{f(1/\overline{\lambda})}. 
	\end{equation*} 
The polynomial $f^{\vee}$ is defined by 
	\begin{equation*}
	f^{\vee }(\lambda) =  \overline{f\big(\overline{\lambda}\big)}. 
	\end{equation*}
\end{definition} 
\begin{remark} {\em  One can see that \\
(1) $$f^{\sim n}(\lambda)=\lambda^{n} \overline{f^(1/\overline \lambda)}=\lambda^{n} f^{\vee }(1/\lambda).$$
(2) If $f$ is a polynomial of degree $k$, then, for $n\geq k$, $\big(f^{\sim n}\big)^{\sim n}(\lambda)=f(\lambda)$.}
\end{remark}

Algebraic and geometric aspects of rational $\Gamma$-inner functions were studied in \cite{ALY15}. We are going to use some results from the paper.

\begin{proposition}{\em \cite[Proposition 2.2] {ALY15}} \label{G}
Let $h=(s,p)$ be a rational $\Gamma$-inner function \index{$\Gamma$-inner function} of degree $n$. Then there exist polynomials $E$ and $D$ such that
\begin{enumerate}
		\item $\deg(E) ,\deg(D) \leq n$,
		\item $E^{\sim n} = E$,
		\item $D(\lambda) \neq 0 $ on $\overline{\mathbb{D}}$,
		\item $\left| E(\lambda) \right| \leq 2 \left| D(\lambda) \right|$ on $\overline{\mathbb{D}}$,
		\item $s=\dfrac{E}{D}$ on $\overline{\mathbb{D}}$,
		\item $p=\dfrac{D^{\sim n}}{D}$ on $\overline{\mathbb{D}}$.
\end{enumerate}
Furthermore, $E_{1}$ and $D_{1}$ is a second pair of polynomials satisfying conditions {\em (1)--(6)} if and only if there exists a nonzero $t \in \mathbb{R}$ such that 
	\begin{equation*}
	E_{1}= tE \qquad \text{and} \qquad D_{1}=tD.
	\end{equation*} 
	Conversely, if $E$ and $D$ are polynomials which satisfy conditions {\em (1), (2), (4)}, $D(\lambda) \neq 0$ on $\mathbb{ D}$, and $s$ and $p$ are defined by equations {\em (5)} and {\em (6)}, then $h=(s,p)$ is a rational $\Gamma$-inner function of degree less than or equal to $n$.
\end{proposition}

The royal variety $\mathcal{R}_{\Gamma}$ \index{$\mathcal{R}_{\Gamma}$} of the symmetrised bidisc is 
\begin{equation*}
\mathcal{R}_{\Gamma}= \{(s,p) \in \mathbb{ C}^{2}: s^{2}=4p \}.
\end{equation*}
\begin{definition} {\em \cite[Page 7]{ALY15}} \label{Royal}
	Let $h=(s,p)$ be a $\Gamma$-inner function of degree $n$. Let $E$ and $D$ be as in Proposition {\em \ref{G}}. The royal polynomial $R_{h}$ \index{$R_{h}$}of $h$ is defined by \index{royal polynomial of $h$}
	\begin{equation*}
	R_{h}(\lambda)= 4D(\lambda) D^{\sim n}(\lambda) - E(\lambda)^{2}.
	\end{equation*}
\end{definition}

\begin{definition} {\em \cite[Definition 3.6]{ALY15}} \label{def5555}
	Let $h$ be a rational $\Gamma$-inner function such that \newline $h(\overline{\mathbb{D}}) \nsubseteq \mathcal{R}_{\Gamma} \cap \Gamma$. Let $R_{h}$ be the royal polynomial of $h$, and let  $\sigma$ be a zero of $R_{h}$ of order $\ell$. We define the {\em multiplicity} $\# \sigma$ \index{$\# \sigma$ } of $\sigma$ (as a royal node of $h$)  by 
		$$
	\# \sigma =
	\begin{cases}
	\ell             \quad \quad \quad \ \ \       \text{if} \ \sigma \in \mathbb{D},\\
	\frac{1}{2} \ell \quad \quad \quad \  \text{if} \ \sigma \in \mathbb{T} .
	\end{cases}
	$$
We define the type of $h$ to be the ordered pair $(n,k)$, where $n$ is the sum of the multiplicities of the royal nodes of $h$ that lie in $\mathbb{\overline{D}}$, and $k$ is the sum of the multiplicities of the royal nodes of $h$ that lie in $\mathbb{T}$. We define $\mathcal{R}_{\Gamma}^{n,k}$ \index{$\mathcal{R}_{\Gamma}^{n,k}$} to be the collection of rational $\Gamma$-inner functions $h$ of type $(n,k)$. 
\end{definition}

\begin{theorem} {\em \cite[Theorem 3.8]{ALY15}} Let $h \in \mathcal{R}_{\Gamma}^{n,k}$ be nonconstant.  Then $\deg(h)=n$.
\end{theorem}

\section{Rational $\mathbb{\overline{E}}$-inner functions} \label{E-inner_functions}
                                                              
In this section we give a definition of  the degree of a rational tetra-inner function $x$ by means of the fundamental group $\pi_{1}$. Recall that the rational inner functions on $\mathbb{ D}$ of degree $n$ are exactly the finite Blaschke products of degree $n$. Similar to this description of rational inner functions on $\mathbb{  D}$ we give an algorithm for the construction of
 all rational $\mathbb{ \overline{ E}}$-inner functions on $\mathbb{ D}$ in Theorem \ref{min}. In \cite{ALY15}, the authors describe all rational $\Gamma$-inner functions (see Proposition \ref{G}). We use this description and the connection between $\Gamma$-inner functions and $\mathbb{ \overline{ E}}$-inner functions to describe all rational $\mathbb{ \overline{ E}}$-inner functions on $\mathbb{ D}$.

\subsection{Relations between $\mathbb{\overline{E}}$-inner functions and $\Gamma$-inner functions} \label{sece}

\begin{definition} 
An {\em $\mathbb{\overline{E}}$-inner or tetra-inner function} is a map $f: \mathbb{D} \rightarrow \overline{\mathbb{E}}$ that is analytic and is such that the radial limit 
\begin{equation*} \label{tetrainner}
 \lim_{r \rightarrow 1^{-}} f (r \lambda)
\end{equation*}
exists and belongs to b$\mathbb{\overline{E}}$ for almost all $\lambda \in \mathbb{T}$ with respect to Lebesgue measure.
\end{definition}

\begin{remark} 
{\em	 Let $x: \mathbb{D} \rightarrow \mathbb{\overline{E}}$ be a rational $\mathbb{\overline{E}}$-inner function. Since $x$ is rational and bounded on $\mathbb{D}$ it has no poles in $\mathbb{ \overline{ D}}$ and hence $x$ is continuous on $\mathbb{ \overline{ D}}$.
		Thus one can consider the continuous function
		\begin{equation*}
		\tilde{x}: \mathbb{T} \rightarrow b\mathbb{\overline{E}}, \qquad \mbox{where} \quad \tilde{x}(\lambda)= \lim_{r\rightarrow 1^{-}} x(r\lambda) \quad \mbox{for all} \ \lambda \in \mathbb{ T}.
		\end{equation*}
In future we will use the same notation $x$ for both continuous functions $x$ and $\tilde{x}$.
}
\end{remark}

\begin{lemma} \label{x3inner}
Let $x=(x_{1},x_{2},x_{3})$ be an $\mathbb{\overline{E}}$-inner function. Then
\begin{enumerate}
	\item $x_{1}(\lambda)=\overline{x_{2}(\lambda)}x_{3}(\lambda)$, $|x_{2}(\lambda)|\leq 1$ and $|x_{3}(\lambda)|=1$ for almost all $\lambda \in \mathbb{ T}$;
	\item $x_{3}$ is an inner function on $\mathbb{D}$.
\end{enumerate}
\end{lemma}

\begin{proof}
(1) By the definition of $\mathbb{\overline{E}}$-inner function  $$x(\lambda) = (x_{1}(\lambda), x_{2}(\lambda), x_{3}(\lambda)) \in b\mathbb{\overline{E}}, \quad \text{for almost every} \ \lambda \in \mathbb{T}$$   and, by Theorem \ref{imp}, 
$$x_{1}(\lambda) = \overline{x_{2}}(\lambda)x_{3}(\lambda), \quad |x_{3}(\lambda)|=1 \quad \text{and} \quad |x_{2}(\lambda)| \leq 1 \quad \text{for almost all } \lambda \in \mathbb{T}. $$
(2) Since $$ x_{3}: \mathbb{D} \rightarrow \mathbb{\overline{D}} \quad \text{and}, \quad  \text{for almost all} \ \lambda \in \mathbb{T}, \ |x_{3}(\lambda)|=1,$$  $x_{3}$ is an inner function.
\end{proof}

\begin{remark} \label{x3innerblaschke} {\em
Let $x=(x_{1}, x_{2}, x_{3} )$ be a rational $\mathbb{\overline{E}}$-inner function. By Lemma {\ref{x3inner}}, $x_{3}$ is an inner function on $\mathbb{D}$, and so $x_{3}$ is a finite Blaschke product. }
\end{remark}

In \cite{TB} the author shows that there is a relation between points in the symmetrised bidisc and the tetrablock as follows.

\begin{lemma} {\em \cite[Lemma 3.2]{TB}} \label{TB1}
A point $x=(x_{1},x_{2},x_{3}) \in \mathbb{C}^{3}$ belongs to the tetrablock if and only if the pair $(x_{1}+zx_{2},zx_{3})$ is in the symmetrised bidisc $\mathbb{G}$ for every $z \in \mathbb{T}$.
\end{lemma}

\begin{proof} 
By Proposition \ref{sp} (1), $(s,p) \in \mathbb{G}$ if and only if 
\begin{equation} \label{s-sp}
| s-\overline{s}p| < 1-|p|^{2}.
\end{equation}
Suppose that $x=(x_{1},x_{2},x_{3}) \in \mathbb{E}$, $s_{z}=x_{1}+zx_{2}$ \; and $p_{z}=zx_{3}$.
\begin{eqnarray*}
|s_{z}-\overline{s_{z}}p_{z}|&=&|x_{1}+zx_{2}-\overline{(x_{1}+zx_{2})} zx_{3}|\\
&=&|x_{1}-\overline{x_{2}}x_{3}+z(x_{2}-\overline{x_{1}}x_{3})|  \\
&\leq&|x_{1}-\overline{x_{2}}x_{3}| + |x_{2}-\overline{x_{1}}x_{3}|, \qquad \text{since} \ |z|=1, \\
&<&1-|x_{3}|^{2}=1-|p_{z}|^{2}, \quad \text{by Theorem \ref{tetbnd}}.
\end{eqnarray*}
Hence $(s_{z},p_{z}) \in \mathbb{G}$.

Conversely, let $x=(x_{1}, x_{2}, x_{3}) \in \mathbb{C}^{3}$ and, for $z \in \mathbb{T}$, let  
 \begin{equation} \label{925}
s_{z}= x_{1}+zx_{2} \quad \text{and} \quad p_{z}= zx_{3}.
\end{equation}
Suppose for all $z \in \mathbb{ T}$, we have $(s_{z},p_{z}) \in \mathbb{G}$.
We want to show that $x=(x_{1}, x_{2}, x_{3}) \in \mathbb{E}$. Let us prove that
\begin{eqnarray*}
|x_{1}-\overline{x_{2}}x_{3}| + |x_{2}-\overline{x_{1}}x_{3}| < 1-|x_{3}|^{2}.
\end{eqnarray*} 
By assumption for all $z \in \mathbb{ T}$, $|s_{z}- \overline{s_{z}}p_{z}| < 1-|x_{3}|^{2}$. By equations \eqref{925}, we have
\begin{equation} \label{eq953}
|x_{1}-\overline{x_{2}}x_{3} + z\big(x_{2}-\overline{x_{1}}x_{3}\big)| < 1-|x_{3}|^{2}, \quad \text{for all} \ z \in \mathbb{ T}.
\end{equation} 
Let $$ \begin{cases}
z=e^{i\theta} \hspace{4cm} \ \ \theta \in (0,2\pi]; \\
w_{1}= x_{1}- \overline{x_{2}}x_{3} = |w_{1}|e^{i\theta_{1}} \qquad \theta_{1}  \in (0,2\pi]; \\
w_{2}=x_{2}-\overline{x_{1}}x_{3}  = |w_{2}|e^{i\theta_{2}} \qquad \theta_{2}  \in (0,2\pi].
\end{cases}$$
Now substitute $z$, $w_{1}$ and $w_{2}$ in inequality \eqref{eq953} 
$$\big| |w_{1}|e^{i\theta_{1}} + e^{i\theta} (|w_{2}|e^{i \theta_{2}}) \big| < 1-|x_{3}|^{2}.$$
This implies that
$$\big| |w_{1}|e^{i\theta_{1}} + |w_{2}|e^{i(\theta + \theta_{2})} \big| < 1-|x_{3}|^{2}, \quad \text{for all} \ e^{i\theta}. $$
We can choose $\theta$ such that $\theta + \theta_{2} =\theta_{1}$, that is, $\theta = \theta_{1} - \theta_{2}$. Hence 
$$ \big| |w_{1}|e^{i\theta_{1}} + |w_{2}|e^{i \theta_{1}} \big| = |e^{i\theta_{1}}| (|w_{1}| + |w_{2}|) = |w_{1}| + |w_{2}| = |x_{1}-\overline{x_{2}}x_{3}| + |x_{2}-\overline{x_{1}}x_{3}| < 1-|x_{3}|^{2}. $$
By Theorem \ref{tetbnd}, $(x_{1}, x_{2}, x_{3}) \in \mathbb{E }.$
\end{proof}

\begin{lemma} \label{TB2}
A point $x=(x_{1} ,x_{2}, x_{3}) \in \mathbb{C}^{3}$ belongs to the closed tetrablock if and only if for every $a \in \mathbb{\overline{D}}$,  $(ax_{1} + \overline{a}x_{2}, x_{3}) \in \Gamma$.
\end{lemma}

\begin{proof}
Suppose $x=(x_{1},x_{2},x_{3}) \in \mathbb{\overline{E}}$. Consider $(s_{a},p_{a}) = (ax_{1}+\overline{a}x_{2}, x_{3})$.
By Proposition \ref{sp} (2), $(s,p) \in \Gamma$ if and only if \begin{align}
| s-\overline{s}p| \leq 1-|p|^{2} \; \text{and} \; |s|\le 2.
\end{align}
\begin{eqnarray} \nonumber \label{eq998}
|s_{a}-\overline{s_{a}}p_{a}|&=&|ax_{1}+\overline{a}x_{2}-\overline{(ax_{1}+\overline{a}x_{2})}x_{3}|\\ \nonumber
&=&|a(x_{1}-\overline{x_{2}}x_{3})+\overline{a}(x_{2}-\overline{x_{1}}x_{3})|\\ \nonumber
&\leq& |a(x_{1}-\overline{x_{2}}x_{3})|+|\overline{a}(x_{2}-\overline{x_{1}}x_{3})|,  \\ \nonumber
&\leq&|x_{1}-\overline{x_{2}}x_{3}| + |x_{2}-\overline{x_{1}}x_{3}|, \qquad \text{since} \ |a| \leq 1, \\
&\leq& 1-|x_{3}|^{2}, \qquad \text{by Theorem \ref{tetbnd2}}.
\end{eqnarray}
Thus,
 $|s_{a}-\overline{s_{a}}p_{a}| \leq 1-|x_{3}|^{2}=1-|p_{a}|^{2}$ and
$|s_{a}| = |ax_{1} + \overline{a}x_{2}| \le 2.$
Hence $(s_{a},p_{a}) \in \Gamma $.

Conversely, let $x=(x_{1}, x_{2}, x_{3}) \in \mathbb{C}^{3}$.
	Suppose for every $a \in \mathbb{\overline{D}}$, we have $(s_{a},p_{a}) \in \Gamma$ where 
	\begin{equation} \label{926}
	s_{a}= ax_{1}+\overline{a}x_{2} \quad \text{and} \quad p_{a}=x_{3}.
	\end{equation}
	By equations \eqref{eq998} and \eqref{926}, we have
	\begin{equation} \label{eq927} 
	|s_{a}-\overline{s_{a}}p_{a}| =|a(x_{1}-\overline{x_{2}}x_{3}) + \overline{a}\big(x_{2}-\overline{x_{1}}x_{3}\big)|\leq 1-|p_{a}|^{2} , \quad \text{for all} \ a \in \mathbb{ \overline{D}}. 
	\end{equation} 
	Take $a \in \mathbb{T}$, then $$ \begin{cases}
	a=e^{i\theta} \hspace{4cm} \ \ \theta \in (0,2\pi]; \\
	w_{1}= x_{1}- \overline{x_{2}}x_{3} = |w_{1}|e^{i\theta_{1}} \qquad \theta_{1}  \in (0,2\pi]; \\
	w_{2}=x_{2}-\overline{x_{1}}x_{3}  = |w_{2}|e^{i\theta_{2}} \qquad \theta_{2}  \in (0,2\pi].
	\end{cases}$$
	Substitute $a$, $w_{1}$ and $w_{2}$ into inequality \eqref{eq927}, to get
	\begin{eqnarray*}
	|a(x_{1}-\overline{x_{2}}x_{3}) + \overline{a}\big(x_{2}-\overline{x_{1}}x_{3}\big)|= \big| e^{i\theta}|w_{1}|e^{i\theta_{1}} + e^{-i\theta}|w_{2}|e^{i\theta_{2}}\big|&=&\big| |w_{1}| e^{i(\theta + \theta_{1})} + |w_{2}|e^{i(\theta_{2}-\theta)}  \big| \\
	&\leq& 1-|x_{3}|^{2},
	\end{eqnarray*} for every $\theta \in (0,2\pi]$. Now choose $\theta = \dfrac{\theta_{2}-\theta_{1}}{2}$ to get
	\begin{eqnarray*}
	|x_{1}- \overline{x_{2}}x_{3}|+|x_{2}-\overline{x_{1}}x_{3}|&=& |w_{1}|  + |w_{2}| \\
	&=& |e^{i(\frac{\theta_{2}+\theta_{1}}{2})}| \big(|w_{1}|  + |w_{2}| \big) \\
	&=& \big| |w_{1}| e^{i(\frac{\theta_{2}+\theta_{1}}{2})} + |w_{2}|e^{i(\frac{\theta_{2}+\theta_{1}}{2})}\big| \\
	 &=& \big| |w_{1}|e^{i(\frac{\theta_{2}-\theta_{1}}{2} + \theta_{1})} + |w_{2}|e^{i(\theta_{2}-\frac{\theta_{2}-\theta_{1}}{2})} \big|   \leq 1-|x_{3}|^{2}.
	\end{eqnarray*}
Therefore $x=(x_{1}, x_{2},x_{3}) \in \mathbb{ \overline{ E}}$.

\end{proof}

\begin{lemma} \label{TB3}
	Let $s,p \in \mathbb{C}$ be such that $|s| \leq 2$ and $|p| \leq 1$. The pair $(s,p)$ belongs to $\Gamma$ if and only if $\left(\tfrac{1}{2}s,\tfrac{1}{2}s,p\right) \in \mathbb{ \overline{E }}$.
\end{lemma}

\begin{proof}	By Theorem \ref{tetbnd2},	
	$$\left(\tfrac{1}{2}s,\tfrac{1}{2}s,p\right)  \in \mathbb{\overline{E}} \quad  \Leftrightarrow \quad 
	 |\tfrac{1}{2}s-\tfrac{1}{2}\overline{s}p|+|\tfrac{1}{2}s-\tfrac{1}{2}\overline{s}p| \leq 1 - |p|^{2} . $$
Thus 
	 \begin{eqnarray} \nonumber
	(\tfrac{1}{2}s,\tfrac{1}{2}s,p) \in \mathbb{\overline{E}} &\Leftrightarrow& 2|\tfrac{1}{2}s-\tfrac{1}{2}\overline{s}p| \leq 1 -|p|^{2} \quad  \\ \nonumber
	&\Leftrightarrow& |s-\overline{s}p| \leq 1 -|p|^{2}.
	\end{eqnarray}
By assumption $|s| \leq 2$, hence by Proposition \ref{sp} (2), 
	 $$(\tfrac{1}{2}s,\tfrac{1}{2}s,p) \in \mathbb{\overline{E}}  \Leftrightarrow  (s,p) \in \Gamma .$$
\end{proof}

\begin{lemma} \label{12}
Let $x=(x_{1}, x_{2}, x_{3})$ be a rational $\mathbb{\overline{E}}$-inner function. Then

\begin{enumerate}

\item $h_{1}(\lambda)=\big(x_{1}(\lambda)+x_{2}(\lambda),x_{3}(\lambda) \big)$, for $\lambda \in \mathbb{D}$, is a rational $\Gamma$-inner function;
  \item $h_{2}(\lambda)= (ix_{1}(\lambda)-ix_{2}(\lambda),x_{3}(\lambda))$, for $\lambda \in \mathbb{D}$, is a rational $\Gamma$-inner function.

\end{enumerate}
\end{lemma}

\begin{proof}$(1)$
By Lemma \ref{TB1}, for all $\lambda \in \mathbb{D}$, $x(\lambda) \in \mathbb{E}$  implies that $$\big(x_{1}(\lambda)+x_{2}(\lambda),x_{3}(\lambda)\big) \in \mathbb{G}.$$ Consider $h_{1}=(s_{1},p_{1})$ where
$$  s_{1}(\lambda)=x_{1}(\lambda)+x_{2}(\lambda)  \quad \text{and}  \quad p_{1}(\lambda)=x_{3}(\lambda), \quad
 \text{for} \ \lambda \in \mathbb{D}.$$ 
It is obvious that $h_{1}$ is a rational function from $\mathbb{D}$ to $\mathbb{G}$. By assumption, $x$ is an $\mathbb{\overline{E}}$-inner function. Thus  $x(\lambda) \in b\mathbb{\overline{E}}$ for almost every $\lambda \in \mathbb{T}$. By Theorem \ref{imp} and Lemma \ref{ddb}, for almost all $\lambda \in \mathbb{T}$,
 \begin{equation} \label{dist-bound}
x_{2}(\lambda)=\overline{x_{1}(\lambda)}x_{3}(\lambda), \quad x_{1}(\lambda)=\overline{x_{2}(\lambda)}x_{3}(\lambda) ,\quad |x_{3}(\lambda)|=1 \quad \text{and} \, \, |x_{2}(\lambda)| \leq 1. 
\end{equation}
It is clear that $$|p_{1}(\lambda)|=|x_{3}(\lambda)|=1 \qquad \text{for} \ \lambda \in \mathbb{ T},$$  and, for almost all $\lambda \in \mathbb{T}$,
 \begin{eqnarray*} 
 |s_{1}(\lambda)|&=& |x_{1}(\lambda)+x_{2}(\lambda)| \\ 
  &\leq& |x_{1}(\lambda)|+|x_{2}(\lambda)| \leq 2.
 \end{eqnarray*}
 Since, for almost all $\lambda \in \mathbb{T}$, $x_{2}(\lambda)=\overline{x_{1}(\lambda)}x_{3}(\lambda)$, we have
\begin{eqnarray*} 
\overline {s_{1}(\lambda)} p_{1}(\lambda) &=&[\overline{ x_{1}(\lambda)} +  \overline{x_{2}(\lambda)}]x_{3}(\lambda)\\
&=&\overline{ x_{1}(\lambda)}x_{3}(\lambda) +  \overline{x_{2}(\lambda)}x_{3}(\lambda), \quad \text{by equations \eqref{dist-bound},}\\
&=& x_{1}(\lambda) + x_{2}(\lambda)= s_{1}(\lambda).
\end{eqnarray*}
Hence $s_{1}(\lambda)=\overline{s_{1}(\lambda)}p_{1}(\lambda)$ for almost every $\lambda \in {\mathbb{T}}.$
Therefore, by Proposition \ref{sp} (3), $h_{1}$ is a rational $\Gamma$-inner function. \\

\noindent $(2)$ Following the same steps as (1), let $h_{2}(\lambda)=(s_{2}(\lambda),p_{2}(\lambda)),$ where $$s_{2}(\lambda)=ix_{1}(\lambda)-ix_{2}(\lambda) \qquad  \text{and} \qquad  p_{2}(\lambda)=x_{3}(\lambda), \lambda \in \mathbb{D}.$$ By Lemma \ref{TB2}, $h_{2}$ is rational function from $\mathbb{D}$ to $\mathbb{G}$. Since $x$ is an $\mathbb{\overline{E}}$-inner function, $x(\lambda) \in b\mathbb{\overline{E}}$ for almost all $\lambda \in \mathbb{T}$. 
\noindent By Proposition \ref{sp}, to prove that $h_{2}$ is a rational $\Gamma$-inner function we need to show that
$$|p_{2}(\lambda)|=1, \, \, |s_{2}(\lambda)| \leq 2 \, \, \text{and} \, s_{2}(\lambda)= \overline{s_{2}(\lambda)}p_{2}(\lambda) \quad \text{for almost every } \lambda \in \mathbb{T}.$$
By Theorem \ref{imp}, for almost all $\lambda \in \mathbb{T}$, $|p(\lambda)|=|x_{3}(\lambda)|=1$ and 
\begin{eqnarray*}
|s_{2}(\lambda)|  \leq |ix_{1}(\lambda)|+|ix_{2}(\lambda)| \leq 2.
\end{eqnarray*}
By Lemma \ref{ddb}, $x_{2}(\lambda)=\overline{x_{1}(\lambda)}x_{3}(\lambda)$ for almost all $\lambda \in \mathbb{T}$.
Hence, for almost all  $\lambda \in \mathbb{T}$, 
\begin{eqnarray*}
\overline {s_{2}(\lambda)} p_{2}(\lambda) &=&[\overline{ ix_{1}(\lambda)} -  \overline{ix_{2}(\lambda)}]x_{3}(\lambda)\\
&=& i\big(\overline{x_{2}(\lambda)}\big)x_{3}(\lambda) - i\big(\overline{x_{1}(\lambda)}\big)x_{3}(\lambda), \quad \text{ by equations \eqref{dist-bound}},\\
&=& ix_{1}(\lambda) - ix_{2}(\lambda) = s_{2}(\lambda).
\end{eqnarray*}
Hence $s_{2}(\lambda)=\overline{s_{2}(\lambda)}p_{2}(\lambda)$ for almost every $\lambda \in {\mathbb{T}}.$
Therefore, by Proposition \ref{sp}, $h_{2}$ is a rational $\Gamma$-inner function. 
\end{proof}

\begin{lemma} \label{lem16}
	Let $x=(x_{1},x_{2},x_{3})$ be a rational $\mathbb{\overline{E}}$-inner function. Then 
	\begin{equation*}
	x_{1}(\lambda)=x_{2}^{\vee}(1/ \lambda)x_{3}(\lambda) \quad \text{for all} \ \ \lambda \in \mathbb{C} \setminus \{0 \}.
	\end{equation*}
\end{lemma}

\begin{proof}	By Theorem \ref{imp}, for all $\lambda \in \mathbb{ T}$,
	$$ x_{1}(\lambda)=\overline{x_{2}(\lambda)}x_{3}(\lambda).$$ For $\lambda \in \mathbb{ T}$, we have $|\lambda|=1$, that is, $\lambda \overline{\lambda}=1$, and so 
	\begin{equation*}
	\overline{x_{2}(\lambda)}= x_{2}^{\vee}(\overline{\lambda}) = x_{2}^{\vee}(\tfrac{1}{\lambda}).
	\end{equation*}
	Therefore, for all $\lambda \in \mathbb{T}$, 
	\begin{equation*}
	x_{1}(\lambda)=x^{\vee}_{2}(1/\lambda) x_{3}(\lambda).
	\end{equation*}
	Since $x_{1},x_{2},x_{3}$ are rational functions,
	\begin{equation*}
	x_{1}(\lambda)=x^{\vee}_{2}(1/\lambda) x_{3}(\lambda) \qquad \text{for all} \ \lambda \in \mathbb{C}\setminus \{0 \}.
	\end{equation*}
\end{proof}

\begin{proposition}
	Let $x=(x_{1}, x_{2}, x_{3})$ be a rational $\mathbb{ \overline{E }}$-inner function
	\begin{enumerate} 
		\item If $a  \in \mathbb{C}  \cup \{\infty\}$ is a pole of $x_{3}$ of multiplicity $k \geq 0$ and $\frac{1}{\overline{a}}$ is a zero of $x_{2}$ of multiplicity $\ell \geq 0$, then $a$ is a pole of $x_{1}$ of multiplicity at least $k- \ell$.
		\item If $a  \in \mathbb{C}  \cup \{\infty\}$ is a pole of $x_{1}$ of multiplicity $k \geq 1$, then $a$ is a pole of $x_{3}$ of multiplicity at least $k$.
	\end{enumerate}
	
\end{proposition}

\begin{proof}
	(1) By Lemma \ref{lem16}, we have
	\begin{equation} \label{eq555}
	x_{1}(\lambda)= x_{2}^{\vee}(1/\lambda)x_{3}(\lambda) \qquad \mbox{for} \ \lambda \in \mathbb{C}\setminus \{0 \}.
	\end{equation}
	Since $x_{3}$ is a rational inner function, $x_{3}$ cannot have any pole in $\mathbb{ \overline{ D}}$. Hence $|a| >1$ and so $|\tfrac{1}{a}| < 1$. We know that $x_{2}^{\vee}$ is analytic in $\mathbb{D}$, so $\tfrac{1}{a}$ cannot be a pole of $x_{2}^{\vee}$. By equation \eqref{eq555},
	\begin{equation*}
	(\lambda - a)^{k-\ell -1} x_{1}(\lambda)= (\lambda - a)^{k-\ell -1} x_{2}^{\vee}(1/\lambda)x_{3}(\lambda).
	\end{equation*}
	Take the limit for both sides as $\lambda$ goes to $a$:
	\begin{equation*}
	\lim_{\lambda \rightarrow a} (\lambda - a)^{k-\ell -1} x_{1}(\lambda)= \lim_{\lambda \rightarrow a} (\lambda - a)^{k-\ell -1}  x_{2}^{\vee}(1/\lambda)x_{3}(\lambda).
	\end{equation*}
	The right hand side goes to $\infty$, therefore $x_{1}$ has a pole of multiplicity at least $k-l$ at $a$. \newline
	
	Now suppose that $\infty$ is a pole of $x_{3}$ of multiplicity $k$ and $0$ is a zero of $x_{2}$ of multiplicity $\ell$. By equation \eqref{eq555}, for all $\lambda \in \mathbb{C} \backslash \{0\}$, we have
	$$x_{1}(\tfrac{1}{\lambda})= x^{\vee}_{2}(\lambda)x_{3}(\tfrac{1}{\lambda}).$$
	Multiply both sides by $\dfrac{\lambda^{k-1}}{\lambda^{\ell}}$ to obtain the equation
	\begin{equation} \label{1212}
	 \dfrac{\lambda^{k-1}}{\lambda^{\ell}} x_{1}(\tfrac{1}{\lambda}) = \dfrac{\lambda^{k-1}}{\lambda^{\ell}}  x^{\vee}_{2}(\lambda)x_{3}(\tfrac{1}{\lambda}).
	 \end{equation}
	Since $x_{2}^{\vee}$ is analytic at $0$ and has a zero of multiplicity $\ell >0$ at $0$, we have
	$$ \lim_{\lambda \rightarrow 0} \frac{x_{2}^{\vee}(\lambda)}{\lambda^{\ell}} = c, \qquad \text{where} \ c \in \mathbb{C}\backslash \{0\}.$$
	Since by assumption, $x_{3}(\lambda)$ has a pole of multiplicity $k$ at $\infty$,
	$$\lim_{\lambda \rightarrow 0} \lambda ^{k-1} x_{3}(\tfrac{1}{\lambda})= \infty.$$
	Hence by equation \eqref{1212},
	$$ \lim_{\lambda \rightarrow 0} \lambda ^{k-\ell-1} x_{1}(\tfrac{1}{\lambda})= \infty.$$ It follows that $x_{1}(\frac{1}{\lambda})$ has a pole of multiplicity at least $k -\ell$ at $0$. That is, $x_{1}(\lambda)$ has a pole of multiplicity at least $k-\ell$ at $\infty$.\newline \\
	(2) Let $a \in \mathbb{C}$ be a pole of $x_{1}$ of multiplicity $k \geq 1$. Then $|a| > 1$. This implies $|\tfrac{1}{a}| < 1$. Therefore $x_{2}^{\vee}$ is analytic at $\tfrac{1}{a}$. Now 
	\begin{equation*}
	\lim_{\lambda \rightarrow a} (\lambda - a)^{k -1} x_{1}(\lambda) = \infty.
	\end{equation*}
	Thus $a$ is a pole of $x_{3}$ of multiplicity at least $k$.
	
	If $\infty$ is a pole of $x_{1}$ of multiplicity $k \geq 1$. Then $0$ is a pole of $x_{1}(\frac{1}{\lambda})$ of multiplicity $k$, that is, 
	$$\lim_{\lambda \rightarrow 0} \lambda^{k -1} x_{1}(\tfrac{1}{\lambda}) = \infty. $$
	By relation \eqref{1212},
	$$\lambda^{k-1} x_{1}(\tfrac{1}{\lambda}) = \lambda^{k-1}  x^{\vee}_{2}(\lambda)x_{3}(\tfrac{1}{\lambda}). $$
	Since $x_{2}^{\vee}$ is analytic at $0$, $0$ cannot be a pole of $x_{2}^{\vee}$ and thus 
	$$\lim_{\lambda \rightarrow 0} x_{2}^{\vee}(\lambda)= x_{2}^{\vee}(0).$$
	Therefore
	$$\lim_{\lambda \rightarrow 0} \lambda^{k-1}x_{3} (\tfrac{1}{\lambda})= \infty.$$
	 This completes the proof that $x_{3}$ has a pole of multiplicity at least $k$ at $\infty$. 
\end{proof}

\subsection{The degree of a rational $\mathbb{\overline{E}}$-inner function}
\label{degree_E-inner}

Let us define the notion of the degree $\deg(x)$ of  a rational $\mathbb{\overline{E}}$-inner function $x$ by means of fundamental groups.

\begin{definition} \label{degdef}
	The degree $\deg(x)$ of a rational $\mathbb{\overline{E}}$-inner function $x$ is defined to be $x_{*}(1)$, where $ x_{*}:\mathbb{Z} =\pi_{1}(\mathbb{T}) \rightarrow \pi_{1}(b\mathbb{\overline{E}})$ is the homomorphism of fundamental groups induced by $x$ when $x$ is regarded as a continuous map from $\mathbb{T}$ to  $b\mathbb{\overline{E}}$. 
\end{definition}

We will assume that $\deg(x)$ is a non-negative integer.

\begin{lemma}
	$b\mathbb{\overline{E}}$ is homotopic to $\mathbb{T}$ and $\pi_{1}(b\mathbb{\overline{E}})= \mathbb{Z}$.
\end{lemma}
\begin{proof}
	 The maps
	\begin{align*}
	f&:b\mathbb{\overline{E}} \rightarrow \mathbb{T},\ \text{defined by} \ f(x_{1}, x_{2},x_{3})=x_{3},\\
	g&: \mathbb{T}  \rightarrow b\mathbb{\overline{E}}, \ \text{defined by} \  g(z)=(0, 0, z),
	\end{align*}
	satisfy 
	$$(g \circ f)(x_{1}, x_{2},x_{3})=g\big(f(x_{1}, x_{2},x_{3}) \big)=g(x_{3})=(0, 0,x_{3})$$
	and
	$$(f \circ g)(z)=f(0,0,z)=z, $$
	 that is, $f \circ g = \id_{\mathbb{T}}$. If $(x_{1}, x_{2},x_{3}) \in b\mathbb{\overline{E}}$ and $0 \leq t \leq1$, then $(tx_{1}, tx_{2},x_{3}) \in b\mathbb{\overline{E}}$. Let $I=[0,1]$. Consider the map 
	\begin{equation*}
h: b\mathbb{\overline{E}} \times I \rightarrow b\mathbb{\overline{E}}, 
	\end{equation*}
	which is defined by
	\begin{equation*}
	h(x_{1}, x_{2},x_{3},t)=(tx_{1}, tx_{2},x_{3}).
	\end{equation*}
	One can see that
	\begin{equation*}
	h(x_{1}, x_{2},x_{3},0)=(0x_{1}, 0x_{2},x_{3})=(0, 0,x_{3})=(g \circ f)(x_{1}, x_{2}, x_{3}) \ \text{and}
	\end{equation*}
	\begin{equation*}
	h(x_{1}, x_{2},x_{3},1)=(1x_{1}, 1x_{2},x_{3})=(x_{1}, x_{2},x_{3})
	= {\id_{b\mathbb{\overline{E}}}} (x_{1}, x_{2},x_{3}).
	\end{equation*}
	Therefore $h$ defines a homotopy between $g \circ f$ and $\id_{b\mathbb{\overline{E}}}$, that is, $g \circ f \simeq \id_{b\mathbb{\overline{E}}}$. Hence $b\mathbb{\overline{E}}$ is homotopically equivalent to $\mathbb{ T}$ and 
	it follows that $\pi_{1}(b\mathbb{\overline{E}})=\pi_{1}(\mathbb{T})=\mathbb{Z}$.
\end{proof}
\begin{lemma} \label{zeros}
	Let $B$ be a finite Blaschke product. Then the degree of $B$ is equal to $B_{*}(1)$.
\end{lemma}

\begin{proof} Since $B$ is a finite Blaschke product, it can be written as $$B(\lambda)= e^{i\theta} \prod_{j=1}^{N} \frac{\lambda -\alpha_{j}}{1-\overline{\alpha_{j}} \lambda}, \qquad \text{where} \ \alpha_{j} \in \mathbb{ D}, \ j=1 \dots N, \ \text{and} \ \theta \in [0,2\pi).$$ 
	One can consider the map, $B: \mathbb{T} \rightarrow \mathbb{T}$, and $$B_{*}: \pi_{1}(\mathbb{T}) = \mathbb{Z} \rightarrow \pi_{1}(\mathbb{T}) = \mathbb{Z}.$$ Now $1 \in  \pi_{1}(\mathbb{T})$ is the homotopy class of $\id_{\mathbb{T}}$ and 
	$B_{*}(1)$ is equal to the homotopy class of $B \circ \id_{\mathbb{T}} = B$,  when $B$ is regarded as a continuous map from $\mathbb{T}$ to $\mathbb{T}$. Therefore  $B_{*}(1)=n(\gamma,a)$, where $n(\gamma,a)$ is the winding number of $\gamma$ about $a$, which lies inside $\gamma= \{B(e^{it}) : 0 \leq t \leq 2\pi \}$. Thus, one can see that
	\begin{eqnarray} \nonumber
	n(\gamma,a) &=& \frac{1}{2\pi i}\int_{\gamma}^{} \frac{dz}{z-a}  \\ \nonumber
	&=& \frac{1}{2\pi i}\int_{\mathbb{T}}^{} \frac{B'(z)dz}{B(z)-a}.
	\end{eqnarray}
	By the Argument Principle, \cite[Theorem 18]{Ahl}, the integral
	\begin{equation*}
	\frac{1}{2\pi i}\int_{\mathbb{T}}^{}\frac{B'(z)}{B(z)-a}dz,
	\end{equation*}
	is equal to the number of zeros of $B$ in $\mathbb{D}$.
	It is clear that $B$ has $N$ zeros, counting multiplicities, and has degree  $N$.
	Therefore the number of zeros of $B$ is equal to the winding number of $\gamma$ about $a$, and it is equal to $N$. \newline 
\end{proof}

\begin{proposition} \label{deg}
	For any rational $\mathbb{\overline{E}}$-inner function $x=(x_{1},x_{2},x_{3})$, $\deg(x)$ is the degree $\deg(x_{3})$ {\rm (}in the usual sense{\rm )} of the finite Blaschke product $x_{3}$.
\end{proposition}
\begin{proof}
	Since $x$ is a rational $\mathbb{\overline{E}}$-inner function, $x_{3}$ is an inner function, and so $x_{3}$ is a finite Blaschke product. Two $\mathbb{\overline{E}}$-inner functions $x=(x_{1},x_{2},x_{3})$ and $y=(0,0,x_{3})$ are homotopic if there exists a continuous mapping $f: \mathbb{T} \times I \rightarrow b\mathbb{\overline{E}}$ such that
	\begin{equation*}
	f(\lambda,0)=y(\lambda) \qquad \text{and} \qquad f(\lambda,1)=x(\lambda), \; \lambda \in \mathbb{T}.
	\end{equation*}
	Let $$x^{t}(\lambda)=\big(tx_{1}(\lambda),tx_{2}(\lambda) ,x_{3}(\lambda)\big) \  \text{for} \ \lambda \in \mathbb{D} \ \text{and} \ t \in [0,1].$$ 
	Since $x(\lambda) \in b\mathbb{\overline{E}}$, for all $\lambda \in \mathbb{T}$, by Theorem \ref{imp} (1),
	$$x_{1}(\lambda)=\overline{x_{2}(\lambda)}x_{3}(\lambda) \quad \text{and} \quad |x_{3}(\lambda)|=1.$$
	Hence for all $\lambda \in \mathbb{T}$,
	
	\begin{equation*}
	tx_{1}(\lambda)=t\overline{x_{2}(\lambda)}x_{3}(\lambda).
	\end{equation*}
	Therefore,  $$x^{t}(\lambda)=\big(tx_{1}(\lambda),tx_{2}(\lambda) ,x_{3}(\lambda)\big) \in b\mathbb{\overline{E}} \ \ \text{for} \ \lambda \in \mathbb{T}.$$
	Hence $x^{t}$ is a homotopy between $x=x^{1}$ and $(0,0,x_{3})=x^{0}$. \newline
	It follows that the homomorphism
	$$x_{*}: \pi_{1}(\mathbb{T}) = \mathbb{Z} \rightarrow \pi_{1}(b\mathbb{\overline{E}})=\mathbb{Z}$$ coincides with $(x^{0})_{*}=(0,0,x_{3})_{*}$. By Lemma \ref{zeros}, $(x_{3})_{*}(1)=\deg x_{3}$, since $x_{3}$ is a finite Blaschke product. Therefore $(0,0,x_{3})_{*}(1)$ is the degree of the finite Blaschke product $x_{3}$. 
\end{proof}

\subsection{Description of rational $\mathbb{ \overline{ E}}$-inner functions}\label{desc_E-inner_funct}
\begin{theorem} \label{min}
Let $x=(x_{1}, x_{2}, x_{3})$ be a rational $\mathbb{ \overline{ E}}$-inner function of degree $n$. Then there exist polynomials $E_{1}, E_{2}, D$ such that
\begin{enumerate}
     \item $\deg(E_{1}) ,\deg(E_{2}), \deg(D) \leq n$,
       \item $D(\lambda) \neq 0 $ on $\overline{\mathbb{D}}$,
       \item $x_{3}=\frac{D^{\sim n}}{D}$ on $\overline{\mathbb{D}}$, 
         \item $x_{1}= \frac{E_{1}}{D}$ on $\overline{\mathbb{D}}$,
          \item $x_{2}= \frac{E_{2}}{D}$ on $\overline{\mathbb{D}}$,
           \item $|E_{i}(\lambda)| \leq |D(\lambda)|$ on $\overline{\mathbb{D}}$, \ for \  $i=1,2$,  
            \item $E_{1}(\lambda)= E_{2}^{\sim n}(\lambda), \ \mbox{for} \ \lambda \in \mathbb{\overline {D}}$.            
\end{enumerate}

Conversely, if $E_{1}, E_{2}$ and $D$ satisfy conditions {\em (1), (6) and (7)}, $D(\lambda) \neq 0$ on $\mathbb{D}$ and $x_{1}, x_{2}$ and $x_{3}$ are defined by equations {\em (3)--(5)}, then $x=(x_{1}, x_{2},x_{3})$ is a rational $\mathbb{\overline{E}}$-inner function of degree at most $n$.\\
Furthermore, a triple of polynomials $E_{1}^{1}, E_{2}^{1}$ and $D^{1}$ satisfies relations {\em (1)--(7)} if and only if there exists a real number $t \neq 0$ such that
\begin{equation*}
E_{1}^{1}=tE_{1}, \quad E_{2}^{1}=tE_{2} \quad \text{and} \quad D^{1}=tD.
\end{equation*}

\end{theorem}

\begin{proof}
By assumption $x= \big(x_{1}, x_{2}, x_{3}\big)$ is a rational $\mathbb{ \overline{ E}}$-inner function. By Lemma \ref{12} (1), $h_{1}= (s,p)$ where $s=x_{1}+x_{2} , p=x_{3}$ is a rational $\Gamma$-inner function. Since $x_{3}: \mathbb{D} \rightarrow \mathbb{D}$ is an inner function, it is a finite Blaschke product and, by \cite[Corollary 6.10]{ALY13}, it can be written in the form
\begin{equation*}
x_{3}(\lambda) = c\frac{\lambda^{k}D^{\sim (n-k)}(\lambda)}{D(\lambda)},
\end{equation*}

\noindent where $|c|=1, 0 \leq k \leq n$ and $D$ is a polynomial of degree $n-k$ such that $D(0)=1$. By Proposition \ref{G}, there exist polynomials $E,D$ such that 

\begin{enumerate}
     \item $\deg(E) ,\deg(D) \leq n$,
       \item $E^{\sim n} = E$,
         \item $D(\lambda) \neq 0 $ on $\overline{\mathbb{D}}$,
          \item $\left| E(\lambda) \right| \leq 2 \left| D(\lambda) \right|$ on $\overline{\mathbb{D}}$,
            \item $s=\frac{E}{D}$ on $\overline{\mathbb{D}}$,
             \item $p=\frac{D^{\sim n}}{D}$ on $\overline{\mathbb{D}}$.
\end{enumerate}
Hence 
\begin{equation} \label{nm}
x_{1}+x_{2}=s=\frac{E}{D} \quad and \quad x_{3}=p=\frac{D^{\sim n}}{D}.
\end{equation}
\noindent By Lemma \ref{12} (2), $h_{2}=(s_{2},p_{2})$, where $s_{2}=ix_{1}-ix_{2}$, $p_{2}=x_{3}=p_{1}$ is a rational $\Gamma$-inner function. By Proposition \ref{G}, for $h_{2}=(s_{2},p_{2})$, there exist polynomials $G, D$ such that
\begin{enumerate}
     \item $\deg(G) ,\deg(D) \leq n$,
       \item $G^{\sim n} = G$,
         \item $D(\lambda) \neq 0 $ on $\overline{\mathbb{D}}$,
          \item $\left| G(\lambda) \right| \leq 2 \left| D(\lambda) \right|$ on $\overline{\mathbb{D}}$,
            \item $s_{2}=ix_{1}-ix_{2}=\frac{G}{D}$ on $\overline{\mathbb{D}}$,
             \item $p_{2}=x_{3}=\frac{D^{\sim n}}{D}$ on $\overline{\mathbb{D}}$.
\end{enumerate}
\noindent Therefore, by equation (5),
\begin{equation} \label{eq1}
 x_{1}-x_{2}=-\frac{iG}{D}.
\end{equation}
\noindent By relation (\ref{nm}),
\begin{equation} \label{eq2}
x_{1}+x_{2}=\frac{E}{D} \ .
\end{equation} 
\noindent Add equations (\ref{eq1}) and (\ref{eq2}) to get 
$$x_{1}=\frac{\frac{1}{2}(E-iG)}{D}.$$
\noindent Substitution of  $x_{1}$ in equation (\ref{eq2}) gives 
\begin{equation*}
x_{2}=\frac{\frac{1}{2}(E+iG)}{D}.
\end{equation*}
\noindent Define the polynomials $E_{1}$ and $E_{2}$ by $$E_{1}=\frac{1}{2}(E-iG), \qquad E_{2}=\frac{1}{2}(E+iG).$$ Since the degrees of both polynomials  $E, G$ are at most $n$, $\deg(E_{1}), \deg(E_{2}) \leq n$. Thus, for $\lambda \in \mathbb{ \overline{ D}}$,
\begin{equation*}
x_{1}(\lambda)=\frac{ E_{1}(\lambda)}{D(\lambda)} \qquad \text{and} \qquad 
x_{2}(\lambda)=\frac{ E_{2}(\lambda)}{D(\lambda)}.
\end{equation*}  
Since $x$ is an $\mathbb{ \overline{ E}}$-inner function, for $\lambda \in \mathbb{\overline {D}}$,
\begin{eqnarray*} \nonumber
|x_{1}(\lambda)| \leq 1 \quad &\text{and}& \quad |x_{2}(\lambda)|\leq 1, \\ \text{and so} \qquad
|E_{1}(\lambda)| \leq |D(\lambda)| \qquad &\text{and}& \quad |E_{2}(\lambda)| \leq |D(\lambda)|.
\end{eqnarray*}
Hence $|E_{i}(\lambda)| \leq |D(\lambda)|$ on $\overline{\mathbb{D}}$, where $i=1,2$. Therefore
conditions (1)--(6) of Theorem \ref{min} are satisfied. \\

By assumption, $x$ is a rational $\mathbb{\overline{E}}$-inner function. Thus, for all $\lambda \in \mathbb{T}$,
\begin{eqnarray} \label{iner}
x_{1}(\lambda)= \overline{x_{2}}(\lambda)x_{3}(\lambda)
\Leftrightarrow \frac{E_{1}(\lambda)}{D(\lambda)} &=& \frac{\overline{E_{2}(\lambda)}}{\overline{D(\lambda)}} \times \frac{D^{\sim n}(\lambda)}{D(\lambda)} \nonumber \\ 
\Leftrightarrow \frac{E_{1}(\lambda)}{D(\lambda)} &=& \frac{\overline{E_{2}(\lambda)}}{D^{\vee}(1/\lambda)} \times  \frac{\lambda^{n} D^{\vee}(1/\lambda)}{D(\lambda)}, \   \text{since}  \ \ \overline{D}(\lambda)=D^{\vee}(\overline{\lambda})= D^{\vee}(1/\lambda). \nonumber \\  
\Leftrightarrow \frac{E_{1}(\lambda)}{D(\lambda)} &=& \frac{\lambda^{n} \overline{E_{2}(\lambda)}}{D(\lambda)} \nonumber \\
\Leftrightarrow \frac{E_{1}(\lambda)}{D(\lambda)} &=& \frac{E^{\sim n}_{2}(\lambda)}{D(\lambda)} \nonumber \\
\Leftrightarrow E_{1}(\lambda) &=& E^{\sim n}_{2}(\lambda) .
\end{eqnarray}
Hence $E_{1}(\lambda)= E_{2}^{\sim n}(\lambda)$ for all $\lambda \in \mathbb{T}$, and therefore on $\mathbb{\overline{ D}}$. Thus equation  (7) of Theorem \ref{min}
 is proved. \\

{\bf Let us prove the converse statement}. Let $E_{1}, E_{2}$ and $D$ satisfy relations (1), (6) and (7) of Theorem \ref{min} and $D(\lambda) \neq 0$ on $\mathbb{D}$, and $x_{1},x_{2},x_{3}$ be defined by equations (3)--(5), that is, $$x_{1}= \dfrac{E_{1}}{D}, \quad x_{2}= \dfrac{E_{2}}{D} \quad  \text{and} \quad x_{3}= \dfrac{D^{\sim n}}{D}.$$ 
Let us show that  $x=(x_{1},x_{2},x_{3})$ is a rational $\mathbb{\overline{E}}$-inner function. By Theorem \ref{imp}, we have to prove that $x: \mathbb{D} \rightarrow \mathbb{E}$ and the following conditions are satisfied.
\begin{enumerate}
     \item $|x_{3}(\lambda)|=1$ for almost all $\lambda$ on $\mathbb{T}$, that is, $x_{3}$ is inner,
       \item $|x_{2}| \leq 1$ on $\overline{\mathbb{D}}$,
         \item $x_{1}(\lambda)= \overline{x_{2}(\lambda)}x_{3}(\lambda)$ \ for almost all $\lambda \in \mathbb{T}$.\\
         
\end{enumerate}
\noindent $(1)$ Firstly, if $D$ has no zeros on the unit circle, then $D$ and $D^{\sim n}$ have no common factor. Therefore, $x_{3}(\lambda)=\dfrac{D^{\sim n}(\lambda)}{D(\lambda)}$ maps $\mathbb{T}$ to $\mathbb{T}$. Hence, $x_{3}$ is an inner function and 
\begin{equation*}
\deg(x_{3})= \deg \Bigg(\frac{D^{\sim n}}{D}\Bigg)\ = \max \{\deg(D^{\sim n}), \deg(D) \}=n.
\end{equation*}
Second case: if $D$ has the zeros $a_{1}, ...,a_{\ell}$ on $\mathbb{T}$ then $D$ and $D^{\sim n}$ have the common factor $\prod_{i=1}^{\ell} (\lambda - a_{i})$ and hence $x_{3}=\dfrac{D^{\sim n}}{D}$ is inner and 
\begin{equation*}
\deg(x_{3})= \deg \Bigg(\frac{D^{\sim n}}{D}\Bigg)\ \leq n-\ell.
\end{equation*}
$(2)$ \noindent By assumption (6), $$|E_{2}(\lambda)| \leq |D(\lambda)| \quad \text{for all} \ \lambda \in \overline{\mathbb{D}}.$$ This implies $|\frac{E_{2}(\lambda)}{D(\lambda)}| \leq 1$ and hence $|x_{2}(\lambda)| \leq 1$.\\

\noindent $(3)$ By assumption (7), $E_{1}(\lambda)= E_{2}^{\sim n}(\lambda)$, for almost all $\lambda \in \mathbb{T}$  and by the equality (\ref{iner}), $x_{1}(\lambda)=\overline{x_{2}(\lambda)}x_{3}(\lambda), \quad \text{for almost all} \ \lambda \in \mathbb{T}$.\\

\textbf{Let us show that $x=(x_{1}, x_{2}, x_{3})=\bigg(\dfrac{E_{1}}{D},\dfrac{E_{2}}{D},\dfrac{D^{\sim n}}{D}\bigg)$ maps $\mathbb{D}$ to $\mathbb{E}$}, that is, \newline $x(\lambda)=(x_{1}(\lambda), x_{2}(\lambda), x_{3}(\lambda)) \in \mathbb{E}$ for all $\lambda \in \mathbb{D}$. 
By Theorem \ref{tetbnd2}, for $\lambda \in \mathbb{D}$,
\begin{equation*}
x(\lambda) \in \overline{\mathbb{E}} \Leftrightarrow \|\Psi(.,x(\lambda))\|_{H^{\infty}} \leq 1,
\end{equation*}
where $\Psi(z,x) \ = \ \dfrac{x_{3}z-x_{1}}{x_{2}z-1}$. Note that, for every $z \in \mathbb{D}$,
 \begin{eqnarray*}
 \Psi(z,x) : &\mathbb{D}& \rightarrow \mathbb{C} \\ 
 :&\lambda& \rightarrow \Psi(z,x(\lambda))
 \end{eqnarray*}
   is analytic on $\mathbb{ D}$ because $x_{i} , i=1, 2, 3$, are analytic functions on $\mathbb{ D}$, and $|x_{2}(\lambda)|\leq 1$ and $x_{2}(\lambda)z \neq 1$ for all $\lambda \in \mathbb{ D}$. We have shown above that, for almost all $\lambda \in \mathbb{T}$, $x(\lambda) \in b\overline{\mathbb{E}}$. Thus, by Theorem \ref{imp} (2),
\begin{equation*}
x(\lambda) \in b\mathbb{\overline{E}} \enspace \text{if and only if } \enspace \Psi(.,x(\lambda) ) \ \text{is an automorphism of} \ \mathbb{D}.
\end{equation*}
By the maximum principle, for all $z, \lambda \in \mathbb{D}$, $|\Psi(z,x(\lambda) )|< 1$. Thus by Theorem \ref{imp}, $x(\lambda) \in \mathbb{E}$ for all $\lambda \in \mathbb{D}$.\newline

Suppose that $t$ is a nonzero real number and
\begin{equation*}
 E_{1}^{1}=tE_{1}, \quad E_{2}^{1}=tE_{1} \quad \text{and} \quad D^{1}=tD.
\end{equation*}
Then it is clear that $E_{1}^{1}, E_{2}^{1}$ and $D_{1}$ satisfy conditions (1)--(7). Conversely, let $E_{1}^{1}, E_{2}^{1}$ and $D^{1}$ be a second triple that satisfies relations (1)--(7). Then 
\begin{equation} \label{eql11}
x_{1}= \dfrac{E_{1}}{D}= \dfrac{E_{1}^{1}}{D^{1}} \quad \text{on} \  \overline{\mathbb{D}},
\end{equation}
\begin{equation} \label{eql12}
x_{2}= \dfrac{E_{2}}{D}= \dfrac{E_{2}^{1}}{D^{1}} \quad \text{on} \  \overline{\mathbb{D}},
\end{equation}
\begin{equation} \label{eql13}
x_{3}= \dfrac{D^{\sim n}}{D}= \dfrac{D^{1\sim n}}{D^{1}} \quad \text{on} \  \overline{\mathbb{D}}.
\end{equation}

Suppose that $D(\lambda )=a_{0} + a_{1}\lambda + ... + a_{k}\lambda^{k}$ \ where $a_{0} \neq 0$ and $k \leq n$. Then \newline
\begin{eqnarray*}
D^{\sim n}(\lambda)&=&\lambda^{n}\overline{D\big( 1/\overline{\lambda}}\big)\\
&=&\lambda^{n}\bigg(\overline{ a_{0} + \dfrac{a_{1}}{\overline{\lambda}} + ... + \dfrac{a_{k}}{\overline{\lambda}}} \bigg) \\
&=&\overline{a_{0}}\lambda^{n} + \overline{a_{1}}\lambda^{n-1}+ ... + \overline{a_{k}}\lambda^{n-k} .
\end{eqnarray*}
Thus, for all $\lambda \in \mathbb{D}$, 
\begin{equation*}
x_{3}=\dfrac{D^{\sim n}(\lambda)}{D(\lambda)}=\dfrac{\lambda^{n-k}\big( \overline{a_{0}}\lambda^{k} + \overline{a_{1}}\lambda^{k-1}+ ... + \overline{a_{k}} \big)}{a_{0} + a_{1}\lambda + ... + a_{k}\lambda^{k}}.
\end{equation*}
Therefore, $x_{3}$ has a zero of multiplicity $(n-k)$ at $0$ , has $k$ poles in $\mathbb{C}$, counting multiplicity, and has degree $n$. Hence the poles of $x_{3}$ in $\{z \in \mathbb{C}: |z| >1 \}$, $n$ and $k$ are determined by $x_{3}$. Thus polynomials $D$ and $D^{1}$ have the same degree $k$ and the same finite number of zeros in $\{z \in \mathbb{C}: |z| >1 \}$, counting multiplicity. Hence there exists $t \in \mathbb{C}, t\neq 0$ where 
\begin{equation} \label{D}
D^{1}=tD \quad on \  \overline{\mathbb{D}}. 
\end{equation}
By equality (\ref{eql13}), for $\lambda \in \overline{\mathbb{D}}$
\begin{equation*}
x_{3}= \dfrac{D^{\sim n}}{D}= \dfrac{D^{1\sim n}}{D^{1}}= \dfrac{\overline{t}D^{\sim n}}{tD}
\end{equation*}
Thus $t=\overline{t}$, and so, $t \in \mathbb{R}\backslash \{0\}$. By the equalities (\ref{eql11}) and \eqref{D}
\begin{equation*}
x_{1}= \dfrac{E_{1}}{D}= \dfrac{E_{1}^{1}}{D^{1}}=\dfrac{E_{1}^{1}}{tD}, \qquad \text{on} \ \mathbb{ \overline{ D}}.
\end{equation*}
This implies that $E_{1}^{1}=tE_{1}$.
By the equalities (\ref{eql12}) and \eqref{D}
\begin{equation*}
x_{2}= \dfrac{E_{2}}{D}= \dfrac{E_{2}^{1}}{D^{1}}=\dfrac{E_{2}^{1}}{tD}, \qquad \text{on} \ \mathbb{ \overline{ D}}. 
\end{equation*}
Thus $E_{2}^{1}=tE_{2}$.
\end{proof}
\begin{remark} \label{dimE} {\rm 
For a fixed polynomial $D$ of degree $n$, the set of polynomials $E_1$ satisfying the conditions of Theorem \ref{min} is a subset of a real vector space of dimension $2n+2$.
Hence the set of rational $\mathbb{ \overline{E }}$-inner functions  of degree $n$ with $x_{3}=\frac{D^{\sim n}}{D}$
 is a subset of a ($2n +2$)-dimensional real space of rational functions.
}
\end{remark}

\begin{lemma} \label{lemma55}
	Let $$x=(x_{1}, x_{2},x_{3})=\bigg(\frac{E_{1}}{D}, \frac{E_{2}}{D},\frac{D^{\sim n}}{D}\bigg)$$ be a rational $\mathbb{ \overline{E }}$-inner function. Then, for $\lambda \in \mathbb{T}$, 
	$$|E_{1}(\lambda)|= |E_{2}(\lambda)|, \quad  \text{and so} \quad
	 |x_{1}(\lambda)|= |x_{2}(\lambda)|.$$ 
\end{lemma}
\begin{proof}
By Theorem \ref{min} (7), for all $\lambda \in \mathbb{T}$, 
\begin{eqnarray*}
	E_{1}(\lambda)&=&E_{2}^{\sim n}(\lambda)= \lambda^{n} \overline{E_{2}(1/\overline{\lambda})}.  	
\end{eqnarray*}
Thus, since $\lambda \overline{\lambda}=1$,
\begin{eqnarray*} 
	|E_{1}(\lambda)| &=& |\lambda^{n} \overline{E_{2}(1/\overline{\lambda})}| \\ 
	&=& |E_{2}(1/\overline{\lambda})| = |E_{2}(\lambda)| .
\end{eqnarray*}
Therefore, for all $\lambda \in \mathbb{T}$, $$|x_{1}(\lambda)|= |x_{2}(\lambda)|.$$
\end{proof}
\begin{example}
		{\em Let $x=(x_{1},x_{2}, x_{3})$ be a rational $\mathbb{ \overline{ E}}$-inner function such that $x_{3}(\lambda)= \lambda$. Clearly, 
\begin{equation*}
			\deg(x)=\deg(x_{3})=1.
			\end{equation*}
By Theorem \ref{min}, there exist polynomials $E_{1}, E_{2}, D$ such that
			\begin{equation*}
		D(\lambda)=1, \qquad 	\deg(E_{1}) \leq 1, \qquad \deg(E_{2}) \leq 1,
			\end{equation*}
\[
E_{1}(\lambda)= E_{2}^{\sim n}(\lambda), \qquad  
  |E_{i}(\lambda)|\leq |D(\lambda)| =1, \ \ i=1,2, \ \mbox{for all } \ \lambda \in \mathbb{\overline {D}},\]
and
\[
x= \left( \frac{E_{1}}{D}, \frac{E_{2}}{D},\frac{D^{\sim 1}}{D}\right).
\]
Therefore the function 
			\begin{equation*}
			x(\lambda)= \big(a_{1}+a_{2}\lambda, \overline{a_{2}} + \overline{a_{1}}\lambda, \lambda \big)
			\end{equation*}
is rational $\mathbb{ \overline{ E}}$-inner for $a_{1}, a_{2} \in \mathbb{ \overline{ D}}$ such that 
$$|a_{1}+a_{2}\lambda| \leq 1 \qquad \text{and} \qquad |\overline{a_{2}} + \overline{a_{1}}\lambda| \leq 1 \qquad \text{for all} \ \lambda \in \overline{\mathbb{D}}.$$
In particular, one can choose $a_{1}=1$ and $a_{2}=0$ to get the rational $\mathbb{ \overline{ E}}$-inner function
			\begin{equation*}
			x(\lambda)= (1,\lambda, \lambda).
			\end{equation*}}
	\end{example}
\begin{example} \textbf{\em $\mathbb{ \overline{E }}$-inner functions} \label{ex} \\
\rm Suppose that $\mathbb{B}_{2 \times 2}$ = $\{A \in \mathbb{C}^{2 \times 2} : \|A\| <1\}$\index{$\mathbb{B}_{2 \times 2}$}. Let us construct an analytic map from the open unit disc $\mathbb{D}$ to $\mathbb{B}_{2 \times 2}$. Consider nonconstant inner functions $\varphi, \psi \in H^{\infty}(\mathbb{D})$ and the diagonal matrix 
\begin{equation*}
h(\lambda)=\begin{bmatrix}\varphi(\lambda) &0\\0& \psi(\lambda) \end{bmatrix} \quad \quad \text{for }\; \lambda \in \mathbb{D}.
\end{equation*}
Note  $\|h(\lambda)\| = \max \{|\varphi(\lambda)|,|\psi(\lambda)| \} < 1$ \; for $\lambda \in \mathbb{D}$ and $h:\mathbb{D} \rightarrow \mathbb{B}_{2 \times 2}$ is analytic. \\
By Theorem \ref{tetbnd}, for all $\lambda \in \mathbb{ D}$,
\begin{equation*}
\big(\varphi(\lambda), \psi(\lambda), \det h(\lambda)\big) \in \mathbb{E },
\end{equation*}
and $\varphi(\lambda)\psi(\lambda)= \det h(\lambda)$. Recall that such points are called triangular points of $\mathbb{E}$. 
 However, we are seeking more interesting and general examples. To get such examples we make use of the singular value decomposition. \newline

Let $U,V$ be $2 \times 2$ \index{unitary matrices} unitary matrices. Then $h_{1}: \mathbb{D} \rightarrow \mathbb{C}^{2 \times 2}$ defined by
\begin{equation*}
h_{1}=UhV 
\end{equation*}
maps $\mathbb{D}$ to $\mathbb{B}_{2 \times 2}$. 
For example, if
 $$U=\begin{bmatrix}\dfrac{1}{\sqrt {2}} &\dfrac{1}{\sqrt {2}}\\-\dfrac{1}{\sqrt {2}}& \dfrac{1}{\sqrt {2}} \end{bmatrix} \quad \text{and} \quad V=I=\begin{bmatrix}1 &0\\0& 1 \end{bmatrix} $$ then $U$ is unitary and we obtain  

\begin{eqnarray*}
h_{1}(\lambda)&=&Uh(\lambda)I\\
&=& \dfrac{1}{\sqrt {2}} \begin{bmatrix}1  &1\\-1& 1 \end{bmatrix} \begin{bmatrix}\varphi(\lambda) &0\\0& \psi(\lambda) \end{bmatrix}  \begin{bmatrix}1 &0\\0& 1 \end{bmatrix} \\
&=&\dfrac{1}{\sqrt{2}}\begin{bmatrix}\varphi(\lambda) &\psi(\lambda) \\-\varphi(\lambda)& \psi(\lambda) \end{bmatrix}.
\end{eqnarray*}
Define $x: \mathbb{D} \rightarrow \mathbb{\overline{E}}$ by $x=\pi \circ h_{1}$, where $\pi$ is given in Definition \ref{pi12}. Then, for $\lambda \in \mathbb{D},$
\begin{eqnarray*} 
x(\lambda)&=& \pi\big(h_{1}(\lambda)\big) \\
&=& \pi \bigg(\dfrac{1}{\sqrt{2}}\begin{bmatrix}\varphi(\lambda) &\psi(\lambda) \\-\varphi(\lambda)& \psi(\lambda) \end{bmatrix}\bigg)\\
&=& \bigg(\dfrac{\varphi (\lambda)}{\sqrt{2}},\dfrac{\psi (\lambda)}{\sqrt{2}}
, \varphi(\lambda)\psi(\lambda)\bigg).
\end{eqnarray*}
Note that this $x(\lambda)$ is not a triangular point unless either $\varphi (\lambda) =0$ or $\psi (\lambda)
=0.$ 
\vspace{3mm}

Let us show that this function $x$ is 
 $\mathbb{ \overline{E }}$-inner.  By Theorem \ref{imp}(1), for $\lambda \in \mathbb{ T}$, since $\varphi , \psi$ are inner functions,
\begin{eqnarray*} 
 \overline{x_{2}(\lambda)} x_{3}(\lambda)&=& \overline{\bigg(\dfrac{\psi (\lambda)}{\sqrt{2}} \bigg) }\varphi(\lambda)\psi(\lambda) 
= \dfrac{  \varphi(\lambda) }{\sqrt{2}}  \overline{\psi (\lambda)} \psi(\lambda) \\ 
&=&  \dfrac{  \varphi(\lambda) }{\sqrt{2}}  |\psi (\lambda)|^{2} 
=  \dfrac{  \varphi(\lambda) }{\sqrt{2}} = x_{1}(\lambda).
\end{eqnarray*} 
Since $|\psi(\lambda)| < 1 $ for $\lambda \in \mathbb{ D}$, this implies that $\Big|\dfrac{\psi(\lambda)}{\sqrt{2}}\Big| < 1$. Thus $|x_{2}(\lambda)| <1$. Finally, for $\lambda \in \mathbb{T}$, since $\varphi , \psi$ are inner functions,
\begin{eqnarray*}
|x_{3}(\lambda)| &=& \big| \varphi(\lambda)\psi(\lambda)\big| \\
&=& | \varphi(\lambda)| |\psi(\lambda)|= 1 .
\end{eqnarray*}
Therefore $x$ is an $\mathbb{ \overline{E }}$-inner function. 

\end{example}
\begin{remark} {\em
	In the previous example if we choose the functions  $ \varphi$ and $\psi$ to be in the Schur class but not to be inner functions then one can check that we obtain an analytic function $x : \mathbb{D} \rightarrow \mathbb{ \overline{E }}$ which is not an $\mathbb{ \overline{E }}$-inner function.
}
\end{remark}

\begin{proposition} \label{gamtet}
	Let $(s,p)$ be a $\Gamma$-inner function. Then $x=\big( \frac{s}{2}, \frac{s}{2}, p \big)$ is an $\mathbb{ \overline{ E}}$-inner function. 
\end{proposition}
\begin{proof}
	By Lemma \ref{TB3}, for every $\lambda \in \mathbb{D}$, $x(\lambda)=\big( \frac{s}{2}(\lambda), \frac{s}{2}(\lambda), p (\lambda)\big) \in \mathbb{ \overline{ E}}$. It is easy to see that $x$ is in
 the set of analytic functions  ${\rm Hol}(\mathbb{ D}, \mathbb{ \overline{ E}})$ from $\D$ to the tetrablock $\mathbb{ \overline{ E}}$. By Proposition \ref{sp} (3), for almost all $\lambda \in \mathbb{T}$,
	\begin{equation*}
	|p(\lambda)|= 1, \qquad |s(\lambda)| \leq 2 \qquad \text{and} \quad s(\lambda)-\overline{s(\lambda)}p(\lambda)=0.
	\end{equation*}
	Thus, for almost all $\lambda \in \mathbb{T}$, 
	\begin{equation*}
	|p(\lambda)|=1, \qquad \dfrac{s(\lambda)}{2} = \dfrac{\overline{s(\lambda)}}{2}p(\lambda) \qquad \text{and} \qquad \dfrac{|s(\lambda)|}{2} \leq 1.
	\end{equation*} 
	 Hence $x$ is $\mathbb{ \overline{ E}}$-inner.
\end{proof}
See \cite{ALY13} for many examples of $\Gamma$-inner functions.

\subsection{Superficial $\mathbb{\overline{E}}$-inner functions }\label{superficial}

  In this subsection, we study $\mathbb{ \overline{ E}}$-inner functions $x$ such that $x(\lambda) $ lies in the topological boundary $\partial\mathbb{ \overline{ E}}$ of $\mathbb{ \overline{ E}}$ for all $\lambda \in \mathbb{ D}$. 

The topological boundary of $\mathbb{E}$ is denoted by $\partial \mathbb{\overline{E}}$. 

\begin{lemma} \em{\cite{AWY}} \label{topbon} 
	Let $x =(x_{1}, x_{2}, x_{3})\in \mathbb{C}^{3}$. Then $x \in \partial \mathbb{E}$ if and only if
	\begin{equation*}
	|x_{1}-\overline{x}_{2}x_{3}|+|x_{2}-\overline{x}_{1}x_{3}|
 = 1 - |x_{3}|^{2}, \quad |x_{1}| \leq 1 \quad \text{and} \quad |x_{2}| \leq 1.
	\end{equation*}
\end{lemma}

Here we show that,
for any inner function $x_{3}$ and $\beta_{1}, \beta_{2} \in \mathbb{ C}$ such that $|\beta_{1}|+|\beta_{2}| = 1$, the function 
$x= (\beta_{2}+ \overline{\beta}_{1}x_{3}, \beta_{1}+ \overline{\beta}_{2}x_{3},x_{3})$ is $\mathbb{ \overline{ E}}$-inner and has the property that it maps $\mathbb{  D}$ to $\partial\mathbb{ \overline{ E}}$.  Recall the definition of superficial function in the set of analytic functions  ${\rm Hol}(\mathbb{ D}, \Gamma)$ from $\D$ to the symmetrised bidisc $\Gamma$ from \cite{ALY13}.
\begin{definition}
	An analytic function $h: \mathbb{D} \rightarrow \Gamma$ is  {\em superficial}  if $h(\mathbb{D}) \subset \partial \Gamma$ \index{$\partial \Gamma$}.
\end{definition}
One can define a similar notion for functions in   ${\rm Hol}(\mathbb{ D}, \mathbb{ \overline{ E}})$.
 
\begin{definition}
	An analytic function $x: \mathbb{D} \rightarrow \mathbb{\overline{E}}$ is  {\em superficial} if $x(\mathbb{D}) \subset \partial \mathbb{\overline{E}}$.
\end{definition}

We also consider relations between superficial $\Gamma$-inner functions and superficial $\mathbb{ \overline{ E}}$-inner functions. 

\begin{proposition} {\em \cite[ Proposition 8.3]{ALY13}} \label{prop123}
	A $\Gamma$-inner function h is superficial if and only if there is an $\omega \in \mathbb{T}$ and an inner function $p$ such that $h=(\omega p+\overline{\omega},p)$.
\end{proposition}

\begin{proposition} \label{super}
	An analytic function $x:\mathbb{D} \rightarrow \mathbb{\overline{E}}$ such that
\begin{equation*}
	x(\lambda)=(\beta_{1}+\overline{\beta_{2}}x_{3}(\lambda),\beta_{2}+\overline{\beta_{1}}x_{3}(\lambda),x_{3}(\lambda)), \quad \lambda \in \mathbb{D}, 
	\end{equation*}
where $x_{3}$ is an inner function and $|\beta_{1}|+|\beta_{2}|=1$ is $\mathbb{ \overline{ E}}$-inner and superficial.
\end{proposition}
\begin{proof}
By Lemma \ref{topbon}, we need to show that, for $\lambda \in \mathbb{D}$, 
	\begin{equation*}
	x(\lambda)=(\beta_{1}+\overline{\beta_{2}}x_{3}(\lambda),\beta_{2}+\overline{\beta_{1}}x_{3}(\lambda),x_{3}(\lambda)). 
	\end{equation*}
	is in $\partial \mathbb{E}$. Here
	$$x_{1}(\lambda)=\beta_{1}+\overline{\beta_{2}}x_{3}(\lambda) , \quad x_{2}(\lambda)=\beta_{2}+\overline{\beta_{1}}x_{3}(\lambda).$$
Note, for $\lambda \in \mathbb{D}$, 
	\begin{eqnarray} \label{Eqn1}
	|\big(x_{1}-\overline{x}_{2}x_{3}\big)(\lambda)|&=&\bigg|\beta_{1}+\overline{\beta}_{2}x_{3}(\lambda)-\Big(\overline{\beta_{2}+\overline{\beta_{1}}x_{3}(\lambda)}\Big)x_{3}(\lambda)\bigg| \\ \nonumber	 
	&=&\bigg|\beta_{1}\Big(1-|x_{3}(\lambda)|^{2}\Big)\bigg|.
	\end{eqnarray}
We also have
	\begin{eqnarray} \label{Eqn2}
	|(x_{2}-\overline{x}_{1}x_{3})(\lambda)|&=&\bigg|\beta_{2}+\overline{\beta}_{1}x_{3}(\lambda)-\Big(\overline{\beta_{1}+\overline{\beta_{2}}x_{3}(\lambda)}\Big)x_{3}(\lambda)\bigg| \\ \nonumber	 
	&=&\bigg|\beta_{2}\Big(1-|x_{3}(\lambda)|^{2}\Big)\bigg|.
	\end{eqnarray}
Note that by equations \eqref{Eqn1} and \eqref{Eqn2}, for all $\lambda \in \mathbb{D}$,
	\begin{eqnarray*}
		|(x_{1}-\overline{x}_{2}x_{3})(\lambda)|+ |(x_{2}-\overline{x}_{1}x_{3})(\lambda)|&=&\bigg|\beta_{1}\Big(1-|x_{3}(\lambda)|^{2}\Big)\bigg|+\bigg|\beta_{2}\Big(1-|x_{3}(\lambda)|^{2}\Big)\bigg| \\
		&=& \Big(|\beta_{1}|+|\beta_{2}|\Big)\big(1-|x_{3}(\lambda)|^{2}\big)=1-|x_{3}(\lambda)|^{2}.
	\end{eqnarray*}
	
	\noindent By Theorem \ref{tetbnd2} and Lemma \ref{topbon}, for  $\lambda \in \mathbb{ D}$, the point $$x(\lambda)=\big(x_{1}(\lambda),x_{2}(\lambda),x_{3}(\lambda)\big)$$ lies in $\partial \mathbb{E}$.
Let us check that $x$ is $\mathbb{ \overline{ E}}$-inner. Clearly, for almost all $\lambda \in \mathbb{T}$,
\begin{eqnarray*}
\overline{x_{2}(\lambda)}x_{3}(\lambda)&=& \big(\overline{\beta_{2}+ \overline{\beta_{1}}x_{3}(\lambda)}\big)x_{3}(\lambda) \\
&=&  \overline{\beta_{2}}x_{3}(\lambda)+ \beta_{1}\overline{x_{3}(\lambda)}x_{3}(\lambda) \\
&=& \beta_{1} + \overline{\beta_{2}}x_{3}(\lambda) = x_{1}(\lambda).
\end{eqnarray*}	
We also have, for almost all $\lambda \in \mathbb{ T}$, 
\begin{equation*}
|x_{2}(\lambda)|= |\beta_{2} + \overline{\beta_{1}}x_{3}(\lambda)| \leq |\beta_{2}|+ |\beta_{1}x_{3}(\lambda)| = |\beta_{2}|+ |\beta_{1}| =1.
\end{equation*}
Since $x_{3}$ is inner, for almost all $\lambda \in \mathbb{T}$, $|x_{3}(\lambda)|=1$. 
Therefore $x(\lambda) \in b\mathbb{ \overline{ E}}$, for almost all $\lambda \in \mathbb{ T}$, and hence $x$ is $\mathbb{ \overline{ E}}$-inner.		
\end{proof}

\begin{lemma}
	Let $x: \mathbb{D} \rightarrow \mathbb{\overline{E}}$ be such that $x(\lambda\big)=\big(\beta_{1}+\overline{\beta_{2}}x_{3}(\lambda),\beta_{2}+\overline{\beta_{1}}x_{3}(\lambda),x_{3}(\lambda)\big)$, where $x_{3}$ is a non-constant rational inner function and $|\beta_{1}|+|\beta_{2}|=1$. Then \newline $  \Psi_{\omega}(x(\lambda))=k \ \text{for all} \  \lambda \in \mathbb{D}$, where 
	\begin{equation*}
	\omega = \frac{\overline{\beta_{2}}}{|\beta_{2}|}, \quad k=\frac{\beta_{1}}{|\beta_{1}|} \qquad \mbox{on} \ \mathbb{ T}.
	\end{equation*}
\end{lemma}
\begin{proof}
By  definition,
	\begin{eqnarray*}
		\Psi_{\omega}(x(\lambda))&=&\frac{x_{3}(\lambda)\omega-x_{1}(\lambda)}{x_{2}(\lambda)\omega-1} \\
		&=& \frac{x_{3}(\lambda)\omega-(\beta_{1}+\overline{\beta}_{2}x_{3}(\lambda))}{(\beta_{2}+\overline{\beta}_{1}x_{3}(\lambda))\omega-1} \\
		&=& \frac{x_{3}(\lambda)\omega-\beta_{1}-\overline{\beta}_{2}x_{3}(\lambda)}{\beta_{2} \omega+\overline{\beta}_{1} \omega x_{3}(\lambda)-1}.
	\end{eqnarray*}
Thus, for all $\lambda \in \mathbb{D}$,  
	\begin{eqnarray*}
		\Psi_{\omega}(x(\lambda))=k &\Leftrightarrow& \frac{x_{3}(\lambda)\omega-\beta_{1}-\overline{\beta}_{2}x_{3}(\lambda)}{\beta_{2} \omega+\overline{\beta}_{1} \omega x_{3}(\lambda)-1}=k \\
		&\Leftrightarrow&  x_{3}(\lambda)\omega-\beta_{1}-\overline{\beta}_{2}x_{3}(\lambda)=k[\beta_{2} \omega+\overline{\beta}_{1} \omega x_{3}(\lambda)-1] \\
		&\Leftrightarrow& x_{3}(\lambda)[\omega-\overline{\beta}_{2}-k\overline{\beta}_{1} \omega] +[k-\beta_{1}-k\beta_{2}\omega]=0 . 
	\end{eqnarray*}
Since $x_{3}$ is a nonconstant rational inner function, this  statement is equivalent to 
	\begin{equation*} 
	\omega-\overline{\beta}_{2}-k\overline{\beta}_{1} \omega = 0 \quad \text{and} \quad k-\beta_{1}-k\beta_{2}\omega = 0.
	\end{equation*} 
\noindent Multiply both sides of the first equation by 
$\overline{\omega}$ and the second equation by $\overline{k}$. We get
	\begin{equation} \label{eqn15}
	\systeme{
		\overline{\beta_{1}}k+\overline{\beta_{2}}\overline{\omega}=1,
		\beta_{1}\overline{k}+\beta_{2}\omega=1.}
	\end{equation}	
Since $|\beta_{1}|+|\beta_{2}|=1$, it
	 is easy to see that
	$$\omega = \frac{\overline{\beta_{2}}}{|\beta_{2}|} \quad \mbox{and} \quad k=\frac{\beta_{1}}{|\beta_{1}|}$$
	satisfy equation \eqref{eqn15}, and so 
	\begin{equation*}
	\Psi_{\omega}(x(\lambda))=k \qquad \mbox{for all} \ \lambda \in \mathbb{  D}.
	\end{equation*}
\end{proof}
\begin{lemma}
For any inner function $x_{3}: \mathbb{  D} \rightarrow \mathbb{ \overline{ D}}$, there are $x_{1}, x_{2}: \mathbb{  D} \rightarrow \mathbb{ \overline{ D}}$ such that the function $x: \mathbb{D} \rightarrow \mathbb{ \overline{E}}$ defined by $x=(x_{1}, x_{2}, x_{3})$ is a superficial $\mathbb{ \overline{ E}}$-inner function, but $h=(x_{1}+x_{2}, x_{3}): \mathbb{  D} \rightarrow \Gamma$ is not a superficial $\Gamma$-inner function.
\end{lemma}

\begin{proof}
	By Proposition \ref{super}, for any $\beta_{1}, \beta_{2} \in \mathbb{ C}$ such that $|\beta_{1}|+ |\beta_{2}|=1$, the function $x: \mathbb{ D} \rightarrow \mathbb{ \overline{ E}}$ defined by
	\begin{equation*}
	x=\Big(\beta_{1}+\overline{\beta_{2}}x_{3},\beta_{2}+\overline{\beta_{1}}x_{3},x_{3}\Big)
	\end{equation*}
	is a superficial $\mathbb{ \overline{ E}}$-inner function.
	By Proposition \ref{prop123}, $h:\mathbb{D} \rightarrow \Gamma$ is superficial if and only if there exists an $\omega \in \mathbb{T}$ such that  $h=(\omega p+ \overline{\omega},p)$. Note that, for $x_{1}= \beta_{1}+\overline{\beta_{2}}x_{3}$ and $x_{2}=\beta_{2}+\overline{\beta_{1}}x_{3}$, 
	\begin{eqnarray*}
		h(\lambda)=(x_{1}+x_{2}, x_{3})(\lambda)&=&\Big(\beta_{1}+\overline{\beta}_{2}x_{3}(\lambda)+\beta_{2}+\overline{\beta}_{1}x_{3}(\lambda),x_{3}(\lambda) \Big) \\
		&=&\Big((\overline{\beta}_{1}+\overline{\beta}_{2})x_{3}(\lambda)+(\beta_{1}+ \beta_{2}),x_{3}(\lambda) \Big), \quad \lambda \in \mathbb{D}.
	\end{eqnarray*}
	One can see that there are some $\beta_{1}, \beta_{2} \in \mathbb{C}$ with $|\beta_{1}|+|\beta_{2}|=1$, but $\beta_{1}+\beta_{2} \notin \mathbb{T}$. For example, take
	\begin{eqnarray*}
		\beta_{1} = i\frac{1}{2},\quad \beta_{2}=-i\frac{1}{2}.
	\end{eqnarray*}
Then
	$|\beta_{1}|+|\beta_{2}|=\frac{1}{2}+\frac{1}{2}=1,$ but $\beta_{1} + \beta_{2} = i\frac{1}{2} -i\frac{1}{2}=0 \notin \mathbb{T}.$
	 Thus, $h$ is not a superficial $\Gamma$-inner function for $\beta_{1} = \frac{i}{2}$ and $\beta_{2}=\frac{-i}{2}$.
\end{proof}

\section{The construction of rational $\mathbb{ \overline{ E}}$-inner functions} \label{construct}

The formula for a Blaschke product is an explicit representation of a rational inner function in terms of its zeros and one other parameter (a unimodular complex number). In this chapter we aim to find a comparable representation for rational $\mathbb{ \overline{ E}}$-inner functions. The first question is: what is the tetrablock analogue of the zeros of an inner function? We shall show that one satisfactory choice consists of the royal nodes of an $\mathbb{ \overline{ E}}$-inner function  $x$ together with the zeros of $x_{1}$ and $x_{2}$. 
We construct a rational $\mathbb{ \overline{ E}}$-inner function $x$ from its royal nodes and the zeros of $x_{1}$ and $x_{2}$. We show that there exists a $3$-parameter family of rational $\mathbb{ \overline{ E}}$-inner functions with prescribed zero sets of $x_{1}$, $x_{2}$ and prescribed royal nodes. We also prove that a nonconstant rational $\mathbb{ \overline{ E}}$-inner function $x$ of degree $n$ either maps $\mathbb{ \overline{ D}}$ to the royal variety of $\mathbb{\overline{E} }$ or $x(\mathbb{ \overline{ D}})$ meets the royal variety exactly $n$ times.

\subsection{The royal polynomial of an $\mathbb{ \overline{ E}}$-inner function} \label{Sec5.2}

We define the {\em royal variety} \index{royal variety for $\mathbb{\overline{E}}$} for $\mathbb{\overline{E}}$ to be 
\begin{equation*}
\mathcal{R}_{\mathbb{\overline{E}}} = \big\{ (x_{1}, x_{2}, x_{3}) \in \C^3 : x_{1}x_{2}=x_{3} \big\}.
\end{equation*}
It was shown in \cite{You08} that 
$\mathcal{R}_{\mathbb{\overline{E}}} \cap \mathbb{E}$ is the orbit of 
$ \{(0,0,0)\}$ under the group of biholomorphic automorphisms of $\mathbb{E}$.
By Theorem \ref{min}, for a rational $\mathbb{\overline{E}}$-inner function $x=(x_{1}, x_{2}, x_{3})$, there are polynomials $E_{1}, E_{2},D$ such that
$$x_{1}= \dfrac{E_{1}}{D}, \quad x_{2}= \dfrac{E_{2}}{D}, \quad x_{3}= \dfrac{D^{\sim n}}{D}.$$
Thus, for $\lambda \in \mathbb{ \overline{ D}}$,  
\begin{equation*}
\big(x_{3} -x_{1}x_{2} \big)(\lambda)= \bigg[\dfrac{D^{\sim n}}{D} - \dfrac{E_{1}}{D} \dfrac{E_{2}}{D}  \bigg](\lambda).
\end{equation*}
The {\em royal polynomial} \index{royal polynomial of $x$} of the rational $\mathbb{\overline{E}}$-inner function $x$ is defined to be 
\begin{eqnarray*} 
R_{x}(\lambda)&=&D^{2}(\lambda) \bigg[\dfrac{D^{\sim n}}{D}-\dfrac{E_{1}E_{2}}{D^{2}} \bigg](\lambda) \nonumber \\
&=& \big[D^{\sim n}D-E_{1}E_{2}\big](\lambda). 
\end{eqnarray*}
\begin{definition} {\em \cite[Definition 3.4]{ALY15} \label{n-sym}}
	We say a polynomial $f$ is {\em $n$-symmetric}  if $\deg(f) \leq n$ and $f^{\sim n}=f$. 
\end{definition}

\begin{definition} {\em \cite[Definition 3.4]{ALY15} \label{n-zeros}}
	For any $E \subset \mathbb{C}$, the number of zeros of $f$ in $E$, counted with multiplicities, is denoted by $\mathrm{ord}_{E}(f)$ \index{$\mathrm{ord}_{E}(f)$} and $\mathrm{ord}_{0}(f)$ means the same as $\mathrm{ord}_{\{0\} }(f)$. \index{$\mathrm{ord}_{\{0\} }(f)$}
\end{definition}

\begin{proposition} \label{propit}
	Let $x$ be a rational $\mathbb{\overline{E}}$-inner function of degree $n$ and let $R_{x}$ be the royal polynomial of $x$. Then, for $\lambda \in \mathbb{T}$,
	\begin{enumerate}
		\item $\lambda^{-n}R_{x}(\lambda) = |D(\lambda)|^{2}   -  |E_{2}(\lambda)|^{2}$  and
		\item $\lambda^{-n}R_{x}(\lambda) = |D(\lambda)|^{2}   -  |E_{1}(\lambda)|^{2}$.
	\end{enumerate}	
\end{proposition}
\begin{proof}
(1) For $\lambda \in \mathbb{T}$,
\begin{eqnarray} \nonumber \label{eqn22}
\lambda^{-n}R_{x}(\lambda) &=& \lambda^{-n} [ D^{\sim n} D - E_{1}E_{2}] (\lambda) \\ \nonumber
&=& \lambda^{-n} [\lambda^{n}\overline{D(1/\overline{\lambda})} D(\lambda)  - E_{2}^{\sim n}(\lambda) E_{2}(\lambda)] ,\quad  \text{since} \ \ E_{1}(\lambda)=E_{2}^{\sim n}(\lambda) \ \text{on} \ \mathbb{T}  \\ \nonumber 
&=& \lambda^{-n} [\lambda^{n}\overline{D(\lambda)} D(\lambda)  -\lambda^{n} \overline{E_{2}(1/\overline{\lambda})} E_{2}(\lambda)], \quad \text{since} \ \ E_{2}(1/\overline{\lambda}) = E_{2}(\lambda) \ \text{on} \ \mathbb{T} \\ \nonumber
&=& \lambda^{-n} [\lambda^{n}\overline{D(\lambda)} D(\lambda)  -\lambda^{n} \overline{E_{2}(\lambda)} E_{2}(\lambda)]   \\
&=&  |D(\lambda)|^{2}   -  |E_{2}(\lambda)|^{2}.  
\end{eqnarray}
(2) Since $x$ is rational $\mathbb{ \overline{E }}$-inner function, by Lemma \ref{lemma55},
\begin{equation} \label{eqn2}
|E_{1}(\lambda)|= |E_{2}(\lambda)| \quad \text{for} \ \lambda \in \mathbb{T}.
\end{equation}
By equations \eqref{eqn22} and \eqref{eqn2},
\begin{equation*} 
\lambda^{-n}R_{x}(\lambda) = |D(\lambda)|^{2}   -  |E_{1}(\lambda)|^{2} \quad \text{for} \ \lambda \in \mathbb{T}.
\end{equation*}
\end{proof}

\begin{lemma}  \label{E1E2}
	Let $E_{1}$ and $E_{2}$ be two polynomials such that \ $\deg E_{1}, \deg E_{2} \leq n$. Then 
	\begin{equation*}
	E_{1}(\lambda)= E_{2}^{\sim n}(\lambda) \ \text{for all}  \ \lambda \in \T \quad \text{if and only if}  \quad E_{2}(\lambda) = E_{1}^{\sim n}(\lambda) \ \text{for all}  \ \lambda \in \T.
	\end{equation*} 
\end{lemma}
\begin{proof}
	Suppose that $E_{1}(\lambda)= E_{2}^{\sim n}(\lambda)$ for all $\lambda \in \T$. Then by definition,
	$$E_{1}(\lambda)= E_{2}^{\sim n}(\lambda) = \lambda^{n}\overline{E_{2}(1/\overline{ \lambda})}, \qquad \lambda \in \T.$$
	Therefore, for all $\lambda \in \T$, 
	\begin{eqnarray*} 
	E_{1}(\lambda) = \lambda^{n}\overline{E_{2}(1/\overline{ \lambda})} \ \text{for all} \ \lambda \in \T &\Leftrightarrow &   
	(1/\lambda^{n}) E_{1}(\lambda)= \overline{E_{2}(1/\overline{ \lambda})} \ \text{for all} \ \lambda \in \T \\  
	&\Leftrightarrow & (1/\overline{\lambda})^{n}\overline{ E_{1}(\lambda)}= E_{2}(1/\overline{ \lambda}) \ \text{for all} \ \lambda \in \T \\  
	&\Leftrightarrow & \lambda^{n}\overline{ E_{1}(1/\overline{\lambda})} = E_{2}( \lambda) \ \text{for all} \ \lambda \in \T \\
	&\Leftrightarrow &  E_{1}^{\sim n}(\lambda)= E_{2}(\lambda)  \ \text{for all} \ \lambda \in \T. 
	\end{eqnarray*}
The converse is obvious.
\end{proof}

\begin{definition}
	A nonzero polynomial $R$ is {\em $n$-balanced} \index{$n$-balanced} if 
	\begin{enumerate}
		\item $\deg(R) \leq 2n$, 
		\item $R$ is $2n$-symmetric, and 
		\item $\lambda^{-n}R(\lambda) \geq 0$ for all $\lambda \in \mathbb{T}$.
	\end{enumerate}
\end{definition}
For completeness, we shall say that zeros of the zero polynomial have infinite order.

\begin{proposition} \label{symm} 
	Let $x$ be a rational $\mathbb{\overline{E}}$-inner function of degree $n$ and let $R_{x}$ be the royal polynomial of $x$. Then $R_{x}$ is $2n$-symmetric, $\lambda^{-n} R_{x}(\lambda)\geq 0$ for all $\lambda \in \mathbb{T}$, and the zeros of $R_{x}$ on $\mathbb{T}$ have even order or infinite order.
\end{proposition}
\begin{proof} To show that $R_{x}$ is $2n$-symmetric we have to prove that $R_{x}^{\sim {2n}}(\lambda)=R_{x}(\lambda)$, for $\lambda \in \T$.
Recall that 
  \begin{equation*}
  R_{x}(\lambda) = D^{\sim n}(\lambda)D(\lambda)-E_{1}(\lambda)E_{2}(\lambda), \quad \text{where} \ x= \bigg( \frac{E_{1}}{D}, \frac{E_{2}}{D}, \frac{D^{\sim n}}{D}   \bigg),
  \end{equation*}
By Theorem \ref{min} (7) and Lemma \ref{E1E2}, for $ \lambda \in \T$,
  \begin{equation*}
  E_{1}(\lambda)= E_{2}^{\sim n}(\lambda) = \lambda^{n} \overline{E_{2}(1/\overline{\lambda})}, \qquad  E_{2}(\lambda)= E_{1}^{\sim n}(\lambda) = \lambda^{n} \overline{E_{1}(1/\overline{\lambda})}.
  \end{equation*}
Now 
  \begin{eqnarray*}
  	R_{x}^{\sim {2n}}(\lambda)
&=&\lambda^{2n} \overline{R_{x}(1/ \overline{\lambda})} \\
&=& \lambda^{2n}\overline{ \big[D^{\sim n}(1/ \overline{\lambda})D(1/ \overline{\lambda})-E_{1}(1/\overline{\lambda}) E_{2}(1/ \overline{\lambda})\big]} \\
&=& \lambda^{n} \overline{D^{\sim n}(1/ \overline{\lambda})} \  \lambda^{n}\overline{D(1 / \overline{\lambda})} - \lambda^{n} \overline{E_{1}(1/\overline{\lambda})}\ \lambda^{n} \overline{E_{2}(1/\overline{\lambda})} \\ 
  			&=& D(\lambda) D^{\sim n}(\lambda)- E_{2}(\lambda) E_{1}(\lambda)= R_{x}(\lambda).
  \end{eqnarray*}
Hence $R_{x}$ is $2n$-symmetric.

Clearly, if $x(\overline{\mathbb{D}}) \subseteq \mathcal{R}_{\mathbb{\overline{E}}}$, the royal polynomial $R_{x}$ is identically zero. Thus the zeros of $R_{x}$ on $\mathbb{T}$ have infinite order. 

On the other hand, if $x(\mathbb{ \overline{D}}) \nsubseteq \mathcal{R}_{\overline {\mathbb{E}}} $, by Proposition \ref{propit}, for $\lambda \in \mathbb{T}$,  
 \begin{equation} \label{eqn8888}
 \lambda^{-n}R_{x}(\lambda) = |D(\lambda)|^{2}   -  |E_{2}(\lambda)|^{2}.
 \end{equation}
  By Theorem \ref{min} (6),
  \begin{equation} \label{eqn8889}
  |D|^{2}   -  |E_{2}|^{2} \geq 0 \ \  \text{on} \  \mathbb{T}.
  \end{equation}
  By equations \eqref{eqn8888} and \eqref{eqn8889}, $\lambda^{-n}R_{x}(\lambda) \geq 0$ on $\mathbb{ T}$. By the Fej\'er-Riesz theorem \cite[Section 53]{RN}, there exists an analytic polynomial $P(\lambda) = \sum_{i=0}^{n}b_{i}\lambda^i$ of degree $n$ such that $P$ is outer and 
  \begin{equation*}
  \lambda^{-n}R_{x}(\lambda)= |P(\lambda)|^{2} \qquad \text{for all} \ \lambda \in \mathbb{T}.
  \end{equation*} 
Hence if $\sigma \in \mathbb{T}$ is a zero of $R_{x}$, then $\sigma$ is a zero of even order. Therefore in the case $x(\overline{\mathbb{D}}) \nsubseteq \mathcal{R}_{\overline {\mathbb{E}}} $, the zeros of $R_{x}$ that lie in $\mathbb{T}$ have even order.
\end{proof}

\begin{lemma} \label{nbalanced}
	Let $x$ be a rational $\;\overline{\mathbb{E}}$-inner function of degree $n$. Then the royal polynomial $R_{x}$ of $x$ is either $n$-balanced or identically zero. 
\end{lemma}
	\begin{proof}	If $x(\overline{\mathbb{D}}) \subseteq \mathcal{R}_{\overline{ \mathbb{E}}}$ then, by the definition of the royal variety,
		\begin{equation*}
		x_{1}(\lambda) x_{2}(\lambda) = x_{3}(\lambda) \quad \quad \text{for all} \ \lambda \in \overline{ \mathbb{D}}.
		\end{equation*}
		Thus $$E_{1}(\lambda)E_{2}(\lambda)= D(\lambda)D^{\sim n}(\lambda)  \quad \quad \text{for all} \ \lambda \in \overline{\mathbb{ D}}.$$ Therefore the royal polynomial $R_{x}$ is identically zero. 
		
		If $x(\overline{\mathbb{D}}) \nsubseteq \mathcal{R}_{\overline{ \mathbb{E}}}$ then, by Proposition \ref{symm},  the royal polynomial $R_{x}$ of $x$ is $2n$-symmetric and $\lambda^{-n}R_{x}(\lambda) \geq 0$ for all $\lambda \in \mathbb{T}$. Clearly, $\deg(R_{x}) \leq 2n$.
	 Hence $R_{x}$ is $n$-balanced. 
	\end{proof}

\subsection{Rational $\mathbb{\overline{E}}$-inner functions of type $(n,k)$}\label{(n,k)type}

\begin{definition} \label{def11}
	Let $x=(x_{1}, x_{2}, x_{3})$ be a rational $\;\mathbb{\overline{E}}$-inner function such that \newline $x(\overline{\mathbb{D}}) \nsubseteq \mathcal{R}_{\overline{ \mathbb{E}}} $. Let $R_{x}$ be the royal polynomial of $x$. Let $\sigma$ be a zero of $R_{x}$ of order $\ell$. We define the multiplicity $\# \sigma$ \index{$\# \sigma$ } of $\sigma$ (as a royal node of $x$)  by 
	$$
	\# \sigma =
	\begin{cases}
	\ell \quad \quad \quad \ \ \  \text{if} \ \sigma \in \mathbb{D},\\
	\frac{1}{2} \ell \quad \quad \quad \text{if} \ \sigma \in \mathbb{T} .
	\end{cases}
	$$
We define the type of $x$ to be the ordered pair $(n,k)$, where $n$ is the sum of the multiplicities of the royal nodes of $x$ that lie in $\mathbb{\overline{D}}$, and $k$ is the sum of the multiplicities of the royal nodes of $x$ that lie in $\mathbb{T}$. 
\end{definition}

\begin{definition} \label{R}
	Let $\mathcal{R}^{n,k}$ \index{$\mathcal{R}^{n,k}$} denote the collection of rational $\mathbb{\overline{E}}$-inner functions of type $(n,k)$.
\end{definition}

\begin{remark} {\em \cite[Equations (3.2) and (3.3)]{ALY15} } \label{degr}
{\em	  For any $m$-symmetric polynomial $f$, the following two relations hold
\begin{enumerate}
\item \begin{equation*}
			\deg(f)=m-\mathrm{ord}_{0}(f).
			\end{equation*} 
\item Since $f$ is $m$-symmetric, if $\alpha \in \mathbb{ D}\backslash \{0\}$ is a zero of $f$, then $\frac{1}{\overline{\alpha}}$ is also a zero of $f$. Thus
\begin{equation*}
			\mathrm{ord}_{0}(f) +2\mathrm{ord}_{ \mathbb{D} \backslash \{0 \} }(f)+\mathrm{ord}_{\mathbb{T}}(f)= \deg(f).
			\end{equation*}			
\end{enumerate} 
}
\end{remark}

\begin{theorem} \label{Thm23}
Let $x \in \mathcal{R}^{n,k}$ be nonconstant. Then the degree of $x$ is equal to $n$. 
\end{theorem}
\begin{proof}	Let $R_{x}$ be the royal polynomial of $x$. By assumption $x \in \mathcal{R}^{n,k}$ and is nonconstant. Hence $n \geq 1$ and $x(\mathbb{\overline{\mathbb{D}}}) \nsubseteq \mathcal{R} _{\overline {\mathbb{E}}}$. Thus $R_{x}$ is not identically zero.  By Proposition \ref{symm}, $R_{x}$ is $2\deg(x)$-symmetric. By Remark \ref{degr} $(1)$ and (2), it follows that 
$$\deg(R_{x})=2\deg(x)-\mathrm{ord}_{0}(R_{x})$$
and
$$\mathrm{ord}_{0}(R_{x}) +2\mathrm{ord}_{ \mathbb{D} \backslash \{ 0\} }(R_{x})+\mathrm{ord}_{\mathbb{T}}(R_{x})= \deg(R_{x}).$$
Substitute the first equation in the second equation,
$$\mathrm{ord}_{0}(R_{x}) +2\mathrm{ord}_{ \mathbb{D} \backslash \{0\} }(R_{x})+\mathrm{ord}_{\mathbb{T}}(R_{x}) = 2\deg(x)-\mathrm{ord}_{0}(R_{x}),$$
which implies that 
	 $$2\mathrm{ord}_{0}(R_{x}) +2\mathrm{ord}_{ \mathbb{D} \backslash \{0 \} }(R_{x})+\mathrm{ord}_{\mathbb{T}}(R_{x})=2 \deg(x).$$
	 Therefore, by Definition \ref{R},
	 $$n=\mathrm{ord}_{0}(R_{x}) +\mathrm{ord}_{ \mathbb{D} \backslash \{ 0\} }(R_{x})+ \frac{1}{2}\mathrm{ord}_{\mathbb{T}}(R_{x}) =\deg(x).$$
\end{proof}

\begin{theorem} \label{times}
Let $x$ be a nonconstant rational $\mathbb{\overline{E}}$-inner function. Then either \newline $x(\overline{\mathbb{D}}) \subseteq \mathcal{R}_{\mathbb{\overline{E}}}$ or $x(\overline{\mathbb{D}})$ meets $\mathcal{R}_{\mathbb{\overline{E}}}$ exactly $\deg (x)$ times.
\end{theorem}
\begin{proof}
	Suppose that $x$ is a nonconstant rational $\mathbb{\overline{ E}}$-inner function. Then either, $x(\mathbb{ \overline{D }} )\subseteq \mathcal{R}_{\mathbb{\overline{E}}}$ and the royal polynomial $R_{x}$ of $x$ is identically zero, or by Theorem \ref{Thm23}, $x(\mathbb{ \overline{D }})$ meets $\mathcal{R}_{\mathbb{ \overline{E }}}$ exactly $\deg(x)$ times.
\end{proof}

\begin{lemma} {\em\cite[Lemma 4.4]{ALY15}} \label{Lem23} Let $R$ be a nonzero polynomial and let  $n$ be a positive integer.
	For $ \sigma \in \overline{\mathbb{D}}$, let the polynomial $Q_{\sigma}$ be defined by the formula
	$$Q_{\sigma}(\lambda)=(\lambda - \sigma)(1-\overline{\sigma} \lambda).$$ 
	 The polynomial $R$ is $n$-balanced if and only if there are points $\sigma_{1}, \sigma_{2},...,\sigma_{n} \in \overline{\mathbb{D}}$ and $t_{+} >0$ such that
	\begin{equation*}
	R(\lambda)=t_{+} \prod_{j=1}^{n}Q_{\sigma_{j}}(\lambda), \quad \lambda \in \mathbb{C}.
	\end{equation*}
\end{lemma}

\begin{proposition} \label{mult} Let  $x$ be  a rational $\overline{\mathbb{E}}$-inner function. Suppose the royal nodes of $x$ are $\sigma_{1},..., \sigma_{n}$,
	 with repetition according to the multiplicity of the royal nodes as described in Definition {\em\ref{def11}}. Then the royal polynomial $R_{x}$ of $x$, up to a positive multiple, is 
	\begin{equation} \label{eq258}
	R_{x}(\lambda)= \prod_{j=1}^{n} Q_{\sigma_{j}}(\lambda).
	\end{equation}
\end{proposition}

\begin{proof}
	By Lemma \ref{nbalanced}, $R_{x}$ is $n$-balanced. This implies that, by Lemma $\ref{Lem23}$, there exists $t_{+} >0$ and $\eta_{1}, \dots , \eta_{n} \in \mathbb{ \overline{D }}$ such that 
	\begin{equation*}
	R_{x}(\lambda) = t_{+} \prod_{j=1}^{n} Q_{\eta_{j}}(\lambda).
	\end{equation*} 
	By Definition \ref{def11},  the list $\eta_{1}, \dots , \eta_{n}$ coincides, up to a permutation, with the list $\sigma_{1}, \dots , \sigma_{n}$. Therefore $R_{x}$ is given, up to a positive multiple, by equation \eqref{eq258}.
\end{proof}
	
Before we proceed to the next theorem on the construction of a tetra-inner function $x$ from the zeros of $x_{1}$ and $x_{2}$ and royal nodes of $x$, let us prove the following elementary lemma. 	
\begin{lemma} \label{lem22}
	 Let $E_{1}$ and $E_{2}$ be polynomials of degree at most $n$ such that \newline $E_{1}(\lambda)= E_{2}^{\sim n}(\lambda)$, for $\lambda \in \mathbb{ \overline{D}}$. Let  $\alpha_{1}^{1},...,\alpha_{k_{1}}^{1}$ be the zeros of $E_{1}$ in $\mathbb{ \overline{D}}$, and let $\alpha_{1}^{2},...,\alpha_{k_{2}}^{2}$ be the zeros of $E_{2}$ in $\mathbb{ \overline{D}}$, where $k_{1}+k_{2}=n$. Then 
	\begin{equation*}
	E_{1}(\lambda)= t\prod_{j=1}^{k_{1}}(\lambda - \alpha_{j}^{1}) \prod_{j=1}^{k_{2}}(1- \overline{\alpha}_{j}^{2}\lambda),
	\end{equation*}
where $t \in \mathbb{C} \backslash \{0 \}$.
\end{lemma}
\begin{proof}	Since $\alpha_{1}^{1},...,\alpha_{k_{1}}^{1} \in \overline{\mathbb{D}}$ and $\alpha_{1}^{2},...,\alpha_{k_{2}}^{2} \in\overline{\mathbb{D}}$, where $k_{1}+k_{2}=n$, are the zeros of $E_{1}$ and $E_{2}$ respectively, we have 
\begin{equation} \label{eqnn1}
E_{1}(\lambda) = (\lambda - \alpha_{1}^{1})...(\lambda - \alpha_{k_{1}}^{1} ) p_{1}(\lambda)
\end{equation}
and
\begin{equation*} 
E_{2}(\lambda) = (\lambda - \alpha_{1}^{2})...(\lambda - \alpha_{k_{2}}^{2} ) p_{2}(\lambda).
\end{equation*}
where the polynomials $p_{1}(\lambda)$ and $p_{2}(\lambda)$ do not vanish in $\mathbb{\overline{ D}}$. \newline
Since $E_{1}(\lambda)=\lambda^{n} \overline{E_{2}(1/\overline{\lambda)}}$ on $\mathbb{ \overline{D }}$, we have
\begin{eqnarray} \nonumber \label{eqnn3}
E_{1}(\lambda)&=&  \lambda^{n} \overline{(1/\overline{\lambda}- \alpha_{1}^{2})}...\overline{(1/\overline{\lambda}- \alpha_{k_{2}}^{2})}. \overline{p_{2}(1/\overline{\lambda})} \\ 
&=&  \lambda^{n-k_{2}} (1- \overline{\alpha}_{1}^{2}\lambda)...(1- \overline{\alpha}_{k_{2}}^{2}\lambda). \overline{p_{2}(1/\overline{\lambda})}.
\end{eqnarray}
Since $\deg E_{1} \leq n$ and $k_{1}+k_{2}=n$, equations \eqref{eqnn1} and \eqref{eqnn3} implies that $E_{1}$ can be written in the form
\begin{eqnarray*} 
E_{1}(\lambda)&=&t_{1} (\lambda - \alpha_{1}^{1})...(\lambda - \alpha_{k_{1}}^{1} ) \ \  (1- \overline{\alpha}_{1}^{2}\lambda)...(1- \overline{\alpha}_{k_{2}}^{2}\lambda) \\ 
&=&t_{1} \prod_{j=1}^{k_{1}}(\lambda - \alpha_{j}^{1}) \prod_{j=1}^{k_{2}}(1- \overline{\alpha}_{j}^{2}\lambda), \quad \quad \lambda \in \mathbb{\overline{D}},  
\end{eqnarray*}
for some $t_{1} \in \mathbb{C}\backslash \{0\}$, and
\begin{equation*}
E_{2}(\lambda)= t_{2} \prod_{j=1}^{k_{2}}(\lambda - \alpha_{j}^{2}) \prod_{j=1}^{k_{1}}(1- \overline{\alpha}_{j}^{1}\lambda) \quad \lambda \in \mathbb{\overline{D}},
\end{equation*}
for some $t_{2} \in \mathbb{C}\backslash\{0\}$.
Since $E_{2}(\lambda) = \lambda^{n}\overline{E_{1}(1/\overline{\lambda})}$,
\begin{eqnarray} \nonumber
 \lambda^{n}\overline{E_{1}(1/\overline{\lambda})} &=&  \lambda^{n} \bigg( \overline{ t_1 \prod_{j=1}^{k_{1}}(1/\overline{\lambda} - \alpha_{j}^{1}) \prod_{j=1}^{k_{2}}(1- \overline{\alpha}_{j}^{2}1/\overline{\lambda})} \bigg)    \\ \nonumber
&=& \lambda^{n} \overline{t_{1}} \prod_{j=1}^{k_{1}}(1/\lambda - \overline{\alpha}_{j}^{1}) \prod_{j=1}^{k_{2}}(1- \alpha_{j}^{2} 1/\lambda)   \\ \nonumber
&=&  \overline{t_{1}} \prod_{j=1}^{k_{1}}(1 - \overline{\alpha}_{j}^{1}\lambda) \prod_{j=1}^{k_{2}}(\lambda- \alpha_{j}^{2} )   \\ \nonumber
&=&   E_{2}(\lambda)= t_{2} \prod_{j=1}^{k_{2}}(\lambda - \alpha_{j}^{2}) \prod_{j=1}^{k_{1}}(1- \overline{\alpha}_{j}^{1}\lambda), \qquad \lambda \in \mathbb{\overline{D}},
\end{eqnarray}
and so $t_{2} = \overline{t}_{1}$.
\end{proof}

\begin{remark} {\em For the polynomials $E_{1}$ and $E_{2}$ from Lemma {\rm \ref{lem22}}, if $\alpha  \in \mathbb{D}\backslash\{0\}$ and $\alpha$ is a zero of $E_{1}$ then $\frac{1}{\overline{\alpha}}$ is a zero of $E_{2}$.
}
\end{remark}

	\begin{theorem} \label{cons}
	Suppose that $\alpha_{1}^{1},...,\alpha_{k_{1}}^{1} \in \overline{\mathbb{D}}$ and $\alpha_{1}^{2},...,\alpha_{k_{2}}^{2} \in\overline{\mathbb{D}}$, where $k_{1}+k_{2}=n$. Suppose that $\sigma_{1},..., \sigma_{n} \in \overline{\mathbb{D}}$ are distinct from the points of the set $\{\alpha_{j}^{i}, j=1,...,k_{i} , i=1,2\} \cap \mathbb{T}$. Then there exists a rational $\overline{\mathbb{E}}$-inner function $x=(x_{1},x_{2},x_{3}) : \mathbb{D} \rightarrow \overline{\mathbb{E}}$ such that
	\begin{enumerate}
		\item the zeros of $x_{1}$ in $\overline{ \mathbb{D}}$, repeated according to multiplicity, are $\alpha_{1}^{1},...,\alpha_{k_{1}}^{1}$;
		\item the zeros of $x_{2}$ in $\overline{ \mathbb{D}}$, repeated according to multiplicity, are $\alpha_{1}^{2},...,\alpha_{k_{2}}^{2}$;
		\item the royal nodes of $x$ are $\sigma_{1},..., \sigma_{n} \in \overline{\mathbb{D}}$, with repetition according to the multiplicity of the nodes.
	\end{enumerate}
Such a function $x$ can be constructed as follows. Let $t_{+}>0$ and let $t \in \mathbb{C} \backslash \{0\}$. Let $R$ be defined by
\begin{equation*}
R(\lambda)= t_{+} \prod_{j=1}^{n}(\lambda-\sigma_{j})(1-\overline{\sigma}_{j}\lambda), \quad \text{and} 
\end{equation*} 
let $E_{1}$ be defined by
\begin{equation*}
E_{1}(\lambda)=t \prod_{j=1}^{k_{1}}(\lambda-\alpha_{j}^{1}) \prod_{j=1}^{k_{2}}(1-\overline{\alpha}_{j}^{2}\lambda).
\end{equation*}

Then the following statements hold.
\begin{enumerate}  \label{eqn23}
	\item[(i)] There exists an outer polynomial $D$ of degree at most $n$ such that
	\begin{equation} \label{eqn26}
	\lambda^{-n}R(\lambda)+|E_{1}(\lambda)|^{2}=|D(\lambda)|^{2}
	\end{equation}
	for all $\lambda \in \mathbb{T}$.
	
	\item[(ii)] The function $x$ defined by 
	\begin{equation*}
	x=(x_{1},x_{2},x_{3})=\bigg( \frac{E_{1}}{D}, \frac{E_{1}^{\sim n}}{D},\frac{D ^{\sim n}}{D}\bigg)
	\end{equation*}
	is a rational $\ \overline{\mathbb{E}}$-inner function such that the degree of $x$ is equal to $n$ and conditions {\em(1), (2)} and {\em(3)} hold. The royal polynomial of $x$ is $R$.
\end{enumerate}
\end{theorem}
\begin{proof}
	(i) By Lemma \ref{Lem23}, $R$ is $n$-balanced, and so $\lambda^{-n}R(\lambda) \geq 0$ for all $\lambda \in \mathbb{T}$. Therefore
	\begin{equation*} 
		\lambda^{-n}R(\lambda)+ |E_{1}(\lambda)|^{2} \geq 0 \quad \text{for all} \ \lambda \in \mathbb{T}.
	\end{equation*}
	By the Fej\'er-Riesz theorem, there exists an outer polynomial $D$ of degree at most $n$ such that
	\begin{equation} \label{deqn}
	\lambda^{-n}R(\lambda)+ |E_{1}(\lambda)|^{2}=|D(\lambda)|^{2} \quad \text{for all} \ \lambda \in \mathbb{T}.
	\end{equation} 
	
	(ii)  By hypothesis \newline
	 $$\{ \sigma_{j} : 1\leq j \leq n  \} \cap \bigg( \{ \alpha_{j}^{i} : 1 \leq j \leq k_{i}, i=1,2 \} \cap \mathbb{ T} \bigg)= \emptyset .$$ Thus $\lambda^{-n}R(\lambda)$ and $|E_{1}(\lambda)|^{2}$ are non-negative trigonometric polynomials on $\mathbb{T}$ with no common zero. Therefore
	\begin{equation*}
	\lambda^{-n}R(\lambda) + |E_{1}(\lambda)|^{2} > 0 \quad \text{on} \ \mathbb{T}.
	\end{equation*}
	The outer polynomial $D$ is defined by equality \eqref{deqn} for all $\lambda \in \mathbb{T}$.
Hence $D$ has no zero on $\mathbb{T}$, and so $D$ and $D^{\sim n}$ have no common factor. Thus
	\begin{equation*}
	\deg (x_{3})= \deg \bigg(\frac{D^{\sim n}}{D}\bigg)= \max \{\deg(D), \deg(D^{\sim n})\}=n.
	\end{equation*}
	Since \begin{equation*}
		\lambda^{-n}R(\lambda) \geq 0 \quad \text{for all} \ \lambda \in \mathbb{T},
	\end{equation*}
	\begin{equation*}
	|D(\lambda)|^{2}= 	\lambda^{-n}R(\lambda)+ |E_{1}(\lambda)|^{2} \geq |E_{1}(\lambda)|^{2}
	\end{equation*}
	for all $\lambda \in \mathbb{T}$. It follows that  
	\begin{equation*}
	|D(\lambda)| \geq |E_{1}(\lambda)|, \quad \quad \text{for all} \ \lambda \in \mathbb{T}.
	\end{equation*}
	Since $D(\lambda) \neq 0$ on $\mathbb{\overline{D}}$, the function $\frac{E_{1}}{D}$ is analytic in a neighbourhood of $\mathbb{\overline{D}}$. By the Maximum Modulus Principle, we have
	\begin{equation*}
	\frac{|E_{1}(\lambda)|}{|D(\lambda)|} \leq 1 \quad \text{for all} \ \lambda \in \mathbb{\overline{D}}.
	\end{equation*}
	
	\noindent Therefore, by the converse of Theorem \ref{min}, since conditions (1), (2), (6) and (7)  are satisfied, the function 
	\begin{equation*}
	x(\lambda) = \bigg(  \frac{E_{1}(\lambda)}{D(\lambda)}, \frac{E_{1}^{\sim n}(\lambda)}{D(\lambda)} , \frac{D^{\sim n}(\lambda)}{D(\lambda)}  \bigg) \quad \text{for} \ \lambda \in \mathbb{D},
	\end{equation*}
	is a rational $\mathbb{\overline{E}}$-inner function such that $\deg(x)=n$.
	The royal polynomial of $x$ is defined by 
	\begin{equation*}
	R_{x}(\lambda)= D(\lambda) D^{\sim n}(\lambda)-E_{1}(\lambda)E_{2}(\lambda), \quad \lambda \in \mathbb{D},
	\end{equation*}
	where $E_{2}(\lambda)= E_{1}^{\sim n}(\lambda), \lambda \in \mathbb{D}$. 
	By Proposition \ref{propit}, for all $\lambda \in \mathbb{T}$,
	\begin{eqnarray} \nonumber
	\lambda^{-n} R_{x}(\lambda)&=& |D(\lambda)|^{2}-|E_{1}(\lambda)|^{2}.	 
	\end{eqnarray}
	Therefore, by equation \eqref{eqn26},
	\begin{equation*}
   \lambda^{-n} R_{x}(\lambda) = \lambda^{-n} R(\lambda) \quad \quad \text{for all} \ \  \lambda \in \mathbb{T},
	\end{equation*}
	where $E_{2}(\lambda) = E_{1}^{\sim n}(\lambda)$ for $\lambda \in \mathbb{D}$.
	Thus $ R_{x} = R$, that is, the royal polynomial of $x$ is equal to $R$.
 \end{proof}

\begin{remark} {\em (1) By the Fej\'er-Riesz theorem, there exists an outer polynomial $D$
 satisfying equation \eqref{eqn26}, but to find algebraically such $D$ can be difficult for large $n$. \\
(2) The solution $D$ is only identified up to a multiplication by $\overline{\omega} \in \mathbb{ T}$. Thus if we replace $D$ by $\overline{\omega}D$ we obtain a new solution
		\begin{equation*}
		x= \bigg(\omega\frac{E_{1}}{D}, \omega \frac{E_{1}^{\sim n}}{D}, \omega^{2}\frac{D^{\sim n}}{D}\bigg).
		\end{equation*} 
}
\end{remark}
\begin{example} \label{ex1} {\em 
	Let $n=1$, $\alpha_{1}^{2}= \tfrac{1}{2}$ and $\sigma_{1} = 0$. Let us construct a rational $\mathbb{ \overline{ E}}$-inner function $x=(x_{1}, x_{2}, x_{3}): \mathbb{D} \rightarrow \mathbb{ \overline{E }}$ such that $\alpha_{1}^{2}$ is a zero of $x_{2}$ and $\sigma_{1}$ is a royal node of $x$. \newline
	
	As in Theorem \ref{cons}, for $\lambda \in \mathbb{T}$, let
	\begin{eqnarray*}
    R(\lambda) &=& t_{+}\lambda, \qquad \quad t_{+} \  \text{is a positive real number.} \\
    E_{1}(\lambda) &=& t (1-\tfrac{1}{2} \lambda), \qquad \quad t \in \mathbb{C}\backslash \{0\}.
	\end{eqnarray*}
	 
	\noindent The equation \eqref{eqn26} for the polynomial $D$ is the following, for all $\lambda \in \mathbb{ T}$,
	\begin{eqnarray} \label{eqn201} 
	|D(\lambda)|^{2} &=& \lambda^{-1} R(\lambda) +  |E_{1}(\lambda)|^{2} \\ \nonumber
	&=&  \overline{ \lambda} t_{+} \lambda + |t(1-\tfrac{1}{2}\lambda)|^{2}   \\ \nonumber
	&=& t_{+} + \tfrac{5}{4}|t|^{2}- \tfrac{|t|^{2}}{2} \lambda - \tfrac{|t|^{2}}{2} \overline{\lambda}.
	\end{eqnarray}
	Since the degree of $D$ is at most $1$, 
	 $D(\lambda)= a_{1}+ a_{2}\lambda$, where $a_{1}, a_{2} \in \mathbb{C}$ and $\lambda \in \mathbb{ T}$,
	\begin{eqnarray} \label{eqn202}
	D(\lambda) \overline{ D(\lambda)} &=& |a_{1}+ a_{2} \lambda|^{2} \nonumber \\
	&=& (a_{1}+ a_{2} \lambda)(\overline{a}_{1}+ \overline{a}_{2} \overline{\lambda}) = |a_{1}|^{2} + |a_{2}|^{2} + a_{1} \overline{a}_{2} \overline{ \lambda} + \overline{a}_{1}a_{2} \lambda.
	\end{eqnarray}
	Compare equations \eqref{eqn201} and \eqref{eqn202}. We have
	\begin{equation} \label{aa}
	\begin{cases}
	\overline{a}_{1} a_{2} = -\dfrac{|t|^{2}}{2},
	\\
	a_{1} \overline{a}_{2} = -\dfrac{|t|^{2}}{2},
	\\
	|a_{1}|^{2} + |a_{2}|^{2} = t_{+} + \tfrac{5}{4}|t|^{2} .
	\end{cases}
	\end{equation}
	Finally the function $x$ can be written in the form 
	\begin{equation*}
	x = \bigg( \frac{E_{1}}{D}, \frac{E_{1}^{\sim 1}}{D}, \frac{D^{\sim 1}}{D } \bigg),
	\end{equation*}
	where
	\begin{equation*}
	\begin{cases} \vspace{0.3cm}
	x_{1}(\lambda) =\dfrac{E_{1}}{D}(\lambda) = \dfrac{t (1-\frac{1}{2} \lambda)}{a_{1}+ a_{2}\lambda},
	\\ \vspace{0.3cm}
	x_{2}(\lambda) =\dfrac{E_{1}^{\sim 1}}{D}(\lambda) = \dfrac{\overline{t }(\lambda - \frac{1}{2})}{a_{1}+ a_{2}\lambda},
	\\
	x_{3}(\lambda)= \dfrac{D^{\sim 1}}{D }(\lambda) =\dfrac{\overline{ a}_{1}\lambda + \overline{ a}_{2}}{a_{1}+ a_{2}\lambda},
	\end{cases}
	\end{equation*}
where $|a_{2}| < |a_{1}|$ and $a_{1}$, $a_{2}$ are given by solving equations \eqref{aa} as functions of $t_{+}$ and $t$. These formulas give a parametrization of solutions for the above problem.}
\end{example}
For example, for the given $t= \sqrt{2}$ and $t_{+}= \dfrac{7}{4}$, the system 
\eqref{aa}
has  solutions $$a_{1} = -2i, a_{2} = \tfrac{1}{2}i \quad \text{or} \quad a_{1} = -2, a_{2} = \tfrac{1}{2}.$$ 
Hence the functions 
\[
x (\lambda) = \left(\dfrac{\sqrt{2} (1-\tfrac{1}{2} \lambda)}{-2 \omega  (1- \tfrac{1}{4}\lambda)},  \dfrac{\sqrt{2} (\lambda-\tfrac{1}{2})}{-2 \omega  (1- \tfrac{1}{4}\lambda)}, \dfrac{\overline{\omega}(\lambda - \tfrac{1}{4})}{\omega(1-\tfrac{1}{4}\lambda)} \right), \ \ \lambda \in \D,
\] 
where $\omega \in \mathbb{ T}$, 
are rational $\mathbb{ \overline{ E}}$-inner functions such that $\dfrac{1}{2}$ is a zero of $x_{2}$ and $0$ is a royal node of $x$. 

\begin{remark} {\em
	Theorem {\ref{cons}} shows that there exists a $3$-parameter family of rational $\mathbb{ \overline{ E}}$-inner functions with given royal nodes and given zeros of $x_{1}$ and $x_{2}$. It looks at first sight that the construction in Theorem { \ref{cons}} gives us a $4$-parameter family of rational $\mathbb{ \overline{ E}}$-inner functions with the given data. However, the choice of $t_{+}, t, D$ and $\omega$ leads to the same $x$ as the choice $1, t/\sqrt{t_{+}}, D/\sqrt{t_{+}}$ and $\omega$. Theorem {\ref{conv}} tells us that the construction yields all solutions of the problem, and so the family of functions $x$ with the required properties is indeed a $3$-parameter family. }
\end{remark}

\begin{theorem} \label{conv}
	Let $x= (x_{1}, x_{2}, x_{3})$ be a rational $\mathbb{ \overline{E }}$-inner function of degree $n$ such that
	\begin{enumerate} 
		\item the zeros of $x_{1}$ are $\alpha_{1}^{1}, \dots \alpha_{k_{1}}^{1} \in \mathbb{ \overline{D}}$, repeated according to multiplicity,
		\item the zeros of $x_{2}$ are $\alpha_{1}^{2}, \dots \alpha_{k_{2}}^{2} \in \mathbb{ \overline{D}}$, repeated according to multiplicity, where $k_{1}+ k_{2}=n$, and
		\item the royal nodes of $x$ are $\sigma_{1}, \dots , \sigma_{n} \in \mathbb{ \overline{D}}$, repeated according to multiplicity.
	\end{enumerate}
There exists some choice of $t_{+} > 0$, $t \in \mathbb{C} \backslash \{0\}$ and $\omega \in \mathbb{ T}$ such that the algorithm of Theorem {\em \ref{cons}} with these choices creates the function $x$.
\end{theorem}

\begin{proof}
	By Theorem \ref{min}, there are polynomials $E_{1}^{1}, E_{2}^{1}$ and $D^{1}$ such that \begin{enumerate} 
		\item $\deg(E_{1}^{1}) , \deg(E_{2}^{1})$ and $\deg(D^{1}) \leq n$,
		\item $D^{1}(\lambda) \neq 0$ \quad $ \text{on} \ \mathbb{ \overline{D }}$,
		\item $E_{2}^{1}(\lambda)= (E_{1}^{1})^{\sim n}(\lambda)$  
		\item $|E_{i}^{1}(\lambda)|\leq |D^{1}(\lambda)|$ on $\mathbb{ \overline{D }}, i=1,2$ and
		\item $x_{1} = \dfrac{E_{1}^{1}}{D^{1}}, \quad x_{2}=\dfrac{E_{2}^{1}}{D^{1}} \quad \text{and} \quad x_{3}= \dfrac{(D^{1})^{ \sim n}}{D^{1}} \quad \text{on} \  \mathbb{ \overline{D }}$.
	\end{enumerate}
By hypothesis, the zeros of $x_{1}$,  repeated according to multiplicity, are $\alpha_{1}^{1},...,\alpha_{k_{1}}^{1}$, and the zeros of $x_{2}$, repeated according to multiplicity, are $\alpha_{1}^{2},...,\alpha_{k_{2}}^{2}$ where $k_{1}+k_{2}=n$. \newline
By Lemma \ref{lem22}, 
\begin{equation*}
E_{1}^{1}(\lambda)= t\prod_{j=1}^{k_{1}}(\lambda - \alpha_{j}^{1}) \prod_{j=1}^{k_{2}}(1- \overline{\alpha}_{j}^{2}\lambda) \quad \text{for some} \  t \in \mathbb{C}\backslash\{0\} \ \text{and all} \ \lambda \in  \mathbb{ D}. 
\end{equation*}
By hypothesis, $\sigma_{1}, \dots , \sigma_{n}$ are the royal nodes of $x$. Thus, by Proposition \ref{mult}, for the royal polynomial $R^{1}$ of $x$, there exists $t_{+} > 0$ such that 
\begin{equation*}
R^{1}(\lambda)= t_{+} \prod_{j=1}^{n}(\lambda-\sigma_{j})(1-\overline{\sigma}_{j}\lambda).
\end{equation*}
By Proposition \ref{propit}, for $\lambda \in \mathbb{T}$,
\begin{equation*}
\lambda^{-n} R^{1}(\lambda) = |D^{1}(\lambda)|^{2} - |E_{1}^{1}(\lambda)|^{2}.
\end{equation*}
By Theorem \ref{min}, $D^{1}(\lambda) \neq 0$ on $\mathbb{ \overline{D}}$. Hence, for $\lambda \in \mathbb{ T}$
\begin{equation*}
\lambda^{-n} R^{1}(\lambda) + |E_{1}^{1}(\lambda)|^{2} = |D^{1}(\lambda)|^{2}   \neq 0.
\end{equation*}
This implies that $\alpha_{1}^{1}, \dots , \alpha_{k_{1}}^{1} $ and $\alpha_{1}^{2}, \dots , \alpha_{k_{2}}^{2}$ which are on $\mathbb{ T}$ are distinct from $\sigma_{i}$, $i=1, \dots, n$.
By the construction in Theorem \ref{cons}, for $\sigma_{i}$, $i=1, \dots, n$ and $\alpha_{1}^{1}, \dots , \alpha_{k_{1}}^{1}$ and $\alpha_{1}^{2}, \dots , \alpha_{k_{2}}^{2}$, the rational $\mathbb{ \overline{E }}$-inner function $x=(x_{1},x_{2},x_{3})$ can be defined by
\begin{equation*}
\bigg(  \frac{E_{1}}{D}, \frac{E_{1}^{\sim n}}{D} , \frac{D^{\sim n}}{D}  \bigg) \ 
\end{equation*}
for a suitable choice of $t_{+}> 0$, $t \in \mathbb{C} \backslash\{0\}$ and $\omega \in \mathbb{ T}$. Since $E_{1}^{1}$ and $R^{1}$ coincide with  $E_{1}$ and $R$ in the construction of Theorem \ref{cons} for a suitable choice of $t_{+} > 0$ and $t \in \mathbb{C} \backslash \{ 0\}$, 
  $D^{1}$ is a permissible choice for $\omega D$  for some choice $\omega \in \mathbb{T}$, as a solution for equation \eqref{eqn26}. 
Thus the algorithm of Theorem \ref{cons} creates $x=(x_{1}, x_{2}, x_{3})$ for the appropriate choices of $t_{+} > 0$, $t \in \mathbb{C} \backslash \{ 0\}$ and $\omega \in \mathbb{T}$.
\end{proof}

\section{Convex subsets of $\mathbb{ \overline{ E}}$ and extremality} \index{convex domain}  \label{convexity}

	In this section we study convex subsets of $\mathbb{ \overline{ E}}$. We show that, for a fixed $x_{3} \in \mathbb{ \overline{D}}$, the subset $\mathbb{ \overline{E}} \cap \big(  \mathbb{C}^{2} \times \{ x_{3}\} \big)$ is convex. Recall that the distinguished boundaries of the tridisc $\mathbb{ D}^{3}$ and the ball $\mathbb{B}_{3}$ contain no line segments. Thus every inner function in the set of analytic functions ${\rm Hol}(\mathbb{ D}, \mathbb{ D}^{3})$ from $\D$ to $ \mathbb{ D}^{3}$ is an extreme point of ${\rm Hol}(\mathbb{ D}, \mathbb{ D}^{3})$ and every inner function in the set of analytic functions $ {\rm Hol}(\mathbb{ D}, \mathbb{B}_{3})$ from $\D$ to $\mathbb{B}_{3}$ is an extreme point of $ {\rm Hol}(\mathbb{ D}, \mathbb{B}_{3})$. However, this property contrasts sharply with the situation in the tetrablock. Despite the fact that the set $\mathcal{J}$ of rational tetra-inner functions\index{set of rational tetra-inner functions $\mathcal{J}$} \index{$\mathcal{J}$} is not convex, the conventional notion of extreme point of $\mathcal{J}$ is well defined and fruitful. In Theorem \ref{parex}, we prove that for $x \in \mathcal{R}^{n,k}$ with $2k \leq n$, $x$ is not an extreme point. A class of extreme points of the set $\mathcal{J}$ is given  in Proposition \ref{exetremE}.
	
\subsection{Convex subsets in the tetrablock}\label{convex}
	\begin{definition}
		A set $\Omega$ in a vector space is {\em convex} if for all $z,w \in \Omega$ and all $t$ such that $0 \leq t \leq 1$, the point $\;tz+(1-t)w\;$ belongs to $\Omega$. 
	\end{definition}

\begin{proposition} {\em \cite[page 8]{AWY}}
	$\mathbb{\overline{E}}$ is not convex.
\end{proposition}
\begin{proof}	Take $x=(i,1,i)$ and $y=(-1,i,-i)$ in 
$\mathbb{ \overline{E }}$. Let us take  $t=1/2$, then the point $w= tx+(1-t)y= \tfrac{1}{2} (-1+i,1+i,0)$. One can show that $w$ is not in $\mathbb{ \overline{E }}$. Therefore $\mathbb{\overline{E}}$ is not convex. 
\end{proof}
Let us show that the set $\mathbb{ \overline{E }}$ is convex in $(x_{1}, x_{2})$ for a fixed $x_{3} \in \mathbb{ \overline{D }}$, that is, the set 
\begin{equation*}
	\mathbb{\overline{E}} \cap \big( \mathbb{C}^{2} \times \{ x_{3}\} \big) = \{ x=(x_{1}, x_{2}, x_{3}) \in \mathbb{C}^{3} : |x_{1}-\overline{x_{2}}x_{3}| + |x_{2} - \overline{x_{1}}x_{3}| \leq 1 - |x_{3}|^{2}  \}
\end{equation*}
is convex for every $x_{3} \in \mathbb{ \overline{D}}$.
\begin{proposition} \label{propconv}
The following sets are convex:
\begin{enumerate}
\item $\overline{\mathbb{E}} \cap \big(\mathbb{C}^{2} \times \{ x_{3}\} \big)$ for any $x_{3} \in \mathbb{ \overline{D }}$;
\item $b\mathbb{ \overline{E }} \cap \big( \mathbb{C}^{2} \times \{ x_{3}\} \big) $ for any $x_{3} \in \mathbb{ \overline{D }}$.
\end{enumerate}
	\end{proposition}
\begin{proof}	(1) Let $x, y  \in \mathbb{\overline{E}} \cap \big( \mathbb{C}^{2} \times \{ x_{3}\} \big) $, and so, by Theorem \ref{tetbnd2}, $x$ and $y$ satisfy the inequalities 
\begin{equation} \label{EQ18}
         |x_{1}-\overline{x}_{2}x_{3}| +|x_{2}- \overline{x}_{2}x_{3}| \leq 1 - |x_{3}|^{2}
         \end{equation} 
and 
\begin{equation} \label{EQ19}
         |y_{1}-\overline{y}_{2}x_{3}| +|y_{2}- \overline{y}_{1}x_{3}| \leq 1 - |x_{3}|^{2}
         \end{equation}
respectively.
For all $t \in [0,1]$, 
\begin{eqnarray*}
          w=tx + (1-t)y &=& t(x_{1}, x_{2}, x_{3} ) + (1-t)(y_{1}, y_{2}, x_{3}) \\ 
          &=& \big(tx_{1} + (1-t) y_{1} , tx_{2} + (1-t) y_{2}, tx_{3} + (1-t) x_{3} \big)  \\ 
          &=& \big(tx_{1} + (1-t) y_{1} , tx_{2} + (1-t) y_{2}, x_{3} \big).
\end{eqnarray*}
Let us check that the point $w \in \mathbb{C}^{3}$ is in the set 
           $\mathbb{\overline{E}} \cap ( \mathbb{C}^{2} \times \{ x_{3}\})$. By Theorem \ref{tetbnd2}, $w \in \mathbb{ \overline{E }}$ if and only if    
           \begin{equation*}
          \underbrace{|w_{1}- \overline{w}_{2}x_{3}|}+ \underbrace{|w_{2}- \overline{w}_{1}x_{3}|} \leq 1-|x_{3}|^{2}.
           \end{equation*}
          Let us consider the first term on the left hand side
          \begin{eqnarray}  \label{EQ30}
          |w_{1}- \overline{w}_{2}x_{3}|&=& 
          \big|  tx_{1} + (1-t) y_{1} - \big(t\overline{x}_{2} + (1-t) \overline{y}_{2}\big) x_{3}  \big| \\ \nonumber &=&  \big|t(x_{1}- \overline{x}_{2}x_{3})+ (1-t)(y_{1}-\overline{y}_{2}x_{3}) \big| \\ \nonumber
          &\leq& t\big|x_{1}- \overline{x}_{2}x_{3} \big| + (1-t) \big| y_{1}-\overline{y}_{2}x_{3} \big| .
          \end{eqnarray}
          For the second term of the left hand side we have 
          \begin{eqnarray}  \label{EQ31}
          |w_{2}- \overline{w}_{1}x_{3}|&=&
          \big|  tx_{2} + (1-t) y_{2} - \big(t\overline{x}_{1} + (1-t) \overline{y}_{1}\big) x_{3}  \big| \\ \nonumber &=& \big|t(x_{2}- \overline{x}_{1}x_{3})+ (1-t)(y_{2}-\overline{y}_{1}x_{3}) \big| \\ \nonumber
          &\leq&  t\big|x_{2}- \overline{x}_{1}x_{3} \big| + (1-t) \big| y_{2}-\overline{y}_{1}x_{3} \big| .
          \end{eqnarray}
          Add inequalities \eqref{EQ30} and \eqref{EQ31} we get
          \begin{multline*}
          |w_{1}- \overline{w}_{2}x_{3}| + |w_{2}- \overline{w}_{1}x_{3}| 
          \leq \\ 
t\big|x_{1}- \overline{x}_{2}x_{3} \big| + (1-t) \big| y_{1}-\overline{y}_{2}x_{3} \big|  + t\big|x_{2}- \overline{x}_{1}x_{3} \big|+ (1-t) \big| y_{2}-\overline{y}_{1}x_{3} \big|.   
           \end{multline*}
           Therefore, by inequalities \eqref{EQ18} and \eqref{EQ19},
       \begin{eqnarray*}
       |w_{1}- \overline{w}_{2}x_{3}| + |w_{2}- \overline{w}_{1}x_{3}| 
       &\leq& t \big(\underbrace{\big|x_{1}- \overline{x}_{2}x_{3} \big|   + \big|x_{2}- \overline{x}_{1}x_{3} \big|} \big) \\ \nonumber &+&  (1-t) \big( \underbrace{\big|  y_{1}-\overline{y}_{2}x_{3} \big|  + \big| y_{2}-\overline{y}_{1}x_{3} \big|} \big) \\ \nonumber
       &\leq& t(1-|x_{3}|^{2}) + (1-t)(1-|x_{3}|^{2}) \\ \nonumber
       &=& 1-|x_{3}|^{2}.
       \end{eqnarray*}
          Hence for all $t \in [0,1]$, $w= tx + (1-t)y \in \mathbb{\overline{E}} \cap \big( \mathbb{C}^{2} \times \{ x_{3}\} \big)$. Therefore $\mathbb{\overline{E}} \cap \big( \mathbb{C}^{2} \times \{ x_{3}\} \big)$ is convex for any fixed $x_{3} \in \mathbb{\overline{D}}$. \\\\\
          (2) Let $x_{3} \in \mathbb{ \overline{D }}$ and $x, y \in b\mathbb{ \overline{E }} \cap \big(   \mathbb{C}^{2} \times \{ x_{3}\}  \big)$,  where $x=(x_{1}, x_{2}, x_{3})$ and $y=(y_{1}, y_{2}, x_{3})$. Note that, by Theorem \ref{imp} $(1)$,
          \begin{equation} \label{Cond2}
           w \in \mathbb{C}^{3} \ \text{belongs to} \  b\mathbb{\overline{E}} \quad \text{if and only if} \quad w_{1}= \overline{w}_{2}w_{3}, \; |w_{2}| \leq 1 \; \text{and} \; |w_{3}|=1.
          \end{equation}
          Thus we have 
          \begin{equation*}
          x_{1}= \overline{x}_{2}x_{3}, \quad |x_{2}| \leq 1 \quad \text{and} \ \ |x_{3}|=1,
          \end{equation*}
          and
          \begin{equation*}
          y_{1}= \overline{y}_{2}x_{3}, \quad |y_{2}| \leq 1 \quad \text{and} \ \ |x_{3}|=1.
          \end{equation*}
          For $t$ such that $ 0 \leq t \leq 1$, let  $$w=\big(w_{1}, w_{2}, w_{3}\big)= tx + (1-t)y = \big(tx_{1} + (1-t) y_{1} , tx_{2} + (1-t) y_{2}, x_{3}\big ).$$ To prove the convexity of $b\mathbb{ \overline{E }} \cap ( \mathbb{C}^{2} \times \{ x_{3}\} )$, we need to check that, for all $t$ such that $0\leq t \leq 1$,  $w$ lies in $b\mathbb{\overline{E}} \cap \big(  \mathbb{C}^{2} \times \{ x_{3}\} \big)$, that is, it satisfies condition \eqref{Cond2}. \newline
Note that
\begin{eqnarray*} 
          \overline{w}_{2}w_{3} &=&  \big(\overline{tx_{2} + (1-t) y_{2} }\big)x_{3} \\
          &=&t\overline{x}_{2}x_{3} + (1-t) \overline{y}_{2} x_{3} \\ 
          &=& tx_{1} +(1-t)y_{1} = w_{1}
          \end{eqnarray*}
and 
\begin{eqnarray*}
          |w_{2}|&=& |tx_{2} + (1-t) y_{2}| \\
          &\leq& t|x_{2}| + (1-t)|y_{2}| \\
          &\leq& t + 1 - t = 1.
          \end{eqnarray*}
      Obviously, $|w_{3}|=|x_{3}|=1$. Therefore the set $b\mathbb{ \overline{E }} \cap ( \mathbb{C}^{2} \times \{ x_{3}\} )$ is convex for any fixed $x_{3} \in \mathbb{ \overline{D }}$. 
    \end{proof}
\begin{lemma}
		Let $x=(x_{1}, x_{2},x_{3})$, $x^{1}=(x^{1}_{1}, x^{1}_{2},x^{1}_{3})$ and $x^{2}=(x^{2}_{1}, x^{2}_{2},x^{2}_{3})$ be in $b\mathbb{ \overline{E}}$ and satisfy  $x=tx^{1}+(1-t)x^{2}$ for some $t \in (0,1)$. Then $x_{3}=x_{3}^{1}=x_{3}^{2}$.
\end{lemma}
		\begin{proof}
			Since $x$, $x^{1}, x^{2} \in b\mathbb{ \overline{ E}}$, by Theorem \ref{imp},
			\begin{equation*}
			|x_{3}|= 1, \qquad |x_{3}^{1}|=1 \qquad \text{and} \qquad |x_{3}^{2}|=1.
			\end{equation*}  
			By assumption $x_{3}=tx_{3}^{1}+(1-t)x_{3}^{2}$. Since every point of $\mathbb{ T}$ is an extreme point of $\mathbb{ \overline{ D}}$, $x_{3}=x_{3}^{1}=x_{3}^{2}$.
		\end{proof}

\subsection{Extremality in the set of $\mathbb{\overline{ E}}$-inner functions} \label{Sec6.2}

In this section we show that, for a fixed inner function $x_{3}$, the set of  rational $\mathbb{ \overline{ E}}$-inner functions $x=(x_{1},x_{2},x_{3})$ with third component $x_{3}$ is a convex set in $\mathcal{J}$.
 We prove that an $\mathbb{ \overline{ E}}$-inner function $x$ is not an extreme point of the set $\mathcal{J}$ if the number of the royal nodes of $x$ on $\mathbb{ T}$, counted with multiplicity, is less than or equal to half of the degree of $x$. In Proposition \ref{exetremE} we give a class of extreme rational $\mathbb{ \overline{ E}}$-inner functions $x \in \mathcal{R}^{n,k}$ of the set $\mathcal{J}$ for which $2k>n$.	
	\begin{theorem} \label{conv13}
		For a fixed inner function  $x_{3}$, the set of~   
$~\mathbb{\overline{E }}$-inner functions $(x_{1}, x_{2}, x_{3})$  is convex.
	\end{theorem}
	\begin{proof}	For the fixed inner function $x_{3}$, let $x=(x_{1}, x_{2}, x_{3})$ and $y=(y_{1}, y_{2}, x_{3})$ be $\mathbb{ \overline{E }}$-inner functions. For   $0 \leq t \leq 1$ and  $\lambda \in \mathbb{ D}$,
		\begin{equation*}
		\big( tx + (1-t)y\big) (\lambda) = \big( tx_{1} + (1-t) y_{1} , tx_{2} + (1-t) y_{2}, x_{3}\big)  (\lambda).
		\end{equation*}
		The function  $$w(\lambda)=\big( w_{1}, w_{2}, w_{3} \big)(\lambda)= \big( tx_{1} + (1-t) y_{1} , tx_{2} + (1-t) y_{2}, x_{3}\big)  (\lambda), \quad \lambda \in \mathbb{D}, $$ is analytic on $\mathbb{D}$ and, by Proposition \ref{propconv} (1), $w(\mathbb{D}) \subseteq \mathbb{ \overline{E }}$. By Proposition \ref{propconv} (2), since for almost all  $\lambda \in \mathbb{ T}$, $x(\lambda)$ and $y(\lambda)$ are in $b\mathbb{ \overline{E }}$, $w(\lambda)$ has also to be in $b\mathbb{ \overline{E }}$. Thus $w$ is an $\mathbb{ \overline{ E}}$-inner function.
		Therefore the set of $\mathbb{\overline{E }}$-inner functions $x=(x_{1}, x_{2}, x_{3})$ is convex for any fixed  inner function $x_{3}$.
	\end{proof}

\begin{definition}
	A rational $\mathbb{\overline{E }}$-inner function $x$ is an {\em extreme point} \index{extreme point of $\mathcal{J}$} of $\mathcal{J}$ if whenever $x$ has a representation of the form $x=tx^{1}+ (1-t)x^{2}$ for $t \in (0,1)$ and $x^{1}, x^{2}$ are rational $\mathbb{ \overline{E }}$-inner functions, $x^{1}=x^{2}$.
\end{definition}
	We will show below that $\mathcal{J}$ is not convex, however the notion of extreme points still has the usual sense. 
	\begin{lemma} \label{lem14}
		Let $x=(x_{1},x_{2},x_{3})$, $x^{1}=(x_{1}^{1},x_{2}^{1},x_{3}^{1})$ and $x^{2}=(x_{1}^{2},x_{2}^{2},x_{3}^{2})$ be rational $\mathbb{ \overline{E }}$-inner functions. If $x=tx^{1}+ (1-t)x^{2}$ for some $t \in (0,1)$ then $x_{3}=x_{3}^{1}=x_{3}^{2}$.
	\end{lemma} 
	\begin{proof}
		Since $x=tx^{1} + (1-t)x^{2}$, we have 
		\begin{equation*}
		\big(x_{1}, x_{_2}, x_{3}\big)= \big( tx_{1}^{1}, tx_{2}^{1}, tx_{3}^{1}\big) + \bigg( (1-t)x_{1}^{2}, (1-t)x_{2}^{2}, (1-t)x_{3}^{2} \bigg).
		\end{equation*}
		Thus 
		$x_{3}= tx_{3}^{1} + (1-t)x_{3}^{2}$.
		Hence, for every point $\lambda \in \mathbb{T}$, 
		$$x_{3}(\lambda)= tx_{3}^{1}(\lambda) + (1-t)x_{3}^{2}(\lambda).$$
	
By assumption, $x^{1}$ and $x^{2}$ are rational $\mathbb{ \overline{ E}}$-inner functions, and so, by Lemma \ref{x3inner} (2), $x_{3}^{1}$ and $x_{3}^{2}$ are rational inner functions, that is, for all $\lambda \in \mathbb{T}$,
\begin{equation*}
 |x_{3}^{1}(\lambda)|=1 \qquad \text{and} \qquad |x_{3}^{2}(\lambda)|=1 .
\end{equation*}
Every point of $\mathbb{T}$ is an extreme point of $\mathbb{ \overline{D }}$, and therefore,
\begin{equation*}
x_{3}(\lambda)= x_{3}^{1}(\lambda) = x_{3}^{2}(\lambda)
\end{equation*}
for all $\lambda \in \mathbb{ T}$. Since $x^{1}$ and $x^{2}$ are rational functions, $x_{3}=x_{3}^{1}=x_{3}^{2}$.
\end{proof}

\begin{lemma}
	The set of rational $\mathbb{ \overline{ E}}$-inner functions $\mathcal{J}$ is not convex.
\end{lemma}
\begin{proof}
	Suppose that $x^{1}= \big( x_{1}^{1},x_{2}^{1},x_{3}^{1} \big) \in \mathcal{J}$ and $x^{2}=(x_{1}^{2},x_{2}^{2},x_{3}^{2}) \in \mathcal{J}$ such that $x_{3}^{1} \neq x_{3}^{2}$. Then by Lemma \ref{lem14}, $x=tx^{1}+ (1-t)x^{2}$ is not in $\mathcal{J}$ for all $t \in (0,1)$. Therefore $\mathcal{J}$ is not convex.
\end{proof}

For  an inner function $p$ of degree $n$,  let $\mathcal{R}_p^n$ be the set of rational $\mathbb{ \overline{ E}}$-inner functions with third component $p$.

\begin{proposition}\label{prop5.10}  The set $\mathcal{R}_p^n$ is convex for every inner function $p$ of degree $n$. For any collection $S$ of rational  $\mathbb{ \overline{ E}}$-inner functions,
$S$ is convex if and only if there exists an inner function $p$ of degree $n$ such that $S$ is a convex subset of $\mathcal{R}_p^n$.
\end{proposition}
\begin{proof} It follows from Theorem \ref{conv13} and Lemma \ref{lem14}.
\end{proof}

\begin{proposition}\label{prop5.20}
Let $x$ be a rational $\mathbb{\overline{ E}}$-inner function
 of degree $n$. Then $x$ is a convex combination of at most $2n+3$ extreme rational $\mathbb{\overline{ E}}$-inner functions of degree at most $n$.
\end{proposition}
\begin{proof} By assumption, $x=(x_{1},x_{2}, x_{3})$ is a rational $\mathbb{ \overline{ E}}$-inner function of degree $n$, and so  $x_3$ is an inner function of degree $n$. Thus
$x \in \mathcal{R}_{x_3}^n$. By Remark \ref{dimE}, the convex  set 
$\mathcal{R}_{x_3}^n$ is a subset of a  ($2n+2$)-dimensional real subspace of the rational functions. Therefore, by a theorem of Carath$\acute{\rm e}$odory \cite{Ca,St}, $x$ is a convex combination of at most $2n+3$ extreme rational  $\mathbb{ \overline{ E}}$-inner functions of degree at most $n$.
\end{proof}


\begin{definition} \label{order} Let $f$ be 
	a real or complex-valued function on a real interval $I$. We say that $f$ takes a value $y$ to order $m \geq 1$ at a point $t_{0} \in I$ if $f \in C^{m}(I)$, $f(t_{0})=y$, $f^{(j)}(t_{0})=0$ for $j=1, \dots m-1$ and $f^{(m)}(t_{0}) \neq 0$. We say that $f$ vanishes to order $m \geq 1$ at a point $t_{0} \in I$ if $f$ takes the value $0$ to order $m$ at $t_{0}$.
\end{definition}

\begin{lemma} \label{rem1}
	Let $f \in C^{m}(I)$, $f(t_{0})=y$ at $t_{0} \in I$, and let $y \neq 0$. If $f^{2}$ takes the value $y^{2}$ to order $m \geq 1$ at $t_{0}$, then $f$ takes the value $y$ to order $m$ at $t_{0}$.
\end{lemma}
\begin{proof}
	Let $I$ be a real interval and let $f \in C^{m}(I)$. Suppose that $f^{2}$ takes the value $y^{2}$ to order $m$ at $t_{0}$. Then, by Definition \ref{order}, 
	\begin{equation} \label{alleq}
	f^{2}(t_{0}) = y^{2}, \; [f^{2}]^{(1)}(t_{0})= [f^{2}]^{(2)}(t_{0})= \dots= [f^{2}]^{(m-1)}(t_{0}) =0, \;  [f^{2}]^{(m)}(t_{0})\neq 0.
	\end{equation}
One can check that 
\begin{equation*} 
f(t_{0}) = y, \qquad f^{(j)}(t_{0})=0, \quad \text{for} \ j=1,\dots , m-1 \qquad \mbox{and}\quad  f^{(m)}(t_{0})\neq 0.
\end{equation*}
\end{proof}

\begin{definition} \label{analonT}
	A function $f$ is analytic on $\mathbb{ T}$ if there exists a function $g$ analytic in a neighbourhood $U_{\mathbb{ T}}$ of $\mathbb{ T}$ such that $f=g|_{\mathbb{ T}}$.
\end{definition}
\begin{lemma} \label{lem77}
	Let $\tau = e^{it_{0}}$ and let $f(t)= (e^{it}-\tau)^{2v} G(e^{it})$ in a neighbourhood of $t_{0}$ where $G(z)$ is analytic on $\mathbb{ T}$ and $G(\tau) \neq 0$. Then 
	\begin{equation}
	f^{(j)}(t_{0}) =0 \ \ \text{for} \ \ j=0,1, \dots , 2v-1\quad \text{and} \ \ f^{(2v)}(t_{0}) \neq 0.
	\end{equation}

\end{lemma}
\begin{proof}
	Since $G$ is analytic on $\mathbb{ T}$, by Definition \ref{analonT}, there exists $U_{\mathbb{T}}$ a neighbourhood of $\mathbb{T}$ and there exists $\tilde{G}$ analytic on $U_{\mathbb{ T}}$ such that $G=\tilde{G}|_{\mathbb{ T}}$.  Let $z=e^{it}$, $\phi(z)= (z-\tau)^{2v}G(z)$ and $\tilde{\phi}(z)=(z-\tau)^{2v}\tilde{G}(z)$. Define $\gamma(\tau, r)$ to be an anticlockwise circle centred at $\tau$ with radius $r$  $$\gamma(\tau, r)= \{z \in \mathbb{ C} : |z-\tau|=r\},$$
	where $r$ is taken sufficiently small that $\gamma \subset U_{\mathbb{T}}$. Hence the function $\tilde{\phi}$ is analytic inside the curve $\gamma$. 
	Therefore, by Cauchy's integral formula,
	\begin{eqnarray} \nonumber \label{eq1234}
	\tilde{\phi}^{(j)}(\tau) &=& \frac{j !}{2\pi i} \int_{\gamma} \frac{\tilde{\phi}(z)}{(z-\tau)^{j+1}} dz,  \\ \nonumber
	&=& \frac{j !}{2\pi i} \int_{\gamma} \frac{(z-\tau)^{2v}\tilde{G}(z)}{(z-\tau)^{j+1}} dz \\ 
	&=& \frac{j !}{2\pi i} \int_{\gamma} (z-\tau)^{2v-j-1}\tilde{G}(z) dz .
	\end{eqnarray}
	For $j$ such that $0 \leq j \leq 2v-1$, the function $(z-\tau)^{2v-j-1}\tilde{G}(z)$ is analytic on $U_{\mathbb{ T}}$. Therefore, by Cauchy's Theorem, 
	\begin{equation} \label{dev}
	\tilde{\phi}^{(j)}(\tau) = \frac{j !}{2\pi i} \int_{\gamma} (z-\tau)^{2v-j-1}\tilde{G}(z) dz = 0. 
	\end{equation}
	If $j=2v$, then equation \eqref{eq1234} becomes
	\begin{equation*}
	\tilde{\phi}^{(2v)}(\tau) = \frac{(2v) !}{2\pi i} \int_{\gamma} \frac{\tilde{G}(z)}{(z-\tau)} dz.
	\end{equation*}
	By Cauchy's integral formula,
	\begin{equation} \label{neq}
	\tilde{\phi}^{(2v)}(\tau) = \frac{(2v) !}{2\pi i} \int_{\gamma} \frac{\tilde{G}(z)}{(z-\tau)} dz = (2v) ! \ G(\tau) \neq 0.
	\end{equation}
	Hence $\phi^{(2v)}(\tau) \neq 0$ because $\tilde{\phi}$ agrees with $\phi$ on $\mathbb{T}$.
		Note that, \begin{equation*}
		f(t)=(e^{it}-e^{it_{0}})^{(2v)}G(e^{it})= \phi(e^{it}).
		\end{equation*}
By the chain rule, 
\begin{eqnarray*}
\frac{df}{dt}&=& \frac{d\phi}{dz} \ \frac{dz}{dt} \\
 \frac{d^{2}f}{dt^{2}}&=& \frac{d^{2}\phi}{dz^{2}} \bigg( \frac{dz}{dt}\bigg)^{2} + \frac{d\phi}{dz} \frac{d^{2}z}{dt^{2}} \\
\dots &=& \dots\\
\frac{d^{2v-1}f}{dt^{2v-1}}&= &\frac{d^{2v-1}\phi}{dz^{2v-1}} \bigg( \frac{dz}{dt} \bigg)^{2v-1} + \dots + \frac{d\phi}{dz} \ \frac{d^{2v-1}z}{dt^{2v-1}}.
\end{eqnarray*}	
By equation \eqref{dev} and since $\tilde{\phi}$ and $\phi$ agree on $\mathbb{T}$,
 $$\dfrac{d^{j}\tilde{\phi}}{dz^{j}}(\tau)=0, \qquad \text{for} \ j=1, \dots , 2v-1,$$ and so, 
$$\dfrac{d^{j}\phi}{dz^{j}}(\tau)= 0,\qquad  \text{for} \ j=1, \dots , 2v-1.$$ 
Therefore, $f^{(j)}(t_{0})= 0$ for $j=1, \dots, 2v-1$. 
For the $(2v)$th derivative of $f$, we have 
$$ \frac{d^{2v}f}{dt^{2v}}=  \frac{d^{2v}\phi}{dz^{2v}} \bigg( \frac{dz}{dt} \bigg)^{2v} + \dots + \frac{d\phi}{dz}\frac{d^{2v}z}{dt^{2v}}.$$
By equations \eqref{dev} and \eqref{neq},
\begin{equation*}
\dfrac{d^{j}\phi}{dz^{j}}(\tau)= 0,\qquad  \text{for} \ j=1, \dots , 2v-1 \qquad \text{and} \qquad \dfrac{d^{2v}\phi}{dz^{2v}}(\tau)\neq 0.
\end{equation*}
Hence $f^{(2v)}(t_{0}) = \dfrac{d^{2v}\phi}{dz^{2v}} (\tau) \bigg( \dfrac{dz}{dt} \bigg)^{2v}(t_{0})\neq 0$.

	Therefore $f^{(j)}(t_{0}) = 0$ for $j=0, \dots , 2v-1$ and $f^{(2v)}(t_{0}) \neq 0$.
\end{proof}

\begin{lemma} \label{lem13}
	Let $x=(x_{1},x_{2},x_{3})$ be a rational $\mathbb{ \overline{E }}$-inner function. For $\tau \in \mathbb{T}$,
	\begin{enumerate}
	\item $|x_{1}(\tau)|=1 \Leftrightarrow \ \tau \  \text{is a royal node of} \ x$;
	\item $|x_{2}(\tau)|=1 \Leftrightarrow \ \tau \ \text{is a royal node of} \ x$.  \\ 
Moreover,  
	  	\item $\tau = e^{it_{0}}$ is a royal node of $x$ of multiplicity $v$ if and only if $|x_{1}(e^{it})|=1$ to order $2v$ at $t=t_{0}$;
	   \item $\tau = e^{it_{0}}$ is a royal node of $x$ of multiplicity $v$ if and only if $|x_{2}(e^{it})|=1$ to order $2v$ at $t=t_{0}$. 
	\end{enumerate}
\end{lemma}

\begin{proof}
	(1) If $\tau = e^{it_{0}}$ is a royal node of $x$ of multiplicity $v$, by Definition \ref{def11},   \begin{equation} \label{eqn113}
	(x_{3}-x_{1}x_{2})(\lambda)=(\lambda - \tau)^{2v} F(\lambda), 
	\end{equation}
	where $F$ is a rational function, analytic on $\mathbb{ T}$ and $F(\tau) \ne 0$. By Lemma \ref{ddb}, since $x$ is an $\mathbb{ \overline{E }}$-inner function,  $x_{2}=\overline{x_{1}}x_{3}$ on $\mathbb{ T}$. Therefore, for $\lambda \in \mathbb{ T}$,
	\begin{eqnarray} \nonumber \label{eqn114}
	\Big( x_{3}-x_{1}x_{2}\Big)(\lambda)&=& \Big(x_{3}-x_{1}\overline{x_{1}}x_{3}\Big)(\lambda) \\ \nonumber
	&=& x_{3}(\lambda) - x_{3}(\lambda)|x_{1}(\lambda)|^{2} \\ 
	&=& x_{3}(\lambda)(1-|x_{1}(\lambda)|^{2}) .
	\end{eqnarray}
	Therefore, for any $\lambda \in \mathbb{ T}$,
	\begin{equation*}
	|x_{1}(\lambda)|=1 \qquad \iff \qquad (x_{3}-x_{1}x_{2})(\lambda) =0,
	\end{equation*}
	that is, if and only if $\lambda$ is a royal node of $x$. 
Hence 	 $\tau \in \mathbb{ T}$ is a royal node of $x$ if and only if $|x_{1}(\tau)|=1 $.

	(2) Since $x$ is rational $\mathbb{ \overline{E}}$-inner function, by Theorem \ref{imp}, $x_{1}=\overline{x_{2}}x_{3}$ on $\mathbb{ T}$. The rest of the proof is similar to the above proof of (1).\\ 

	(3) Suppose that $\tau=e^{it_{0}}$ is a royal node of $x$ of multiplicity $v \geq 1$. Then on combining equations \eqref{eqn113} and \eqref{eqn114}, we have, for all $t \in \mathbb{R}$,
	\begin{equation*}
	x_{3}(e^{it})(1-|x_{1}(e^{it})|^2)= (e^{it}- \tau)^{2v}F(e^{it}).
	\end{equation*}
	This gives
	\begin{equation*}
	1-|x_{1}(e^{it})|^{2}=(e^{it}-\tau)^{2v} \frac{F(e^{it})}{x_{3}(e^{it})}.
	\end{equation*}
	The rational function $G= \dfrac{F}{x_{3}}$ is analytic on $\mathbb{ T}$ and is not equal to zero at $\tau = e^{it_{0}}$. Thus we have
	\begin{equation*}
	1-|x_{1}(e^{it})|^{2}=(e^{it}-\tau)^{2v} G(e^{it}).
	\end{equation*}
	Since $x$ is rational and $|x_{1}(e^{it_{0}})|=1$, the function $f(t)= 1-|x_{1}(e^{it})|^{2}$ is $C^{\infty}$ on a neighbourhood of $t_{0}$. By Lemma \ref{lem77},
	
	\begin{equation*}
	f^{(j)}(t_{0})=0 \quad \text{for} \  j=0,1, \dots , 2v-1 \quad \text{and} \quad f^{(2v)}(t_{0}) \ne 0.
	\end{equation*}  
	Therefore $f$ takes the value $0$ to order $2v$ at $t_{0}$, which implies, by Lemma \ref{rem1}, $|x_{1}(e^{it})|=1$ to order $2v$ at $t_{0}$.

	(4) The proof of this statement follows from (2) and is similar to the above proof of (3). 
\end{proof}

For  an inner function $p$ of degree $n$ and $k=0,1,\ldots,n$, let
\be\label{eq5.20}
\mathcal{R}_p^{n,k}=\{(x_1, x_2, x_3) \in  \mathcal{R}^{n,k}: x_3=p \}.
\ee

\begin{lemma}\label{convex-nodes} Let 
$x=(x_{1},x_{2}, x_{3}) \in \mathcal{R}_{x_3}^{n,k}$
and let $\tau_1, \tau_2,  \dots, \tau_k \in \T$ be royal nodes of $x$. Suppose $x =  t x^1 + (1-t) x^2$ for some $t$ such that $0 < t < 1$, where 
$x^1=(x_{1}^1,x_{2}^1, x_{3}^1)$ and $x^2=(x_{1}^2,x_{2}^2, x_{3}^2)$ are rational $\mathbb{ \overline{ E}}$-inner functions.
Then  $x_{3} = x_{3}^1= x_{3}^2$,
\[ 
x_1^i(\tau_j) = x_1(\tau_j) \; \; \text{for } \;  j =1,\dots k\; \text{and} \;\: i=1,2,
\]
and
\[ 
x_2^i(\tau_j) = x_2(\tau_j) \; \; \text{for } \;  j =1,\dots k\; \text{and} \;\; i=1,2.
\]
Moreover, $\tau_1, \tau_2,  \dots, \tau_k \in \T$ are royal nodes of $x^1$ and $x^2$.
\end{lemma}
\begin{proof} By Lemma \ref{lem14},  $x_3 = x_3^1= x_3^2$. By Lemma \ref{lem13}, $|x_1(\tau_j)| =1$ and $|x_2(\tau_j)| =1$ at each royal node $\tau_j \in \T$.
By assumption,  
$$x_1(\tau_j) =  t x_1^1(\tau_j) + (1-t) x_1^2(\tau_j)$$ 
for  $t$ such that  $0 < t < 1$  and $|x_1^i(\tau_j)| \le 1$  for $ j =1,\dots,k$ and $i=1,2$. Similarly, 
$$x_2(\tau_j) =  t x_2^1(\tau_j) + (1-t) x_2^2(\tau_j)$$ 
for  $t$ such that  $0 < t < 1$  and $|x_2^i(\tau_j)| \le 1$  for $ j =1,\dots,k$ and $i=1,2$. 
Every point on the circle $\T$ is an extreme point of $\bar{\d}$, and so 
 $ x_1^i(\tau_j) = x_1(\tau_j)$
  and
$ x_2^i(\tau_j) = x_2(\tau_j)$ for  $ j=1,\dots, k$ and $i=1,2$.
Therefore $|x_1^i(\tau_j)| = 1$, $|x_2^i(\tau_j)| = 1$ for  $ j=1,\dots, k$ and $i=1,2$.
By Lemma \ref{lem13}, $\tau_1, \tau_2,  \dots, \tau_k \in \T$ are royal nodes of $x^1$ and $x^2$.

\end{proof}


\begin{lemma} \label{parext}
	Let $n\geq 1$. Any $x=(x_{1}, x_{2}, x_{3}) \in \mathcal{R}^{n,0}$ is not an extreme point of $\mathcal{J}$.
\end{lemma}

\begin{proof}
Since $x$ has no royal nodes on $\mathbb{ T}$, by Lemma \ref{lem13}, for all  $\lambda \in \mathbb{ T}$ ,
\begin{equation*}
|x_{1}(\lambda)|   <1 \qquad \text{and} \qquad |x_{2}(\lambda)|   <1 .
\end{equation*}
Since $\mathbb{ T}$ is compact, the supremum of $x_{1}$ and $x_{2}$ is attained on $\mathbb{T}$, that is, there exist $\lambda_{1}, \lambda_{2} \in \mathbb{ T}$ such that
\begin{equation} \label{ek1}
\sup_{\lambda \in \mathbb{T}} |x_{1}(\lambda)|= |x_{1}(\lambda_{1})| < 1 \qquad \text{and} \qquad \sup_{\lambda \in \mathbb{T}} |x_{2}(\lambda)|= |x_{2}(\lambda_{2})| < 1.
\end{equation}
Choose $\eps_{1} >0$ and $\eps_{2}>0$ such that
\begin{equation} \label{ek2}
|x_{1}(\lambda_{1})| (1 + \eps_{1})<1 \qquad \text{and} \qquad  |x_{2}(\lambda_{2})| (1+ \eps_{2}) <1.
\end{equation}
Take $\eps =\min\{\eps_{1}, \eps_{2}\}$. If $x_{1}(\lambda_{1} )=0$, then $$x_{1}(\lambda)=0 \quad  \text{for all} \ \lambda \in \mathbb{ T}.$$ Likewise, if $x_{2}(\lambda_{2})=0$ then  $$x_{2}(\lambda)=0, \quad \text{for all}  \ \lambda \in \mathbb{ T}.$$ 

  Define $x^{1}$ and $x^{2}$ to be 
\begin{equation*}
x^{1}= \big( (1 + \eps)x_{1}, (1 + \eps) x_{2}, x_{3} \big) \qquad \text{and} \qquad x^{2}= \big( (1 - \eps)x_{1} , (1 - \eps)x_{2} , x_{3}\big).
\end{equation*}
Since $x=(x_{1}, x_{2}, x_{3})$ is a rational $\mathbb{ \overline{ E}}$-inner function, for almost all $\lambda \in \mathbb{T}$,
\begin{equation} \label{ek3}
x_{1}(\lambda)= \overline{x}_{2}(\lambda)x_{3}(\lambda), \quad |x_{2}(\lambda)| \leq 1 \quad \text{and} \quad |x_{3}(\lambda)|=1.
\end{equation}
 Let us check that $x^{1}$ and $x^{2}$ are rational $\mathbb{ \overline{ E}}$-inner functions. By Theorem \ref{imp} (1), this will follow if we show that
 \begin{equation*}
 (1+\eps)x_{1}(\lambda)= (1+\eps)\overline{x}_{2}(\lambda)x_{3}(\lambda), \quad (1+\eps)|x_{2}(\lambda)| \leq 1 \quad \text{and} \quad |x_{3}(\lambda)|=1,
 \end{equation*}
 and $x^{1}(\mathbb{ D}) \subset \mathbb{E}$.
 By equations \eqref{ek3}, we have to show only that $$(1+\eps)|x_{2}(\lambda)|\leq 1 \quad \text{on} \ \mathbb{T}.$$ This statement follows from inequalities \eqref{ek1} and \eqref{ek2}.
 Thus $x^{1}(\mathbb{T}) \subset b\mathbb{ \overline{ E}}$. By Theorem \ref{imp} (2), for almost all $\lambda \in \mathbb{ T}$,
\begin{equation*}
x^{1}(\lambda) \in b\mathbb{ \overline{ E}} \Leftrightarrow \Psi(.,x^{1}(\lambda)) \ \text{is an automorphism of} \ \mathbb{ D}.
\end{equation*}
By the maximum principle, for all $\lambda \in \mathbb{ D}$, $\|\Psi(.,x^{1}(\lambda))\|_{H^{\infty}} < 1$. Therefore, by Theorem \ref{tetbnd}, for all $\lambda \in \mathbb{D}$,  $x^{1}(\lambda) \subseteq \mathbb{E }$. This completes the proof that $x^{1}$ is a rational $\mathbb{ \overline{ E}}$-inner function. \\

In a similar way we can show that $x^{2}$ is a rational $\mathbb{ \overline{ E}}$-inner function.
 Moreover, by Lemma \ref{lem13}, $x^{1}, x^{2}$ have no royal nodes on $\mathbb{ T}$ and therefore $x^{1}, x^{2} \in \mathcal{R}^{n,0}$. However $x=\tfrac{1}{2} x^{1}+ \tfrac{1}{2}x^{2}$, which implies that $x$ cannot be an extreme point of $\mathcal{J}$ since $x^{1} \neq x^{2}$.
\end{proof}

\begin{proposition} \label{prop3}
	Let $x=(x_{1}, x_{2}, x_{3})$ be superficial and $x=tx^{1} + (1-t)x^{2}$ for some $0<t<1$, where 
	$x^{1}=(x_{1}^{1}, x_{2}^{1}, x_{3}^{1})$ and $x^{2}=(x_{1}^{2}, x_{2}^{2}, x_{3}^{2})$ are rational $\mathbb{ \overline{ E}}$-inner functions. Then $x^{1}$ and $x^{2}$ are superficial and $x_{3}= x_{3}^{1}=x_{3}^{2}$.
\end{proposition}
\begin{proof} By Lemma \ref{lem14}, $x_{3}=x_{3}^{1}=x_{3}^{2}$. Suppose, for a contradiction, $x^{1}$ is not superficial. It means there exists $\lambda_{0} \in \mathbb{D}$ such that $x^{1}(\lambda_{0}) \in \mathbb{E}$. We will show that in this case $x(\lambda_{0}) \in \mathbb{E}$, and so $x$ is not superficial.

	 By Theorem \ref{tetbnd}, it is enough to prove that
	\begin{equation} \label{eq611}
	|x_{1}(\lambda_{0})-\overline{ x}_{2}(\lambda_{0})x_{3}(\lambda_{0})| + |x_{2}(\lambda_{0})-\overline{ x}_{1}(\lambda_{0})x_{3}(\lambda_{0})| < 1-|x_{3}(\lambda_{0})|^{2}.
	\end{equation}
	Since $x^{1}(\lambda_{0}) \in \mathbb{E}$ and $x^{2}$ is a rational $\mathbb{ \overline{ E}}$-inner function, this implies that 
	\begin{equation}\label{x^1}
	|x_{1}^{1}(\lambda_{0})-\overline{ x}_{2}^{1}(\lambda_{0})x_{3}^{1}(\lambda_{0})| + |x_{2}^{1}(\lambda_{0})-\overline{ x}_{1}^{1}(\lambda_{0})x_{3}^{1}(\lambda_{0})| < 1-|x_{3} ^{1}(\lambda_{0})|^{2}
	\end{equation}
	and 
	\begin{equation} \label{x^2}
	|x_{1}^{2}(\lambda_{0})-\overline{ x}_{2}^{2}(\lambda_{0})x_{3}^{2}(\lambda_{0})| + |x_{2}^{2}(\lambda_{0})-\overline{ x}_{1}^{2}(\lambda_{0})x_{3}^{2}(\lambda_{0})| \leq 1-|x_{3} ^{2}(\lambda_{0})|^{2}.
	\end{equation}
	Let us begin with the first term on the left hand side of inequality \eqref{eq611}. 
	\begin{eqnarray}  \label{eq612}
	& & |x_{1}(\lambda_{0})- \overline{x}_{2}(\lambda_{0})x_{3}(\lambda_{0})| \nonumber \\ 
&=& \big|  tx_{1}^{1}(\lambda_{0}) + (1-t) x_{1}^{2}(\lambda_{0}) - \big( t\overline{x}_{2}^{1}(\lambda_{0}) + (1-t) \overline{x}_{2}^{2}(\lambda_{0})\big) x_{3}(\lambda_{0})  \big|\nonumber  \\ 
	&\leq & \big| t\big(x_{1}^{1}(\lambda_{0})- \overline{x}_{2}^{1}(\lambda_{0})x_{3}(\lambda_{0})\big) \big| + \big| (1-t) \big(x_{1}^{2}(\lambda_{0})- \overline{x}_{2}^{2}(\lambda_{0})x_{3}(\lambda_{0})\big) \big|.  
	\end{eqnarray}
	The second term on inequality \eqref{eq611}
	\begin{eqnarray}  \label{eq613}
	& &|x_{2}(\lambda_{0})- \overline{x}_{1}(\lambda_{0})x_{3}(\lambda_{0})| \nonumber\\  &=& 
	\big|  tx_{2}^{1}(\lambda_{0}) + (1-t) x_{2}^{2}(\lambda_{0}) - \big(t\overline{x}_{1}^{1}(\lambda_{0}) + (1-t) \overline{x}_{1}^{2}(\lambda_{0})\big) x_{3}(\lambda_{0})  \big| \nonumber  \\ 
	&\leq& \big| t\big(x_{2}^{1}(\lambda_{0})- \overline{x}_{1}^{1}(\lambda_{0})x_{3}(\lambda_{0})\big) \big| + \big| (1-t) \big(x_{2}^{2}(\lambda_{0})- \overline{x}_{1}^{2}(\lambda_{0})x_{3}(\lambda_{0})\big) \big|.  
	\end{eqnarray}
Add inequalities \eqref{eq612} and \eqref{eq613}, and use inequalities \eqref{x^1} and \eqref{x^2} to obtain 
\begin{eqnarray}  \label{eq614} 
& &|x_{1}(\lambda_{0})-\overline{ x}_{2}(\lambda_{0})x_{3}(\lambda_{0})| + |x_{2}(\lambda_{0})-\overline{ x}_{1}(\lambda_{0})x_{3}(\lambda_{0})|\nonumber\\
&\leq&
t \Big( \big|  x_{1}^{1}(\lambda_{0})- \overline{x}_{2}^{1}(\lambda_{0})x_{3}(\lambda_{0}) \big| + \big| x_{2}^{1}(\lambda_{0})- \overline{x}_{1}^{1}(\lambda_{0})x_{3}(\lambda_{0}) \big|\Big)  \nonumber\\ 
 & & + (1-t)\Big( \big|  x_{1}^{2}(\lambda_{0})- \overline{x}_{2}^{2}(\lambda_{0})x_{3}(\lambda_{0}) \big| + \big|  x_{2}^{2}(\lambda_{0})- \overline{x}_{1}^{2}(\lambda_{0})x_{3}(\lambda_{0}) \big| \Big) \nonumber\\ 
	&<& t(1-|x_{3}(\lambda_{0})|^{2}) + (1-t)(1-|x_{3}(\lambda_{0})|^{2}) = 1-|x_{3}(\lambda_{0})|^{2}.
	\end{eqnarray}
\end{proof}
\begin{theorem} \label{parex}
	Let $x \in \mathcal{R}^{n,k}$. If $2k \leq n$, then $x$ is not an extreme point of the set  $\mathcal{J}$ of rational $\mathbb{ \overline{E }}$-inner functions.
\end{theorem}
\begin{proof}
	Let $x \in \mathcal{R}^{n,k}$. By Definition \ref{def11}, $x$ has $n$ royal nodes in $\mathbb{ \overline{D}}$ and $k$ royal nodes that lie in $\mathbb{ T}$. By Theorem \ref{min}, there exist polynomials $E_{1}$, $E_{2}$ and  $D$  of degree at most $n$  such that
	\begin{equation*}
	x= \bigg( \frac{E_{1}}{D},\frac{E_{1}^{\sim n}}{D}, \frac{D^{\sim n}}{D}  \bigg),
	\end{equation*}
where, for all $\lambda \in \mathbb{ \overline{D }}$,  $D(\lambda) \neq 0$  and $E_{2}(\lambda)= E_{1}^{\sim n}(\lambda)$. Let $\tau_{1}, \dots , \tau_{k} \in \mathbb{ T}$ and $\alpha_{k+1}, \dots, \alpha_{n} \in \mathbb{ D}$ be the royal nodes of $x$ in $\mathbb{ \overline{D }}$ repeated according to multiplicity. By Proposition \ref{mult}, the royal polynomial of $x$ is 
\begin{equation*}
	R= r\prod_{j=1}^{k} Q_{\tau_{j}} \prod_{j=k+1}^{n} Q_{\alpha_{j}},
\end{equation*}
for some $r >0$. Thus for all $\lambda \in \mathbb{T}$,
	\begin{eqnarray} \label{eqe49}
	\lambda^{-n}R(\lambda) &=& r\lambda^{-n} \Bigg\{   \prod_{j=1}^{k} (\lambda - \tau_{j})(1-\overline{\tau_{j}}\lambda) \prod_{j=k+1}^{n} (\lambda - \alpha_{j})(1-\overline{\alpha_{j}}\lambda) \Bigg\} \nonumber  \\
	&=& r\prod_{j=1}^{k} |\lambda - \tau_{j}|^{2} \prod_{j=k+1}^{n} |\lambda - \alpha_{j}|^{2}.
	\end{eqnarray}
	By Proposition \ref{propit} and equation \eqref{eqe49}, for all $\lambda \in \mathbb{T}$, 
	\begin{equation} \label{eq65}
	|D(\lambda)|^{2}- |E_{1}(\lambda)|^{2} =\lambda^{-n}R(\lambda)= r \prod_{j=1}^{k} |\lambda - \tau_{j}|^{2} \prod_{j=k+1}^{n} |\lambda - \alpha_{j}|^{2}.
	\end{equation}
	Assume first that $n$ is even and write $n=2m$. This implies that $k \leq m$. Define a polynomial $g$ by
	\begin{equation*}
	g(\lambda)= \overline{\tau}_{1}\dots \overline{ \tau}_{k} \lambda^{m-k} \prod_{j=1}^{k}(\lambda - \tau_{j})^{2}.
	\end{equation*}
Clearly, the polynomial $g$ has degree $m+k \leq n$. Moreover, $g$ is $n$-symmetric since
	\begin{eqnarray*}
	g^{\sim n}(\lambda)= 
\lambda^{n}\overline{ g\big(\frac{1}{\overline{\lambda}} \big)} &=& \lambda^{n} \Bigg\{\overline{ \overline{\tau}_{1}\dots \overline{ \tau}_{k} \frac{1}{(\overline{\lambda})^{m-k}} \prod_{j=1}^{k} (\frac{1}{\overline{\lambda}}- \tau_{j})^{2} }\Bigg\} \\
	&=& \lambda^{2m} \Bigg\{ \tau_{1}\dots  \tau_{k} \frac{1}{\lambda^{m-k}} \prod_{j=1}^{k} (\frac{1}{\lambda}-\overline{ \tau_{j}})^{2} \Bigg\} \\ 
	&=&\overline{\tau}_{1}\dots \overline{ \tau}_{k} \lambda^{m-k} \prod_{j=1}^{k}(\lambda - \tau_{j})^{2} = g(\lambda).
	\end{eqnarray*} 
Let
\begin{equation*}
      E_{1}^{t}= E_{1} + tg \qquad\text{and} \qquad  E_{2}^{t}= E_{1}^{\sim n} + tg \qquad  \text{for} \ t \in \mathbb{R}.
	\end{equation*}
The polynomial $E_{1}^{t}$ has degree at most $n$. We also have, for all $\lambda \in \mathbb{ \overline{D }}$, 
	  \begin{eqnarray} \nonumber \label{eq57}
	  \big( E_{2}^{t}\big)^{\sim n}(\lambda) &=& \big( E_{1}^{\sim n} + tg \big)^{\sim n}(\lambda) \\ 
	  &=& \big( E_{1}^{\sim n} \big)^{\sim n} (\lambda)+ \big( tg \big)^{\sim n}(\lambda)  = \big( E_{1}^{t} + tg\big)(\lambda) = E_{1}^{t}(\lambda).
	  \end{eqnarray}
Note that, on $\mathbb{ T}$, 
	\begin{eqnarray} \nonumber \label{eqn66}
	|D|^{2} - |E_{1}^{t}|^{2} &=& |D|^{2} - |E_{1}+ tg|^{2} \\ \nonumber
	&=& |D|^{2} - (E_{1}+ tg)\overline{(E_{1} + tg) } \\
	&=& |D|^{2} -  |E_{1}|^{2} - t^{2}|g|^{2}-2 \text{Re}(tg\overline{ E}_{1}) . 
	\end{eqnarray}
Let $\|E_{1}\|_{\infty}$ denote the supremum of $|E_{1}|$ on $\mathbb{ T}$. Then, for all $\lambda \in \mathbb{T}$,
\begin{eqnarray} \nonumber \label{eq67}
	\text{Re}(tg(\lambda)\overline{E_{1}}(\lambda) \leq |tg(\lambda)E_{1}(\lambda)|&=& |tE_{1}(\lambda)| \ \bigg| \overline{\tau}_{1} \dots \overline{ \tau}_{j} \lambda^{m-k} \prod_{j=1}^{k} (\lambda - \tau_{j})^{2} \bigg| \\ \nonumber
	&=& |tE_{1}(\lambda)| \prod_{j=1}^{k} |\lambda - \tau_{j}|^{2} \\
	&\leq& |t| \|E_{1}\|_{\infty} \prod_{j=1}^{k} |\lambda - \tau_{j}|^{2}.
\end{eqnarray}
Note that, for all $\lambda \in \mathbb{ T}$,
\begin{equation*}
	|g(\lambda)|^{2} = \big| \overline{\tau}_{1}\dots \overline{ \tau}_{k} \lambda^{m-k} \prod_{j=1}^{k}(\lambda - \tau_{j})^{2}  \big|^{2} = \big| \prod_{j=1}^{k}(\lambda - \tau_{j})^{2} \big|^{2}.
	\end{equation*}
	Combine equations  \eqref{eq65} and \eqref{eqn66} and inequality \eqref{eq67}, for all $\lambda \in \mathbb{ T}$, to get 
	\begin{eqnarray*} 
	& & |D(\lambda)|^{2} - |E_{1}^{t}(\lambda)|^{2} \\ \nonumber &=& |D(\lambda)|^{2} - |E_{1}(\lambda)|^{2} -|t|^{2}|g(\lambda)|^{2} - 2 \text{Re}\big(tg(\lambda)\overline{ E}_{1}(\lambda)\big), \ \text{(by equation} \  \eqref{eqn66}) \ \\
	&=& r \prod_{j=1}^{k} |\lambda - \tau_{j}|^{2} \prod_{j=k+1}^{n} |\lambda -\alpha_{j}|^{2}-|t|^{2}|g(\lambda)|^{2}-2\text{Re}\big(tg(\lambda)\overline{ E}_{1}(\lambda)\big), \ \text{(by equation} \ \eqref{eq65}) \\ 
	&\geq& r \prod_{j=1}^{k} |\lambda - \tau_{j}|^{2} \prod_{j=k+1}^{n} |\lambda -\alpha_{j}|^{2}-|t|^{2}|g(\lambda)|^{2}-2|t| \|E_{1}\|_{\infty} \prod_{j=1}^{k} |\lambda - \tau_{j}|^{2},\  \text{(by inequality} \ \eqref{eq67}) \\ 
	&=& \prod_{j=1}^{k} |\lambda - \tau_{j}|^{2}  r \prod_{j=k+1}^{n} |Q_{\alpha_{j}}(\lambda)|-|t|^{2}\big| \prod_{j=1}^{k} |\lambda - \tau_{j}|^{2} \big|^{2}-2|t| \|E_{1}\|_{\infty} \prod_{j=1}^{k} |\lambda - \tau_{j}|^{2} \\ 
	&\geq& \prod_{j=1}^{k} |\lambda - \tau_{j}|^{2} \Bigg\{ r M-\bigg(|t|^{2}\prod_{j=1}^{k} |\lambda - \tau_{j}|^{2} +2|t| \|E_{1}\|_{\infty}\bigg) \Bigg\} \\ 
	&\geq& \prod_{j=1}^{k} |\lambda - \tau_{j}|^{2} \  \Bigg\{ r M-|t|\big(|t| \ \|g\|_{\infty} +2 \|E_{1}\|_{\infty}\big) \Bigg\}, 
	\end{eqnarray*}
where $M= \inf_{\lambda \in \mathbb{ T}} \prod_{j=k+1}^n |Q_{\alpha_{j}}(\lambda)| > 0$. 

Let us show that, for $|t|$ sufficiently small, $|D(\lambda)|^{2} - |E_{1}^{t}(\lambda)|^{2} \geq 0$ on $\mathbb{ T}$. It suffices to find $|t|$ such that 
$$r M-|t|\big(|t| \ \|g\|_{\infty} +2 \|E_{1}\|_{\infty}\big) > 0,$$ 
or equivalently, 
	\begin{equation*}
	|t| \bigg( |t| + 2 \frac{\|E_{1}\|_{\infty}}{\|g\|_{\infty}}\bigg)- \frac{rM}{\|g\|_{\infty}} <   0.
	\end{equation*}
If we take  $|t| \leq \min \bigg\{ \dfrac{2\|E_{1}\|_{\infty}}{\|g\|_{\infty}}, \dfrac{rM}{8\|E_{1}\|_{\infty}}  \bigg\}$, then 
	\begin{eqnarray} \label{eq68} \nonumber
	|t| \bigg( |t| + 2 \frac{\|E_{1}\|_{\infty}}{\|g\|_{\infty}}\bigg)- \frac{rM}{\|g\|_{\infty}} 
	&\leq& |t| \bigg(2 \frac{\|E_{1}\|_{\infty}}{\|g\|_{\infty}}  + 2 \frac{\|E_{1}\|_{\infty}}{\|g\|_{\infty}}\bigg)- \frac{rM}{\|g\|_{\infty}} \\ \nonumber
	 &\leq& \dfrac{rM}{8\|E_{1}\|_{\infty}} \bigg( 4 \frac{\|E_{1}\|_{\infty}}{\|g\|_{\infty}} \bigg) -\frac{rM}{\|g\|_{\infty}} \\  
	 &=& \frac{rM}{2\|g\|_{\infty}}-\frac{rM}{\|g\|_{\infty}} = -\frac{rM}{2\|g\|_{\infty}} < 0 .
	\end{eqnarray}
Therefore
\begin{equation*}
	|D|^{2}-|E_{1}^{t}|^{2} \geq 0 \qquad \ \text{on} \ \mathbb{T},
	\end{equation*}
	and, by Theorem \ref{cons}, the functions 
\begin{equation*}
	x_{ +t} = \bigg(  \frac{E_{1}^{ +t}}{D},\frac{E_{2}^{ +t}}{D}, \frac{D^{\sim n}}{D}  \bigg) \qquad \text{and} \qquad
	x_{ -t} = \bigg(  \frac{E_{1}^{ -t}}{D},\frac{E_{2}^{ -t}}{D}, \frac{D^{\sim n}}{D}  \bigg)
	\end{equation*}
are rational $\mathbb{ \overline{ E}}$-inner functions. Obviously,
\begin{eqnarray*} 
	\frac{1}{2}x_{+t} + \frac{1}{2}x_{-t}
	&=& \bigg(  \frac{E_{1}^{ +t}+E_{1}^{ -t}}{2D},\frac{{E_{2}^{+t}+ E_{2}^{-t}}}{2D}, \frac{D^{\sim n}}{D}  \bigg)  \\ 
		&=& \bigg(  \frac{E_{1} + tg + E_{1}-tg }{2D},\frac{E_{1}^{\sim n}+ tg + E_{1}^{\sim n}-tg }{2D}, \frac{D^{\sim n}}{D}  \bigg)  \\ 
		&=&  \bigg(  \frac{E_{1}}{D},\frac{E_{1}^{\sim n} }{D}, \frac{D^{\sim n}}{D}  \bigg) =x. 
	\end{eqnarray*}
Hence $x$ is not an extreme point of $\mathcal{J}$. \\
	
	If $n$ is odd, assume $n=2m+1$. This case requires a slight modification. By assumption, $2k \leq n$ thus $2k \leq 2m+ 1 $. This implies that $k \leq m$. Choose $\omega \in \mathbb{ T}$ such that 
	\begin{equation*}
	\omega^{2}= -\overline{\tau}_{1} \prod_{j=1}^{k} \overline{\tau}_{j}^{2}.
	\end{equation*}
	Let 
	\begin{equation*}
	g(\lambda)= \omega \lambda^{m-k}(\lambda - \tau_{1}) \prod_{j=1}^{k}(\lambda -\tau_{j})^{2}, \quad \lambda \in \mathbb{C}.
	\end{equation*}
	Clearly, the polynomial $g$ has degree $m+k+1 \leq n$. Let us check that the $g$ is $n$-symmetric
	\begin{eqnarray*}
	g^{\sim n}(\lambda)= \lambda^{n} 
\overline{g \big(1/\overline{\lambda} \big)} 
&=& \lambda^{n} \bigg(\overline{ \omega\frac{1}{\overline{\lambda}^{m-k}}\big(\frac{1}{\overline{\lambda}}-\tau_{1}\big) \prod_{j=1}^{k} \big( \frac{1}{\overline{\lambda}}-\tau_{j} \big)^{2}} \bigg) \\ 
	&=& \lambda^{n} \bigg( \overline{\omega}\frac{1}{\lambda^{m-k}}\big(\frac{1}{\lambda}-\overline{\tau}_{1}\big) \prod_{j=1}^{k} \big( \frac{1}{\lambda}-\overline{\tau}_{j} \big)^{2} \bigg) \\ 
	&=&  \overline{\omega}\lambda^{m-k} \overline{\tau}_{1}\big( \tau_{1} -\lambda \big) \prod_{j=1}^{k} \overline{\tau}_{j}^{2} \big( \tau_{j}-\lambda \big)^{2}  \\ 
	&=&  \overline{\omega} \ \bigg(- \overline{\tau}_{1} \prod_{j=1}^{k} \overline{\tau}_{j}^{2} \bigg) \ \lambda^{m-k} \big(  \lambda-\tau_{1} \big) \prod_{j=1}^{k}  \big(\lambda - \tau_{j} \big)^{2}  \\ 
	&=&  \overline{\omega} \ \omega^{2} \ \lambda^{m-k} \big(  \lambda-\tau_{1} \big) \prod_{j=1}^{k}  \big(\lambda - \tau_{j} \big)^{2}  \\ 
	&=&  \omega \ \lambda^{m-k} \big(  \lambda-\tau_{1} \big) \prod_{j=1}^{k}  \big(\lambda - \tau_{j} \big)^{2} =g(\lambda).
	\end{eqnarray*}
As in the even case, define the polynomials on $\mathbb{ \overline{D }}$
\begin{equation*}
	E_{1}^{t}= E_{1} + tg \qquad\text{and} \qquad  E_{2}^{t}= E_{1}^{\sim n} + tg \qquad  \text{for} \ t \in \mathbb{R}.
	\end{equation*}
As in equation \eqref{eq57}, for all $\lambda \in \mathbb{ \overline{D }}$, $E_{1}^{t}(\lambda)=\big( E_{2}^{t} \big)^{\sim n}(\lambda)$ and as in  equation \eqref{eqn66}, for all $\lambda \in \mathbb{ T}$,
	\begin{equation} \label{eq91}
	|D(\lambda)|^{2} - |E_{1}^{t}(\lambda)|^{2}= |D(\lambda)|^{2} -  |E_{1}(\lambda)|^{2} - t^{2}|g(\lambda)|^{2}-2 \text{Re}\big(tg(\lambda)\overline{ E_{1}(\lambda)}\big).
	\end{equation}
For all $\lambda \in \mathbb{ T}$, 
	\begin{eqnarray} \label{in10} \nonumber
	\text{Re}(tg\overline{E_{1}(\lambda)} \leq |tgE_{1}(\lambda)|&=& |t E_{1}(\lambda)| \bigg| \omega \lambda^{m-k}(\lambda - \tau_{1}) \prod_{j=1}^{k}(\lambda -\tau_{j})^{2}  \bigg| \\ \nonumber
	&\leq& |t| \ \|E_{1}\|_{\infty} |\lambda - \tau_{1}| \prod_{j=1}^{k}|\lambda -\tau_{j}|^{2}.  
	\end{eqnarray}
Combine equations \eqref{eq65}, \eqref{eq91} and inequality \eqref{in10} to obtain, for all $\lambda \in \mathbb{ T}$,
\begin{eqnarray*} 
	& & |D(\lambda)|^{2} - |E_{1}^{t}(\lambda)|^{2} \\ \nonumber
	 &=& |D(\lambda)|^{2} -  |E_{1}(\lambda)|^{2} - |t|^{2}|g(\lambda)|^{2}-2 \text{Re}\big(tg(\lambda)\overline{ E}_{1}(\lambda)\big)  \\ 
	&=& r \prod_{j=1}^{k} |\lambda - \tau_{j}|^{2} \prod_{j=k+1}^{n} |\lambda -\alpha_{j}|^{2}-|t|^{2}|g(\lambda)|^{2}-2\text{Re}\big(tg(\lambda)\overline{ E}_{1}(\lambda)\big), \ \text{by equation} \ \eqref{eq65} \\  
	&\geq& r \prod_{j=1}^{k} |\lambda - \tau_{j}|^{2} \prod_{j=k+1}^{n} |\lambda -\alpha_{j}|^{2}-|t|^{2}|g(\lambda)|^{2}-2|t| \|E_{1}\|_{\infty} |\lambda -\tau_{1}| \prod_{j=1}^{k} |\lambda - \tau_{j}|^{2}, \ \text{by inequality} \eqref{in10} \\ 
	&=& \prod_{j=1}^{k} |\lambda - \tau_{j}|^{2}  r \prod_{j=k+1}^{n} |Q_{\alpha_{j}}(\lambda)|-|t|^{2} |\lambda - \tau_{1}|^{2} \prod_{j=1}^{k}|\lambda -\tau_{j}|^{4}-2|t| \|E_{1}\|_{\infty} |\lambda -\tau_{1}|\prod_{j=1}^{k} |\lambda - \tau_{j}|^{2} \\ 
	&\geq& \prod_{j=1}^{k} |\lambda - \tau_{j}|^{2} \Bigg\{ r M- \underbrace{|\lambda -\tau_{1}|}_ {\leq 2}\bigg(\underbrace{|t|^{2}|\lambda - \tau_{1}|\prod_{j=1}^{k} |\lambda - \tau_{j}|^{2}}_{\leq |t|^{2} \|g\|_{\infty}} +2|t| \|E_{1}\|_{\infty}\bigg) \Bigg\} \\
	&\geq& \prod_{j=1}^{k} |\lambda - \tau_{j}|^{2} \Bigg\{ r M-2\big(|t|^{2}\|g\|_{\infty} +2|t| \|E_{1}\|_{\infty}\big) \Bigg\}, 
	\end{eqnarray*}
where $M= \inf_{\mathbb{ T}} \prod |Q_{\alpha_{j}}| > 0$. By similar arguments  to those in equations \eqref{eq68}, one can find $|t|$ such that 
\begin{equation*}
	r M-2\big(|t|^{2}\|g\|_{\infty} +2|t| \|E_{1}\|_{\infty}\big) > 0. 
	\end{equation*}
Therefore,
\begin{equation*}
	|D|^{2} - |E_{1}^{t}|^{2} \geq 0, \qquad \text{on} \ \mathbb{ T}.
	\end{equation*}
Hence, by Theorem \ref{cons}, the functions 
	\begin{equation*}
	x_{\pm t} = \bigg(  \frac{E_{1}^{ \pm t}}{D},\frac{\big(E_{1}^{\sim n}\big)^{ \pm t}}{D}, \frac{D^{\sim n}}{D}  \bigg)
	\end{equation*}
are rational $\mathbb{ \overline{ E}}$-inner functions. One can check that $x=\frac{1}{2}x_{+t}+ \frac{1}{2}x_{-t}$ and therefore $x$ is not an extreme point of $\mathcal{J}$. 
\end{proof}

\begin{theorem} {\em \cite[Theorem 5.13]{ALY15}} \label{ext}
	 A rational $\Gamma$-inner function $h \in \mathcal{R}_{\Gamma}^{n,k}$ is extreme in the set of rational $\Gamma$-inner functions if and only if $2k>n$.
\end{theorem}
\begin{proposition} \label{exetremE}
	 Let $x=(x_{1},x_{2}, x_{3}) \in \mathcal{R}^{n,k}$ be a rational $\mathbb{ \overline{ E}}$-inner function such that $x_{1}=x_{2}$ and $2k>n$. Then $x$ is an extreme point of the set  $\mathcal{J}$ of rational $\mathbb{ \overline{ E}}$-inner functions.
\end{proposition}
\begin{proof}
	By Lemma \ref{12} (1), the function $h=(s,p)=(2x_{1},x_{3})$ is $\Gamma$-inner.  By Theorem \ref{min}, there are polynomials $E_{1}, E_{2}, D$ such that $x=\big( \frac{E_{1}}{D}, \frac{E_{2}}{D}, \frac{D^{\sim n}}{D}\big)$. Here, since $x_{1}=x_{2}$, necessarily $ E_{1}= E_{2}$. 
 By Definition \ref{Royal}, the royal polynomial $R_{h}$ of $h$ is 
	\begin{eqnarray*}
	R_{h}(\lambda)&=& D^{2}(\lambda) \bigg(4x_{3} - 4x_{1}^{2} \bigg)(\lambda) \\
	&=& D^{2}(\lambda) \bigg( 4 \frac{D^{\sim n}}{D}- 4\frac{E_{1}^{2}}{D^{2}} \bigg)(\lambda) \\
	&=& 4\big(DD^{\sim n} - E_{1}^{2}\big)(\lambda) = 4R_{x} (\lambda).
 	\end{eqnarray*}
 It is clear that if $x \in \mathcal{R}^{n,k}$, then $h$ has degree $n$ and $k$ royal nodes on $\mathbb{ T}$, counted with multiplicities, such that $2k>n$. 
 Thus, by Theorem \ref{ext}, $h$ is an extreme point of the set of rational $\Gamma$-inner functions. That is, if $h^{1}=(s^{1},p^{1})$ and $h^{2}=(s^{2}, p^{2})$ are $\Gamma$-inner functions such that
	\begin{equation*}
	h=th^{1}+(1-t)h^{2} \qquad \mbox{for some} \ t \in (0,1),
	\end{equation*}
	then $h=h^{1}=h^{2}$. Note that, in this case, we have
	\begin{equation} \label{eqn751}
	\begin{cases}
		s=ts^{1}+(1-t)s^{2}  \qquad \Rightarrow s=s^{1}=s^{2} \\
		p=tp^{1}+(1-t)p^{2} \qquad \Rightarrow p=p^{1}=p^{2}.
	\end{cases}
	\end{equation}
	 Suppose 
	\begin{equation*}
	x=tx^{1} + (1-t)x^{2}, \quad \text{for some} \ t \in (0,1),
	\end{equation*}
	and for rational $\mathbb{ \overline{ E}}$-inner functions $x^{1}=(x_{1}^{1}, x_{2}^{1}, x_{3}^{1})$ and $x^{2}=(x_{1}^{2}, x_{2}^{2}, x_{3}^{2})$.
	This implies that
	\begin{equation*}
	\begin{cases}
	x_{1} = t x_{1}^{1} + (1-t)x_{1}^{2} \\
	 x_{1} = t x_{2}^{1} + (1-t)x_{2}^{2} \\
	  x_{3}=p=tx_{3}^{1}+(1-t)x_{3}^{2}.
	\end{cases}
	\end{equation*}
	Recall that $(s,p)=(2x_{1},x_{3})$, hence 
	\begin{equation} \label{eqn2223}
	\begin{cases}
	s= 2t x_{1}^{1} + 2(1-t)x_{1}^{2} \\
	 s= 2t x_{2}^{1} + 2(1-t)x_{2}^{2} \\
	 p=tx_{3}^{1}+(1-t)x_{3}^{2}.
	\end{cases}
	\end{equation}
	Therefore
	\begin{equation*}
	(s,p)=t\big(2x_{1}^{1},x_{3}^{1}\big)+ (1-t)\big(2x_{1}^{2}, x_{3}^{2}\big)
	\end{equation*}
	and
	\begin{equation*}
	(s,p)=t\big(2x_{2}^{1},x_{3}^{1}\big)+ (1-t)\big(2x_{2}^{2}, x_{3}^{2}\big).
	\end{equation*}
	Since $h$ is an extreme rational $\Gamma$-inner function, we have
	\begin{equation*}
	\begin{cases}
	2x_{1}^{1}= 2x_{1}^{2}=s \\
	 2x_{2}^{1}= 2x_{2}^{2}=s \\
	 x_{3}^{1}= x_{3}^{2}=p.
	\end{cases}
	\end{equation*}
	Therefore $x=x^{1}=x^{2}$. Hence $x$ is extreme in the set $\mathcal{J}$.
\end{proof}

\begin{remark}{\rm It is not clear that a rational $\mathbb{ \overline{ E}}$-inner function $x=(x_{1},x_{2}, x_{3}) \in \mathcal{R}^{n,k}$  such that $2k>n$ and $x_{1}\neq x_{2}$, is an extreme point of the set  $\mathcal{J}$ of rational $\mathbb{ \overline{ E}}$-inner functions.
Could we claim that in Lemma \ref{convex-nodes}, if  $\tau_i \in \T$ is a royal node of $x$ of multiplicity $\nu$, then $\tau_i$ is a royal node of $x^1$ and $x^2$ of the same multiplicity $\nu$? Here  $x =  t x^1 + (1-t) x^2$ for some $t$ such that $0 < t < 1$. If so then $x$ is an extreme point of  $\mathcal{J}$ if and only if $2k>n$.
}
\end{remark}

\bibliography{references}

OMAR  M. O. ALSALHI, Alleith University College, Umm-Alqura University, Al-Leith, Saudi Arabia; e-mail: omsalhi@uqu.edu.sa
\\

ZINAIDA A. LYKOVA,
School of Mathematics, Statistics and Physics, Newcastle University, Newcastle upon Tyne
 NE\textup{1} \textup{7}RU, U.K.;\\
 e-mail: Zinaida.Lykova@newcastle.ac.uk\\

\end{document}